\documentclass[11pt]{article}

\def\?{?\vadjust{\vbox to 0pt{\vss\hbox{\kern\hsize\kern1em\large\bf ?!}}}}

\usepackage[T2A]{fontenc}
\usepackage[english]{babel}
\usepackage[tbtags]{amsmath}
\usepackage{amsfonts,amssymb}
\usepackage{mathrsfs}
\usepackage{bm}
\usepackage{esint}

\newcounter{Th}[section] \newcounter{Lm}[section] \newcounter{Ca}[section] \newcounter{Prop}[section]
\newcounter{ThA}
\newcounter{Problem}[section] \newcounter{Remark}[section] \newcounter{Example}[section]
\newcounter{Def}[section]\newcounter{Assum}[section]

\def\theTh{\arabic{section}.\arabic{Th}}
\def\theThA{\Alph{ThA}}
\def\theLm{\arabic{section}.\arabic{Lm}}
\def\theCa{\arabic{section}.\arabic{Ca}}

\def\theRemark{\arabic{section}.\arabic{Remark}}
\def\theExample{\arabic{section}.\arabic{Example}}
\def\theDef{\arabic{section}.\arabic{Def}}

\newenvironment{Th}[1][\relax]
    {\medspace\refstepcounter{Th}{\bf Theorem \theTh.}\ \it}
    {\rm\medspace}

\newenvironment{Lm}[1][\relax]
    {\medspace\refstepcounter{Lm}{\bf Lemma \theLm.}\ \it}
    {\rm\medspace}

\newenvironment{Ca}[1][\relax]
    {\medspace\refstepcounter{Ca}{\bf Corollary  \theCa.}\ \it}
    {\rm\medspace}

\newenvironment{Remark}[1][\relax]
    {\medspace\refstepcounter{Remark}{\bf Remark \theRemark.}\rm\ }
    {\medspace}

\newenvironment{Def}[1][\relax]
    {\medspace\refstepcounter{Def}{\bf Definition \theDef.}\rm\ }
    {\medspace}

\makeatother

\numberwithin{equation}{section}

\date{}

\begin{document}
\author{\bfseries\large A.~I.~Tyulenev and S.~K.~Vodop'yanov%
\thanks{Steklov Mathematical Institute of the Russian Academy of Sciences, Moscow. E-mails:
tyulenev-math@yandex.ru, tyulenev@mi.ras.ru, vodopis@math.nsc.ru }}

\title{Sobolev $W_{p}^{1}$-spaces on
$d$-thick closed subsets of $\mathbb{R}^{n}$%
\thanks{This work is supported by the Russian Science Foundation under grant 14-50-00005.}}

\maketitle

Let $S \subset \mathbb{R}^{n}$ be a~closed nonempty set such that for some
$d \in [0,n]$ and $\varepsilon > 0$ the~$d$-Hausdorff content  $\mathcal{H}^{d}_{\infty}(S \cap Q(x,r)) \geq \varepsilon r^{d}$
for all cubes~$Q(x,r)$ centered in~$x \in S$ with side length $2r \in (0,2]$. For every $p \in (1,\infty)$, denote by $W_{p}^{1}(\mathbb{R}^{n})$
the classical Sobolev space on $\mathbb{R}^{n}$. We give an~intrinsic characterization of the restriction $W_{p}^{1}(\mathbb{R}^{n})|_{S}$ of the space  $W_{p}^{1}(\mathbb{R}^{n})$ to~the set $S$  provided that $p > \max\{1,n-d\}$. Furthermore, we prove the existence of a bounded linear operator $\operatorname{Ext}:W_{p}^{1}(\mathbb{R}^{n})|_{S} \to W_{p}^{1}(\mathbb{R}^{n})$ such that $\operatorname{Ext}$ is right inverse for the usual trace operator. In particular, for $p > n-1$ we characterize the trace space of the Sobolev space $W_{p}^{1}(\mathbb{R}^{n})$ to the closure $\overline{\Omega}$ of an arbitrary open path-connected set~$\Omega$.
Our results extend those available for $p \in (1,n]$ with much more stringent restrictions on~$S$.

\begin{flushleft}
\end{flushleft}
\def\B{\rlap{$\overline B$}B}
{\small\leftskip=10mm\rightskip=\leftskip\noindent}

\section{Introduction}

For $m \in  \mathbb{N}$, denote by $C^{m}(\mathbb{R}^{n})$ the linear space of all functions on
$\mathbb{R}^{n}$ with continuous partial derivatives up to order~%
$m$ with finite seminorm

$$
\|F|C^{m}(\mathbb{R}^{n})\|:=\sum\limits_{|\alpha|=m}\|D^{\alpha}F|C(\mathbb{R}^{n})\|.
$$

The classical extension problem posed by H. Whitney in 1934 in his famous papers \cite{W1}, \cite{W2}
reads as follows:

\textbf{Classical Whitney Extension Problem.} \textit{For $m \in \mathbb{N}$, how can we tell whether a given function
$f$ defined on an arbitrary closed subset $S \subset \mathbb{R}^{n}$
extends to a $C^{m}(\mathbb{R}^{n})$-function?}

Note that the problem mentioned above appeared to be very complicated. Whitney solved this problem only for the space $C^{m}(\mathbb{R})$, $m \geq 1$ and a similar problem for the Lipschitz space $C^{0,1}(\mathbb{R}^{n})$, $n \geq 1$. After the seminal papers \cite{W1}, \cite{W2} a big progress was made by many mathematicians \cite{Gl}, \cite{Br2}, \cite{Mil} (see also references therein). Only recently C. Fefferman gave complete solutions of the Classical Whitney Extension Problem and other closely related problems \cite{F1}--\cite{F6}.

These papers motivated the study of analogous problems for classical homogeneous $L_{p}^{m}(\mathbb{R}^{n})$ and inhomogeneous $W_{p}^{m}(\mathbb{R}^{n})$ Sobolev-type spaces \textit{in the case} $m \in \mathbb{N}$, $p > n$ \cite{F8}--\cite{F11}, \cite{Is}, \cite{Shv2}, \cite{Shv5}. Note that such problem was completely solved only in the cases $m=1, n \geq 1$, $p > n$ \cite{Shv2} and $m=n=2$, $p > 2$ \cite{Shv5}. Papers \cite{F8}--\cite{F11}, \cite{Is} dealt only with the problem of constructing of the bounded linear operator from the corresponding trace space.

The aim of this paper is to pose correctly and solve \textit{for every} $p \in (1,\infty)$ analog of the Classical Whitney Extension Problem in the context of the first-order Sobolev spaces $W^{1}_{p}(\mathbb{R}^{n})$  \textit{for a sufficiently large class of closed sets} $S \subset \mathbb{R}^{n}$.

\subsection{Main results}

Henceforth $D^{\alpha}f$, where $\alpha$ is a~multi-index, stands for the distributional (Sobolev) partial derivative of
a~function $f \in L^{\text{\rm loc}}_{1}(\mathbb{R}^{n})$. In what follows we set $D^{0}f:=f$.

Let $p \in [1,\infty]$ and $n \in \mathbb{N}$. For an~open set $G \subset \mathbb{R}^{n}$, denote by
$W_{p}^{1}(G)$ the Sobolev space of all equivalence classes of real-valued functions $F \in L_{p}(G)$
whose distributional partial derivatives on~%
$G$ belong to $L_{p}(G)$. Equip this space with the norm

\begin{equation}
\label{eq1.1}
\|F|W^{1}_{p}(G)\|:=\sum\limits_{|\alpha| \le 1}\|D^{\alpha}F|L_{p}(G)\|.
\end{equation}

Similarly, given an open set $G$ for every $p \in [1,\infty]$ and $n \in \mathbb{N}$ we can define the homogeneous Sobolev space $L^{1}_{p}(G)$ of all equivalence classes of real-valued functions $F \in L^{\operatorname{loc}}_{1}(G)$, with the seminorm

\begin{equation}
\label{eq1.2}
\|F|L^{1}_{p}(G)\|:=\sum\limits_{|\alpha| = 1}\|D^{\alpha}F|L_{p}(G)\|.
\end{equation}

By $B(x,r)$ ($Q(x,r)$) we  denote the closed ball (cube) centered at $x \in \mathbb{R}^{n}$ of radius $r > 0$ (with side length $2r$).
We say that $F \in W_{p}^{1,\text{loc}}(\mathbb{R}^{n})$ ($F \in L_{p}^{1,\text{loc}}(\mathbb{R}^{n})$) if and only if $F \in W_{p}^{1}(\operatorname{int}B(0,r))$ ($F \in L_{p}^{1}(\operatorname{int}B(0,r))$)
for every $r \in (0,\infty)$.

In order to formulate our main results we briefly recall (for the reader convenience) the basic notions of geometric measure theory.
For a~detailed exposition, see Ch.~5, Section~1 of~\cite{A} or Ch.~1, Ch.~2 of~\cite{Evans}. We present it here in a slightly different form (see Remark \ref{Rem2.1} for explanations).

Let $0 \le d \le n$, $S$ be a~subset of $\mathbb R^{n}$, and $\delta \in (0,+\infty]$.
Consider the set function
\begin{equation}
\notag
\mathcal{H}^{d}_{\delta} (S)=\inf \sum\limits_{j} r^{d}_{j}
\end{equation}
where the infimum is taken over all countable coverings of~$S$ by balls $B(x_{j},r_{j})$ with arbitrary centers
$x_{j}$ and radii $ r_{j} < \delta$. The~$d$-\textit{Hausdorff content} of a~set~$S$ is defined as $\mathcal{H}^{d}_{\infty} (S)$.
The~$d$-\textit{Hausdorff measure} of a~set~$S$ is defined as
\begin{equation}
\notag
\mathcal{H}^{d}(S) : = \lim\limits_{\delta \to 0}\mathcal{H}^{d}_{\delta} (S).
\end{equation}
The limit exists due to monotonicity of $\mathcal{H}^{d}_{\delta}(S)$ with respect to $\delta$.
It can be shown that for every $S \subset \mathbb{R}^{n}$ there exists a number $d_{0}(S) \in [0,n]$ such that

\begin{equation}
\notag
d_{0}(S)=\sup\{d: \mathcal{H}^{d}(S) = +\infty\} = \inf\{d':\mathcal{H}^{d'}(S)=0\}
\end{equation}
This number is called the \textit{Hausdorff dimension} of~$S$
and denoted by $\operatorname{dim}_{H}S$.

Let $d \in [0,n]$. We say that a set $S \subset \mathbb{R}^{n}$ is \textit{$d$-thick} if and only if there exists $\varepsilon > 0$ such that $\mathcal{H}^{d}_{\infty}(B(x,r) \cap S) \geq \varepsilon r^{d}$ for every $x \in S$ and every $r \in (0,1]$.

Fix a number $d \in [0,n]$ and $p \in (\max\{1,n-d\}, \infty]$. Recall (see Theorem \ref{Th2.4} below for details) that for every function $F \in W_{p}^{1}(\mathbb{R}^{n})$ there exist a set $E_{F} \subset \mathbb{R}^{n}$ and a representative (in the class of equivalent functions) $\widehat{F}$ such that $\mathcal{H}^{d}(E)=0$ and every point $x \in \mathbb{R}^{n} \setminus E_{F}$ is a Lebesgue point of the function $\widehat{F}$. Thus, in what follows we  identify each element $F \in
W_{p}^{1}(\mathbb{R}^{n})$ with an arbitrary chosen such representative. As a result, for every set $S \subset \mathbb{R}^{n}$ with $\operatorname{dim}_{H}S \geq d$ we can define \textit{the $d$-trace of a given element} $F \in W_{p}^{1}(\mathbb{R}^{n})$ to the set $S$ as the \textit{pointwise restriction of} $\widehat{F}$ to $S$.  In what follows, given $d$ and $p$ as above and  an element $F \in W_{p}^{1}(\mathbb{R}^{n})$, by the symbol $F|_{S}$ we will denote the $d$-trace of $F$ to $S$. Clearly, the $d$-trace $F|_{
S}$ is uniquely defined $\mathcal{H}^{d}$-a.e. on $S$. Hence, given a function $f:S \to \mathbb{R}$ we will write $F|_{S}=f$ if and only if $F|_{S}(x)=f(x)$ for $\mathcal{H}^{d}$-almost every $x \in S$.

Using the facts mentioned above, given $p$ and $d$ as above, we can consider  \textit{the $d$-trace space} of the space $W_{p}^{1}(\mathbb{R}^{n})|_{S}$. More precisely,
\begin{equation}
\notag
W_{p}^{1}(\mathbb{R}^{n})|_{S}:=\widetilde{W}_{p}^{1}(\mathbb{R}^{n})|_{S}/\mathring{W}_{p}^{1}(\mathbb{R}^{n})|_{S},
\end{equation}
where
\begin{equation}
\begin{split}
\notag
&\widetilde{W}_{p}^{1}(\mathbb{R}^{n})|_{S}:=\{f:S \to \mathbb{R} | \text{ there exists } F \in W_{p}^{1}(\mathbb{R}^{n}) \text{ such that } F|_{S}=f\},\\
&\mathring{W}^{1}_{p}(\mathbb{R}^{n})|_{S}:=\{f:S \to \mathbb{R} | \text{ there exists } F \in W_{p}^{1}(\mathbb{R}^{n}) \text{ such that } F|_{S}=0\}.
\end{split}
\end{equation}
Equip this space with the usual \textit{quotient-space norm.}

Define the \textit{trace operator} $\operatorname{Tr}|_{S}: W_{p}^{1}(\mathbb{R}^{n}) \to W_{p}^{1}(\mathbb{R}^{n})|_{S}$ which acts as follows
\begin{equation}
\notag
\operatorname{Tr}|_{S}[F]:=F|_{S}.
\end{equation}
Clearly, this operator is \textit{linear and bounded.}

Our main goal in this paper is a study of the following

\textbf{Problem~A}. \textit{Fix parameters $d \in [0,n]$, $p \in (\max\{1,n-d\},\infty]$ and a~closed set
$S \subset \mathbb{R}^{n}$ with $\operatorname{dim}_{H}S \geq d$.
Given a function $f: S \to \mathbb{R}$, how can we decide whether there exists a function $F \in W_{p}^{1}(\mathbb{R}^{n})$ such that the $d$-trace $F|_{S}=f$? Consider the  $W_{p}^{1}(\mathbb{R}^{n})$-norm of all functions $F \in W_{p}^{1}(\mathbb{R}^{n})$ such that $F|_{S}=f$ on $S$. How small can these norms be?}

Furthermore, in this article we consider closely related

\textbf{Problem~B}. \textit{Fix  parameters $d \in [0,n]$, $p \in (\max\{1,n-d\},\infty]$ and a~closed set
$S \subset \mathbb{R}^{n}$ with $\operatorname{dim}_{H}S \geq d$. Does there exist a bounded linear operator $\operatorname{Ext}:W_{p}^{1}(\mathbb{R}^{n})|_{S} \to W_{p}^{1}(\mathbb{R}^{n})$
such that $\operatorname{Tr}|_{S} \circ \operatorname{Ext} = \operatorname{Id}$ on $W_{p}^{1}(\mathbb{R}^{n})|_{S}$?}

\begin{Remark}
In the case where $p=\infty$, the Sobolev space $W_{\infty}^{1}(\mathbb{R}^{n})$
can be identified with the space $\operatorname{LIP}(\mathbb{R}^{n})$ of Lipschitz functions on $\mathbb{R}^{n}$ and it is known that the restriction $\operatorname{LIP}(\mathbb{R}^{n})|_{S}$ coincides with the space $\operatorname{LIP}(S)$ of Lipschitz functions on $S$ and that, furthermore, the
classical Whitney extension operator linearly and continuously maps the space $\operatorname{LIP}(S)$
into the space $\operatorname{LIP}(\mathbb{R}^{n})$ (see e.g., \cite{St2}, Chapter 6). Hence in the sequel \textit{we will deal only with the case} $1 < p <\infty$.
\end{Remark}

Let $d \in [0,n]$. Assume that $\operatorname{dim}_{H}S \geq d$. Let $\{\mu_{k}\}_{k \in \mathbb{N}_{0}}$ be a family of positive Borel measures with $\operatorname{supp}\mu_{k} \subset S$, $k \in \mathbb{N}_{0}$.
We say that $\{\mu_{k}\}_{k \in \mathbb{N}_{0}}$ is  a $d$-\textit{regular system of measures on} $S$ if and only if for some universal constants $C_{1}$, $C_{2}$, $C_{3}$ the following properties hold for every $k \in \mathbb{N}_{0}$:

{\rm (1)} $\mu_{k}(B(x,r)) \le C_{1} r^{d}$ for every $x \in \mathbb{R}^{n}$ and every $r \in (0,2^{-k}]$;

{\rm (2)} $\mu_{k}(B(x,2^{-k})) \geq C_{2} 2^{-k}$ for every $x \in S$;

{\rm (3)} $2^{d-n}\mu_{k}(G) \le \mu_{k-1}(G) \le \mu_{k}(G)$ for every Borel set $G \subset S$.

\begin{Remark}
\label{Rem1.1}
We will show in Corollary \ref{Ca3.1} below that for every closed $d$-thick set $S$ one can construct a $d$-regular system of measures on $S$.
\end{Remark}

Assume that $p \in [1,\infty]$. Let $\mathfrak{m}$ be an arbitrary Borel measure on $\mathbb{R}^{n}$. Given a Borel function $f$, we say that $f \in L^{\text{\rm loc}}_{p}(\mathbb{R}^{n},\mathfrak{m})$
if and only if $f \in L_{p}(B(x,r),\mathfrak{m})$ for all $x \in \mathbb{R}^{n}$ and $r > 0$.

Given a function $f \in L^{\text{\rm loc}}_{1}(\mathbb{R}^{n},\mathfrak{m})$, we set for every Borel set $G \subset \mathbb{R}^{n}$
\begin{equation}
\notag
\fint\limits_{G}f(x)\,d\mathfrak{m}(x):= \frac{1}{\mathfrak{m}(G)}\int\limits_{G}f(x)\,d\mathfrak{m}(x).
\end{equation}

Let $S$ be a closed set in $\mathbb{R}^{n}$ with $\operatorname{dim}_{H} S \geq d$ for some $d \in [0,n]$. Assume that there exists  a $d$-regular system of measures $\{\mu_{k}\}_{k \in \mathbb{N}_{0}}$ on $S$. Let $f \in L^{\text{loc}}_{1}(\mathbb{R}^{n},\mu_{k})$ for every $k \in \mathbb{N}_{0}$. Consider the following \textit{Calderon-type maximal function}. Given a number $t \in [0,1)$, for every $x \in \mathbb{R}^{n}$  we set
\begin{equation}
\notag
f^{\sharp}_{\{\mu_{k}\}}(x,t):=\sup\limits_{r \in (t,1)} \frac{1}{r}\fint\limits_{Q(x,r)}\Bigl|f(y)-\fint\limits_{Q(x,r)}f(z)\,d\mu_{k(r)}(z)\Bigr|\,d\mu_{k(r)}(y),
\end{equation}
where $k(r)$ \textit{is the unique integer number} for which $r \in [2^{-k(r)},2^{-k(r)+1})$.
In what follows we set $f^{\sharp}_{\{\mu_{k}\}}:=f^{\sharp}_{\{\mu_{k}\}}(\cdot,0)$ for brevity.

Fix a closed set $S$ and a parameter  $\lambda \in (0,1)$. For every $j \in \mathbb{N}_{0}$ define the \textit{maximal $2^{-j}$-porous subset} of $S$ as follows. We set for every $j \in \mathbb{N}_{0}$
\begin{equation}
\notag
S_{j}(\lambda):=\{x \in S | \hbox{ there  exists }  y \in Q(x,2^{-j}) \hbox{ for which }  Q(y,\lambda 2^{-j}) \subset \mathbb{R}^{n}\setminus S\}.
\end{equation}
If there exists a number $\lambda \in (0,1)$ such that $S_{j}(\lambda) = S$ for every $j \in \mathbb{N}_{0}$ we say that the set $S$ is \textit{porous}.

Now we are ready to formulate the main result of this paper which gives the solution of Problems \textbf{A} and $\textbf{B}$.


In what follows by $\mathcal{L}_{n}$ we denote the classical $n$-dimensional Lebesgue measure on $\mathbb{R}^{n}$. Recall Remark 1.1.

\begin{Th}
\label{Th1.1}
Let $d \in [0,n]$ and $p \in (\max\{1,n-d\},\infty)$. Let $S$ be a closed $d$-thick set in $\mathbb{R}^{n}$. Let $\{\mu_{k}\}_{k \in \mathbb{N}_{0}}$  be a $d$-regular system of measures on $S$. Then $f$ belongs to the $d$-trace space $W_{p}^{1}(\mathbb{R}^{n})|_{S}$ if and only if there exists a set $S' \subset S$ with $\mathcal{H}^{d}(S \setminus S') = 0$ such that 
\begin{equation}
\label{eq1.3}
\lim\limits_{r \to 0}\fint\limits_{Q(x,r) \cap S} |f(x)-f(z)|\,d\mu_{k(r)}(z) = 0, \quad \forall x \in S',
\end{equation}
and for some $\lambda \in (0,1)$
\begin{equation}
\begin{split}
\notag
&\mathcal{N}_{S,p,\lambda}[f]:=\left(\int\limits_{S}|f(x)|^{p}\,d\mu_{0}(x)\right)^{\frac{1}{p}}+\left(\int\limits_{S}\Bigl(f^{\sharp}_{\{\mu_{k}\}}(x)\Bigr)^{p}\,d\mathcal{L}_{n}(x)\right)^{\frac{1}{p}}\\
&+\left(\sum\limits_{k=0}^{\infty}\int\limits_{S_{k}(\lambda)}\Bigl(f^{\sharp}_{\{\mu_{k}\}}(x,2^{-k})\Bigr)^{p}\,d\mu_{k}(x)\right)^{\frac{1}{p}} < \infty.
\end{split}
\end{equation}
Furthermore
\begin{equation}
\label{eq1.4}
\|f|W_{p}^{1}(\mathbb{R}^{n})|_{S}\|  \sim \mathcal{N}_{S,p,\lambda}[f],
\end{equation}
and there exists a bounded linear operator $\operatorname{Ext}:W_{p}^{1}(\mathbb{R}^{n})|_{S} \to  W_{p}^{1}(\mathbb{R}^{n})$ such that $\operatorname{Tr}|_{S} \circ \operatorname{Ext} = \operatorname{Id}$ on
$W^{1}_{p}(\mathbb{R}^{n})|_{S}$.
\end{Th}

\begin{Remark}
\label{Rem1.2}
Assume that $S = \overline{\Omega}$ for some open path-connected set $\Omega$. It is not difficult to show (see Appendix below) the existence of a sufficiently small $\varepsilon > 0$ such that $\mathcal{H}^{1}_{\infty}(Q(x,r) \cap S) \geq \varepsilon r$ for every $x \in S$ and $r \in (0,1]$. In other words, $\overline{\Omega}$ is $1$-thick. Recall Remark \ref{Rem1.1}.
This  implies that Theorem \ref{Th1.1} provides for every $p > n-1$ a description of the trace space of the Sobolev space $W_{p}^{1}(\mathbb{R}^{n})$ to the closure of $\Omega$.
\end{Remark}

\subsection{Simplifications for sets with porous boundary}

The results of Theorem \ref{Th1.1} can be simplified if $S$ or $\mathbb{R}^{n} \setminus S$
possesses certain "plumpness" properties. More precisely, in this section we restrict ourselves to the case in which the set $S$ has the \textit{porous boundary} $\partial S$.

Let $d \in [0,n]$. Let $S$ be a nonempty closed subset of $\mathbb{R}^{n}$ with $\operatorname{dim}_{H}S \geq d$. Assume that there exists a $d$-regular system of measures $\{\mu_{k}\}_{k \in \mathbb{N}_{0}}$  on $S$. Assume that $f \in L_{1}^{\text{loc}}(S,\mu_{k})$ for every $k \in \mathbb{N}_{0}$.
Define for every $x \in S$ and $r \in (0,1]$ \textit{the normalized with respect to the measure $\mu_{k(r)}$ best approximation of $f$ by constants on $Q(x,r)$}

$$
\mathcal{E}_{\mu_{k(r)}}(f,Q(x,r)):=\inf\limits_{c \in \mathbb{R}}\fint\limits_{Q(x,r) \cap S}|f(y)-c|\,d\mu_{k(r)}(y).
$$

\begin{Remark}
\label{Rem1.3}
It is easy to see that
\begin{equation}
\notag
\mathcal{E}_{\mu_{k(r)}}(f,Q(x,r)) \le  \widetilde{\mathcal{E}}_{\mu_{k(r)}}(f,Q(x,r)) \le 2\mathcal{E}_{\mu_{k(r)}}(f,Q(x,r)),
\end{equation}
where
\begin{equation}
\notag
\widetilde{\mathcal{E}}_{\mu_{k(r)}}(f,Q(x,r)):=\fint\limits_{Q(x,r) \cap S}\Bigl|f(y)-\fint\limits_{Q(x,r) \cap S}f(z)\,d\mu_{k(r)}(z)\Bigr|\,d\mu_{k(r)}(y).
\end{equation}

The exact value of $\mathcal{E}_{\mu_{k(r)}}(f,Q(x,r))$ will not be important for us in the sequel. Hence, we can use
$\widetilde{\mathcal{E}}_{\mu_{k(r)}}(f,Q(x,r))$ instead of $\mathcal{E}_{\mu_{k(r)}}(f,Q(x,r))$ which is easier to compute.
\end{Remark}

Now we present criterion which is simpler to verify in practice. Namely, instead of Calderon-type maximal functions our simplified criterion uses the normalized best approximations.

Given a closed nonempty set $S$, for every $k \in \mathbb{N}_{0}$ consider the set
\begin{equation}
\notag
\Sigma_{k}:=\Sigma_{k}(S):=\{x \in S| \operatorname{dist}(x,\partial S) \le 2^{-k}\}.
\end{equation}

Recall Remark 1.1 and definition of porous sets from the previous subsection.

\begin{Th}
\label{Th1.2}
Let $d \in [0,n]$ and $p \in (\max\{1,n-d\},\infty)$. Let $S$ be a closed $d$-thick set in $\mathbb{R}^{n}$. Let $\{\mu_{k}\}_{k \in \mathbb{N}_{0}}$ be a $d$-regular system of measures on $S$. Assume that $\partial S$ is porous. Then $f$ belongs to the $d$-trace space $W_{p}^{1}(\mathbb{R}^{n})|_{S}$ if and only if there exists a set $S' \subset S$ with $\mathcal{H}^{d}(S \setminus S') = 0$ such that 
\begin{equation}
\notag
\lim\limits_{r \to 0}\fint\limits_{Q(x,r) \cap S} |f(x)-f(z)|\,d\mu_{k(r)}(z) = 0, \quad \forall x \in S',
\end{equation}
and
\begin{equation}
\notag
\begin{split}
&\mathcal{R}_{S,p}[f]:=\left(\int\limits_{S}|f(x)|^{p}\,d\mu_{0}(x)\right)^{\frac{1}{p}}+\left(\int\limits_{S}\bigl(f^{\sharp}_{\{\mu_{k}\}}(x)\bigr)^{p}\,d\mathcal{L}_{n}(x)\right)^{\frac{1}{p}}\\
&+\left(\sum\limits_{k=0}^{\infty}2^{kp(1-\frac{n-d}{p})}\int\limits_{\Sigma_{k}}\Bigl(\mathcal{E}_{\mu_{k}}(f,Q(x,2^{-k}))\Bigr)^{p}\,d\mu_{k}(x)\right)^{\frac{1}{p}} <
\infty.
\end{split}
\end{equation}
Furthermore
\begin{equation}
\label{eq1.5}
\|f|W_{p}^{1}(\mathbb{R}^{n})|_{S}\|  \sim \mathcal{R}_{S,p}[f],
\end{equation}
and there exists a bounded linear operator $\operatorname{Ext}:W_{p}^{1}(\mathbb{R}^{n})|_{S} \to  W_{p}^{1}(\mathbb{R}^{n})$ such that $\operatorname{Tr}|_{S} \circ \operatorname{Ext} = \operatorname{Id}$ on
$W^{1}_{p}(\mathbb{R}^{n})|_{S}$.
\end{Th}

\begin{Remark}
\label{Rem1.4}
In the case in which the set $S$ is \textit{Ahlfors} $d$-\textit{regular} our result coincides with that obtained in \cite{Shv1} ($d=n$) and \cite{Ihn} ($d \in [0,n)$). We will present the details in Section 4.
\end{Remark}

\subsection{Brief overview of previously known results}

Analogs of Problems~\textbf{A} and~\textbf{B} can be posed in the context of Sobolev spaces $W_{p}^{m}(\mathbb{R}^{n})$, $m \in \mathbb{N}$
and even more complicated Besov and Triebel--Lizorkin spaces. Many articles address these problems.
Let us just mention \cite{F8}--\cite{F11}, \cite{Gal},  \cite{Ihn}, \cite{Is}, \cite{Jon1}--\cite{Jon4}, \cite{J}, \cite{Ka}, \cite{Ry1},
\cite{Shv1}--\cite{Shv5} (see also the references there). The reader can also consult the books
\cite{Maz1} and \cite{Maz2}, where many results are collected related to the trace problems for Sobolev spaces
and their applications.

Avoiding a~detailed review of all available results, we consider only recent breakthroughs.
Note that in the context of the first-order Sobolev spaces $W_{p}^{1}(\mathbb{R}^{n})$ ($L_{p}^{1}(\mathbb{R}^{n})$) Problems~\textbf{A} and~\textbf{B}
have been solved either under the condition $p > n$ without restrictions on the closed set~$S$ \cite{Shv2} or for all $p \in (1,\infty]$ under extra regularity assumptions on~%
$S$ \cite{Shv1}, \cite{Jon1}, \cite{Jon2}, \cite{Ihn}. In~particular, all articles cited above avoid the case in which $p \in (n-1,n]$, and~%
$S$ is the closure of an~arbitrary open path-connected subset of $\mathbb{R}^{n}$ (compare with Remark \ref{Rem1.2}).

We would like to note that Rychkov \cite{Ry1} introduced $d$-thick sets and solved analog of the Problem \textbf{B} for Besov spaces $B^{s}_{p,q}(\mathbb{R}^{n})$ and Triebel-Lizorkin spaces $F^{s}_{p,q}(\mathbb{R}^{n})$ under some restrictions on the parameters $s,p,q,d$. Recall that $W_{p}^{m}(\mathbb{R}^{n})=F^{m}_{p,2}(\mathbb{R}^{n})$ for $p \in (1,\infty)$, $m \in \mathbb{N}$. From this fact it follows that results obtained in \cite{Ry1} yield solution to the Problem \textbf{B} for the space $W_{p}^{1}(\mathbb{R}^{n})$, $p \in (1,\infty)$ only in the case $d > n-1$. Clearly, our results partially overlap with \cite{Ry1}. Nevertheless, in the important for applications case in which the set $S$ is the closure of an arbitrary open path-connected set $\Omega \subset \mathbb{R}^{n}$ Problems \textbf{A} and \textbf{B} were not solved in the paper \cite{Ry1} (compare with Remark \ref{Rem1.2})

\subsection{Plan of the paper}

Let us briefly describe the structure of this article.

Section~2 contains the standard definitions, notations, and classical lemmas often used below. 
Moreover include in this section some helpful results useful in what follows. Perhaps, such simple results are not knew but we can not provide the reader with a good reference.

In Section~3 we establish important properties of our key tools: \textit{$d$-regular system of measures}, \textit{Calderon-type maximal functions}, and \textit{porous sets}.

In Section~4 we obtain solutions to  Problems~\textbf{A} and~\textbf{B}
in the case of~$d$-thick closed sets. Namely, we present the proof of Theorem \ref{Th1.1}.

In Section 5 we prove Theorem \ref{Th1.2} which is a refined and simplified version of Theorem \ref{Th1.1} for sets with porous boundary.

In Section 6 we will consider several useful examples illustrating our main Theorems \ref{Th1.1} and \ref{Th1.2}. In particular we show that our results coincide with previously known results in the case of  closed Ahlfors regular sets $S$. Furthermore, we present an explicit construction of a $d$-regular system of measures on an arbitrary sharp closed single cusp. This leads us to a simplified version of Theorem 1.2 for the case of a closed single cusp. We would like to note that even this particular case of Theorem 1.2 is knew
and was never considered in the literature in such generality.

Finally, we decided to include Appendix to make this paper self-contained. In Appendix we present detailed explanations of some simple auxiliary examples upon which the reader will come during reading the paper. For instance, we prove that every open path-connected set is $1$-thick. Furthermore we included the proof of one technical result, which is difficult to find in the literature. We also included refined version of the Frostman-type Lemma which is suitable for our purposes.

\section{Preliminaries}
Our purpose in this section is to collect the required auxiliary material, fix definitions and notation. Furthermore we include some helpful lemmas 
which are useful in what follows. Such lemmas looks like standard but it is difficult to provide the reader with exact references. Hence we include complete proofs
of these statements.

Throughout the paper we use standard notation. By $x=(x_{1},...,x_{n})$ we denote an element of the space $\mathbb{R}^{n}$. Symbols $\alpha$, $\beta$ will be used to denote
multi-indices, i.e. elements of the space $\mathbb{N}^{n}_{0}$.

Throughout the paper, $B(x,r)$ stands for the closed ball (in the standard Euclidean metric) centered at~$x$ of radius $r > 0$. By $Q(x,r)$ we denote the closed cube centered at $x$ with side length $2r \geq 0$ with the edges parallel to the coordinate axes, namely, $Q(x,r):=\prod\limits_{i=1}^{n}[x_{i}-r,x_{i}+r]$.

Let $B=B(x,r)$ ($Q=Q(x,r)$) be a ball (a cube) in $\mathbb{R}^{n}$. Given a number $c > 0$, we write $cB$ ($cQ$) to denote the ball $B(x,cr)$ (the cube $Q(x,cr)$).

In what follows by dyadic cube we mean an arbitrary (half-open) cube $\widetilde{Q}_{k,m}:=\prod\limits_{i=1}^{n}\left[\frac{m_{i}}{2^{k}}, \frac{m_{i}+1}{2^{k}}\right)$, $k \in \mathbb{Z}$,
$m=(m_{1},...,m_{n}) \in \mathbb{Z}^{n}$. Given $k \in \mathbb{Z}$, by $\mathcal{Q}_{k}$ we will denote the mesh of all dyadic cubes with side length $2^{-k}$.

For $E \subset \mathbb{R}^{n}$, denote by $\overline{E}$ and $\operatorname{int}E$ the closure and interior of~%
$E$ in the topology induced by the standard Euclidean norm on~$\mathbb{R}^{n}$ respectively. Recall that \textit{all norms on $\mathbb{R}^{n}$ are equivalent}. Hence, topologies induced by every such norm are coincide with the topology induced by the standard Euclidean norm.

For $A \subset \mathbb{R}^{n}$ and $\delta > 0$, define the~$\delta$-neighbourhood
of~$A$ as $U_{\delta}(A):=\bigcup\limits_{x \in A}\operatorname{int}B_{\delta}(x)$.

Following \cite{Shv1}, it will be convenient for us to measure distances in $\mathbb{R}^{n}$ in the uniform norm
$$
\|x\|=\|x\|_{\infty}:=\max\{|x_{i}|: i=1,..,n\}, \quad x=(x_{1},...,x_{n}) \in \mathbb{R}^{n}.
$$
Given two subsets $A,B \subset \mathbb{R}^{n}$, put
$$
\operatorname{diam}A:=\sup\{\|a-a'\|_{\infty}: a,a' \in A\}, \quad
\operatorname{dist}(A,B):=\inf\{\|a-b\|_{\infty}: a \in A, b \in B\}.
$$

In what follows  the symbol $\mathcal{B}(\mathbb{R}^{n})$  denotes the $\sigma$-\textit{algebra of all Borel subsets of} $\mathbb{R}^{n}$. By \textit{Borel measure on} $\mathbb{R}^{n}$ we mean an arbitrary $\sigma$-additive $\sigma$-finite function $\mathfrak{m}: \mathcal{B}(\mathbb{R}^{n}) \to (0,+\infty]$.

Given a Borel measure $\mathfrak{m}$ and a nonempty set $S \subset \mathbb{R}^{n}$, we define \textit{restriction of} $\mathfrak{m}$ to $S$. More precisely, we set
$\mathfrak{m}\lfloor S (G) := \mathfrak{m}(G \cap S)$ for every Borel set $G$.

In what follows the classical~$n$-dimensional Lebesgue measure of a~Lebesgue measurable set~$A \subset \mathbb{R}^{n}$
will be denoted by~$\mathcal{L}_{n}(A)$.

Assume that $p \in [1,\infty]$. Let $\mathfrak{m}$ be an arbitrary Borel measure on $\mathbb{R}^{n}$. Given a Borel function $f:\mathbb{R}^{n} \to \mathbb{R}$, we say that $f \in L^{\text{\rm loc}}_{p}(\mathbb{R}^{n},\mathfrak{m})$
if and only if $f \in L_{p}(B(0,r),\mathfrak{m})$ for all  $r > 0$.

Given $f \in L^{\text{\rm loc}}_{1}(\mathbb{R}^{n},\mathfrak{m})$ we set for every Borel set $G$ with $\mathfrak{m}(G) < +\infty$

$$
\fint\limits_{G}f(x)\,d\mathfrak{m}(x):= \frac{1}{\mathfrak{m}(G)}\int\limits_{G}f(x)\,d\mathfrak{m}(x).
$$

\subsection{Geometric measure theory}

In this section we briefly recall basic facts from geometric measure theory.

Let $d \in [0,\infty)$, $\delta \in (0,\infty]$ and $S \subset \mathbb{R}^{n}$. Define

\begin{equation}
\label{eq2.1}
\mathcal{H}^{d}_{\delta}(S):=\inf\Bigl\{\sum\limits_{j \in \mathbb{N}}\alpha(d)\Bigl(\frac{\operatorname{diam}G_{j}}{2}\Bigr)^{d}| S \subset \bigcup\limits_{j \in \mathbb{N}}G_{j}, \operatorname{diam}G_{j}
\le \delta\Bigr\},
\end{equation}

where

$$
\alpha(d)=\frac{\pi^{\frac{d}{2}}}{\Gamma(\frac{d}{2}+1)}.
$$

Here $\Gamma(d)=\int\limits_{0}^{\infty}x^{d-1}e^{-x}\,dx$, ($0 < d < \infty$), is the usual gamma function.

\begin{Remark}
\label{Rem2.1}
Using the definition of $\mathcal{H}^{d}_{\delta}(S)$ and the fact that every bounded set $G$ is contained in a ball of diameter $2\operatorname{diam}G$ we see that
$$
\mathcal{H}^{d}_{\delta}(S) \le \widetilde{\mathcal{H}}^{d}_{\delta}(S) \le 2^{d}\mathcal{H}^{d}_{\delta}(S),
$$

where

$$
\widetilde{\mathcal{H}}^{d}_{\delta}(S):=\inf\Bigl\{\sum\limits_{j \in \mathbb{N}}\alpha(d)r_{j}^{d}| S \subset \bigcup\limits_{j \in \mathbb{N}}B(x_{j},r_{j}), 2r_{j} \le \delta\Bigr\}.
$$

Exact value of $\mathcal{H}^{d}_{\delta}(S)$ will not play any essential role in the sequel.
Hence, we will not distinguish between  $\widetilde{\mathcal{H}}^{d}_{\delta}(S)$ and $\mathcal{H}^{d}_{\delta}(S)$. This fact was used in the first section of this article.
\end{Remark}

It is easy to see that given a set $S \subset \mathbb{R}^{n}$ the function $\mathcal{H}^{d}_{\delta}(S)$ \textit{decreases} when $\delta$ \textit{increases.} This fact allows to introduce the following

\begin{Def}
\label{Def2.1}
For $S$ and $d$ as above we call $\mathcal{H}^{d}_{\infty}(S)$  \textit{$d$-Hausdorf content} of the set $S$. We also define
\begin{equation}
\label{eq2.2}
\mathcal{H}^{d}(S):=\lim\limits_{\delta \to 0}\mathcal{H}^{d}_{\delta}(S)=\sup\limits_{\delta > 0}\mathcal{H}^{d}_{\delta}(S)
\end{equation}

and call $\mathcal{H}^{d}(S)$  \textit{$d$-Hausdorf measure} of the set $S$.
\end{Def}

\begin{Def}
\label{Def2.2}
Let $d \in [0,n]$. We say that a~set~$S$ \textit{is Ahlfors~$d$-regular}
if there exist constants $C_{1}, C_{2} > 0$ such that
$$
C_{2}r^{d} \le \mathcal{H}^{d}(Q(x,r) \cap S) \le C_{1}r^{d}
$$
for every cube $Q=Q(x,r)$ with $x \in S$ and $r \in (0,1]$.
\end{Def}

The following definition is a natural generalization (see \textbf{Example 2.1} for explanations) of the previous and is taken from~\cite{Ry1}.

\begin{Def}
\label{Def2.3}
Let $d \in [0,n]$. A~set $S \subset \mathbb{R}^{n}$ is called $d$-\textit{thick} if there exist constants $C'_{1}, C'_{2} > 0$ such that such that
$$
C'_{2}r^{d} \le \mathcal{H}^{d}_{\infty}(Q(x,r) \cap S)  \le C'_{1}r^{d}
$$
for all $x \in S$ and $r \in (0,1]$.
\end{Def}

\begin{Remark}
\label{Rem2.2}
It is not difficult to show (see Section~5 of \cite{A}) that for each $d \in [0,n]$
the conditions $\mathcal{H}^{d}(S) = 0$ and $\mathcal{H}^{d}_{\infty}(S) = 0$ are equivalent.
This obviously implies that every~$d$-thick set is of Hausdorff dimension $\operatorname{dim}_{H}S \geq d$.
\end{Remark}

\textbf{Example 2.1.} In order to illustrate the notion of $d$-thick set we present below several useful examples. Details of these examples see in Appendix.

{\rm (1)} Let $d \in [0,n]$. Every Ahlfors~$d$-regular set $S$ is $d$-thick. The converse is false. Hence, the class of Ahlfors $d$-regular sets is strictly contained in the class of $d$-thick sets.

{\rm (2)} Let $\Omega$ be an arbitrary open path-connected subset of $\mathbb{R}^{n}$. Then $\Omega$ and $\overline{\Omega}$ are $1$-thick.

{\rm (3)} Let $\varepsilon > 0$, $\delta > 0$. The class of all path-connected $(\varepsilon,\delta)$-domains (\cite{J}) is strictly contained in the class of all Ahlfors $n$-regular sets. Hence, every path-connected $(\varepsilon,\delta)$-domain is $n$-thick.

{\rm (4)} Let $\varepsilon > 0$, $\delta > 0$. Let $\Omega$ be an arbitrary $(\varepsilon,\delta)$-domain in $\mathbb{R}^{n-1}$. Let $\varphi_{1}$, $\varphi_{2}$ be continuous functions on $\overline{\Omega}$ such that $\varphi_{2} < 0$ on $\Omega$, $\varphi_{1} > 0$ on $\Omega$  and $\varphi_{1}=\varphi_{2} = 0$ on $\partial \Omega$. Consider the set $G:=\{x=(x',x_{n})| x' \in \overline{\Omega}, x_{n} \in [\varphi_{2}(x'),\varphi_{1}(x')]\}$. Then the set $G$ is $(n-1)$-thick.

\subsection{Fine properties of functions}

In this section we recall several well-known facts about "good" pointwise behavior of functions from Lebesgue and Sobolev spaces.

Given a function $F \in L^{\text{\rm loc}}_{1}(\mathbb{R}^{n})$ and parameters $0 < t < s \le \infty$, put
$$
\operatorname{M}^{< s}_{> t}[F](x):=\sup_{r \in (t,s)}\fint\limits_{B(x,r)}|F(y)|\,d\mathcal{L}_{n}(y), \quad x \in \mathbb{R}^{n}.
$$

In what follows we will use notation $\operatorname{M}$ instead of $\operatorname{M}^{< \infty}_{> 0}$.

\begin{Remark}
\label{Rem2.1'}
Assume that $0 < t' \le t < s \le s'  \le +\infty$. It is clear that
\begin{equation}
\label{eq2.3'}
\operatorname{M}^{< s}_{> t}[F](x) \le \operatorname{M}^{< s'}_{> t'}[F](x), \quad x \in \mathbb{R}^{n}.
\end{equation}

It is obvious that $Q(y,cr) \supset Q(x,r)$ for every $y \in Q(x,t)$ and every $c\geq 2$, $r \geq t$. Hence

\begin{equation}
\notag
\fint\limits_{Q(x,r)}|F(y)|\,d\mathcal{L}_{n}(y) \le C(n,c) \fint\limits_{Q(y,cr)}|F(y)|\,d\mathcal{L}_{n}(y).
\end{equation}

This formula together with \eqref{eq2.3'} gives quite useful estimate

\begin{equation}
\label{eq2.4'}
\operatorname{M}_{> t}[F](x) \le C(n,c) \operatorname{M}_{> ct}[F](y) \le C(n,c) \operatorname{M}_{> t}[F](y), \quad y \in Q(x,t).
\end{equation}
\end{Remark}

The following result is classical. One can find its proof, for example in \cite{St2}, Ch.~1.

\begin{Th}
\label{Th2.1}
Let $p \in (1,\infty)$. Then~%
$\operatorname{M}$ is a~bounded operator from $L_{p}(\mathbb{R}^{n})$ into $L_{p}(\mathbb{R}^{n})$.
\end{Th}

Now we would like to formulate result, which is perhaps not knew. Nevertheless, we can not provide the reader with a good reference. The proof of this result is included into our Appendix.

Recall that a measure $\mathfrak{m}$ on $\mathbb{R}^{n}$ is a Radon measure if $\mathfrak{m}$ is Borel regular and $\mathfrak{m}(K) < \infty$ for each compact set $K \subset \mathbb{R}^{n}$.

\begin{Th}
\label{Th2.2}
Let $d \in [0,n]$  and $\gamma \in (1,\infty)$. Let $\mathfrak{m}$ be a Radon measure on $\mathbb{R}^{n}$ such that for some (universal) constant $C > 0$
\begin{equation}
\label{eq2.3}
\mathfrak{m}(B(x,r)) \le C r^{d}, \quad x \in \mathbb{R}^{n}, r > 0.
\end{equation}

Given $\alpha \geq 0$ and $0 < s < +\infty$, consider the following maximal function
$$
\operatorname{M}^{< s}[f,\alpha](x):=\sup\limits_{r \in (0,s)}r^{\alpha}\fint\limits_{B(x,r)}|f(z)|\,d\mathcal{L}_{n}(z).
$$

If $d \geq n-\alpha$, then the operator $\operatorname{M}^{< s}[\cdot,\alpha]$ is bounded from $L_{\gamma}(\mathbb{R}^{n})$ into $L_{\gamma}(\mathbb{R}^{n},\mathfrak{m})$.
\end{Th}

Now we formulate simple but useful for us result. One can find the proof in section 2.4.3 of \cite{Evans}.

\begin{Th}
\label{Th2.3}
Suppose that $d \in [0,n)$. Then, given a~function $F \in L_{1}^{\operatorname{loc}}(\mathbb{R}^{n})$,
there exists a~set $E_{F} \subset \mathbb{R}^{n}$ with $\mathcal{H}^{d}(E_{F})=0$ such that
$$
\lim\limits_{r \to 0}\frac{1}{r^{d}}\int\limits_{Q(x,r)}|F(y)|\,dy = 0
$$

for every $x \in \mathbb{R}^{n} \setminus E_{F}$.
\end{Th}

The following result helps us to define the trace of a Sobolev function $F$ to a set $S$ with "sufficiently big" Hausdorff dimension.

\begin{Th}
\label{Th2.4}
Suppose that  $q \in (1,\infty)$, $d \in [0,n]$, $d > n-q$.  Assume that $F \in W_{q}^{1,\operatorname{loc}}(\mathbb{R}^{n})$.
Then there exists a~set $E_{F} \subset \mathbb{R}^{n}$ with
$\mathcal{H}^{d}(E_{F}) = 0$ and a representative $\widehat{F}$ of the element $F$ such that every point $x \in \mathbb{R}^{n} \setminus E_{F}$ is a Lebesgue point of the function $\widehat{F}$.
\end{Th}

\begin{Def}
\label{Def2.4}
Let $p \in (1,\infty)$, $d \in [0,n]$, $d > n-p$. Let $S$ be a subset of $\mathbb{R}^{n}$ with $\operatorname{dim}_{H}S \geq d$. Given an element $F \in W_{p}^{1}(\mathbb{R}^{n})$, consider a representative $\widehat{F}$ of the element $F$ which has Lebesgue points $\mathcal{H}^{d}$-almost everywhere on $\mathbb{R}^{n}$. By \textit{a $d$-trace} $F|_{S}$
\textit{of the element} $F$ \textit{to the set} $S$ we mean the \textit{pointwise restriction} of the representative $\widehat{F}$ to the set $S$.
\end{Def}

\begin{Remark}
\label{Rem2.3}
Note the $d$-trace $F|_{S}$ of a given element $F \in W_{p}^{1}(\mathbb{R}^{n})$ is uniquely defined up to the set of $\mathcal{H}^{d}$-measure zero. Indeed, given an element $F$, let  $\widehat{F}_{1}$, $\widehat{F}_{2}$ be representatives of $F$ which have Lebesgue points $\mathcal{H}^{d}$-almost everywhere in $\mathbb{R}^{n}$. Then, pointwise restrictions of the elements $\widehat{F}_{1}$, $\widehat{F}_{2}$ are coincide $\mathcal{H}^{d}$-almost everywhere on $S$. As a result, strictly speaking, the $d$-trace $F|_{S}$ is a class of equivalent functions defined on $S$ modulo coincidence $\mathcal{H}^{d}$-almost everywhere on $S$.
Hence, in what follows, given a function $f:S \to \mathbb{R}$ we will write $F|_{S}=f$ if and only if there is a representative $\widehat{F}$ of the element $F$ such that $f(x)=\widehat{F}(x)$ \textit{for $\mathcal{H}^{d}$-almost every} $x \in S$.
\end{Remark}

Using Definition \ref{Def2.4} and Remark \ref{Rem2.3} we introduce

\begin{Def}
\label{Def2.5}
Let $p \in (1,\infty)$, $d \in [0,n]$, $d > n-p$. Let $S$ be a subset of $\mathbb{R}^{n}$ with $\operatorname{dim}_{H}S \geq d$. Define the \textit{$d$-trace space} of the space $W_{p}^{1}(\mathbb{R}^{n})$ as follows:

Consider the linear spaces
\begin{equation}
\begin{split}
\notag
&\widetilde{W}_{p}^{1}(\mathbb{R}^{n})|_{S}:=\{f:S \to \mathbb{R} | \text{ there exists } F \in W_{p}^{1}(\mathbb{R}^{n}) \text{ such that } F|_{S}=f\},\\
&\mathring{W}^{1}_{p}(\mathbb{R}^{n})|_{S}:=\{f:S \to \mathbb{R} | \text{ there exists } F \in W_{p}^{1}(\mathbb{R}^{n}) \text{ such that } F|_{S}=0\}.
\end{split}
\end{equation}
Call the quotient-space $\widetilde{W}_{p}^{1}(\mathbb{R}^{n})|_{S}/ \mathring{W}^{1}_{p}(\mathbb{R}^{n})|_{S}$ trace space of the space $W_{p}^{1}(\mathbb{R})$ and denote it by $W_{p}^{1}(\mathbb{R}^{n})|_{S}$. Equip the space $W_{p}^{1}(\mathbb{R}^{n})|_{S}$ with the standard quotient-space norm
\begin{equation}
\notag
\|f|W_{p}^{1}(\mathbb{R}^{n})|_{S}\|:=\inf \|F|W_{p}^{1}(\mathbb{R}^{n})\|,
\end{equation}
where the infimum is taken over all elements $F \in W_{p}^{1}(\mathbb{R}^{n})$ such that $F|_{S}=f$.

Furthermore, we define the \textit{trace operator} $\operatorname{Tr}|_{S}: W_{p}^{1}(\mathbb{R}^{n}) \to W_{p}^{1}(\mathbb{R}^{n})|_{S}$ which acts as follows
\begin{equation}
\notag
\operatorname{Tr}|_{S}[F]:=F|_{S} \quad \text{for every } F \in W_{p}^{1}(\mathbb{R}^{n}).
\end{equation}
\end{Def}

Formally speaking we have to add symbol $d$ to our notation $W_{p}^{1}(\mathbb{R}^{n})|_{S}$ of the $d$-trace space. However we will not do it to simplify notation.
Furthermore it will be always clear from the context which $d$ is assumed.

\begin{Remark}
\label{Rem2.4}
Our definition clearly implies that the trace operator $\operatorname{Tr}|_{S}: W_{p}^{1}(\mathbb{R}^{n}) \to W_{p}^{1}(\mathbb{R}^{n})|_{S}$ is \textit{linear and bounded.}
\end{Remark}

\subsection{Whitney decomposition}

The following result is \textit{the Classical Whitney Decomposition Lemma}. Recall that we measure distances in $\mathbb{R}^{n}$ in the uniform norm.

\begin{Lm}
\label{Lm2.1}
For every closed set $S \subset \mathbb{R}^{n}$ there exists a~family of closed dyadic cubes
$W_{S}=\{Q_{\varkappa}\}_{\varkappa \in I} = \{Q(x_{\varkappa},r_{\varkappa})\}_{\varkappa \in I}$ such that

{\rm (1)}
$\mathbb{R}^{n} \setminus S= \bigcup\limits_{\varkappa \in I} Q_{\varkappa}$;

{\rm (2)}
for each $\varkappa \in I$
\begin{equation}
\label{eq2.4}
\operatorname{diam}(Q_{\varkappa}) \le \operatorname{dist}(Q_{\varkappa}, S) \le  4\operatorname{diam}(Q_{\varkappa});
\end{equation}

{\rm (3)}
each point $x \in \mathbb{R}^{n} \setminus S$ is contained in at most
$N=N(n)$ cubes of the family $W_{S}$.
\end{Lm}

\textbf{Proof}. The proof of Lemma~\ref{Lm2.1} is similar to that of Theorem~1 of~\cite{St2}, Ch.~6.

The family of cubes $W_{S}=\{Q_{\varkappa}\}_{\varkappa \in I}=\{Q(x_{\varkappa},r_{\varkappa})\}_{\varkappa \in I}$ constructed in Lemma~\ref{Lm2.1} is called a~\textit{Whitney decomposition} of the open set
$\mathbb{R}^{n} \setminus S$, and the cubes $Q_{\varkappa}$ are called \textit{Whitney cubes}.

Below we also need a~part of a~Whitney decomposition comprised of cubes of small side length.
More precisely, put

\begin{equation}
\notag
\mathcal{I}:=\{\varkappa \in I| r_{\varkappa} \le 1\}, \quad \mathcal{W}_{S}=\{Q_{\varkappa}\}_{\varkappa \in \mathcal{I}}.
\end{equation}

For a~cube $Q \subset \mathbb{R}^{n}$ define $Q^{\ast}:=\frac{9}{8}Q$.

\begin{Lm}
\label{Lm2.2}
For $Q_{\varkappa}, Q_{\varkappa'} \in W_{S}$ with $Q^{\ast}_{\varkappa} \cap Q^{\ast}_{\varkappa'} \neq \emptyset$
the following claims hold:

{\rm (1)}
\begin{equation}
\label{eq2.5}
\frac{1}{4}\operatorname{diam}(Q_{\varkappa})
\le \operatorname{diam}(Q_{\varkappa'})
\le  4\operatorname{diam}(Q_{\varkappa}),
\end{equation}

{\rm (2)}
for each index $\varkappa \in I$ there are at most
$C(n)$ indexes $\varkappa'$ such that
$Q^{\ast}_{\varkappa} \cap Q^{\ast}_{\varkappa'} \neq \emptyset$,

{\rm (3)}
for $\varkappa,\varkappa' \in I$ we have
$Q^{\ast}_{\varkappa} \cap Q^{\ast}_{\varkappa'} \neq \emptyset$
if and only if $Q_{\varkappa} \cap Q_{\varkappa'} \neq \emptyset$.
\end{Lm}

\textbf{Proof.} In essence, the proof is contained in that of Theorem~1 of~\cite{St2}, Ch.~6.
We leave the details to the reader.

The following \textbf{notation} is useful below. Given a~fixed closed set~$S$, for every $\varkappa \in I$ put
$$
b(Q_{\varkappa}):=b(\varkappa):=
\{\varkappa' \in I: Q_{\varkappa} \cap Q_{\varkappa'} \neq \emptyset\}
=\{\varkappa' \in I: Q^{\ast}_{\varkappa} \cap Q^{\ast}_{\varkappa'} \neq \emptyset\}.
$$

Call a~cube $Q_{\varkappa'}$ \textit{neighboring} to a~cube $Q_{\varkappa}$ if $\varkappa' \in b(Q_{\varkappa})$.
Similarly, put $b(x):=\{\varkappa \in I: Q^{\ast}_{\varkappa} \ni x\}$ for every $x \in \mathbb{R}^{n} \setminus S$.

Below we need a~special partition of unity on $\mathbb{R}^{n} \setminus S$.

\begin{Lm}
\label{Lm2.3}
For a~closed set $S \subset \mathbb{R}^{n}$, take a~Whitney decomposition
$\{Q_{\varkappa}\}_{\varkappa \in I}$ of the open set $\mathbb{R}^{n} \setminus S$ constructed in Lemma~\ref{Lm2.1}.
Then there exists a~family of functions $\{\varphi_{\varkappa}\}_{\varkappa \in I}$ with the following properties:

{\rm (1)}
$\varphi_{\varkappa} \in C^{\infty}_{0}(\mathbb{R}^{n} \setminus S)$
for every
$\varkappa \in I$;

{\rm (2)}
$0 \le \varphi_{\varkappa} \le 1$
and
$\operatorname{supp}\varphi_{\varkappa} \subset (Q_{\varkappa})^{\ast}:=\frac{9}{8}Q_{\varkappa}$
for every
$\varkappa \in I$;

{\rm (3)}
$\sum\limits_{\varkappa \in I}\varphi_{\varkappa}(x) = 1$
for all
$x \in \mathbb{R}^{n} \setminus S$;

{\rm (4)}
$\|D^{\alpha} \varphi_{\varkappa}|L_{\infty}(\mathbb{R}^{n})\| \le C(\operatorname{diam}Q_{\varkappa})^{-|\alpha|}$
for every multi-index $\alpha \in \mathbb{N}^{n}_{0}$ and every
$\varkappa \in I$ with the constant $C > 0$ depending only on~$n$.
\end{Lm}

\textbf{Proof}. See \cite{St2}, Ch.~6.

\begin{Def}
\label{Def2.7}
Given a~closed nonempty set~$S$ and $x \notin S$, say that~$\widetilde{x}$ is a~\textit{nearest point} to~%
$x$ or a~\textit{metric projection} of~$x$
to~$S$ whenever $\operatorname{dist}(x,S) = \operatorname{dist}(x,\widetilde{x})$.
\end{Def}

\begin{Remark}
\label{Rem2.5''}
Let $\widetilde{x}$ be a metric projection of $x \in \mathbb{R}^{n} \setminus S$ to $S$. Consider the interval
\begin{equation}
\notag
[x,\widetilde{x}]:=\{y=x+t(\widetilde{x}-x)|t \in [0,1]\}.
\end{equation}
Consider an arbitrary $r \in (0,\|x-\widetilde{x}\|)$. Consider the point $y_{r} = \partial Q(\widetilde{x},r) \cap [x,\widetilde{x}]$. Recall that we measure distance in the uniform norm. Show that $\operatorname{dist}(y_{r},S) = \|y_{r}-\widetilde{x}\|_{\infty} = r$. Clearly $\operatorname{dist}(y_{r},S) \le r$ because $y_{r} \in \partial Q (\widetilde{x},r)$. Assume that $\operatorname{dist}(y_{r},S) < r$. Then there is a point $y' \in S$ such that $\|y_{r}-y'\| < r$. Thus $\operatorname{dist}(x,S) \le \|x-y'\|_{\infty} \le \|x-y_{r}\|_{\infty} + \|y_{r}-y'\|_{\infty} < \|x-\widetilde{x}\|$. This contradicts to the fact that $\|\widetilde{x}-x\|=\operatorname{dist}(x,S)$.
\end{Remark}

\begin{Def}
\label{Def2.8}
Fix a~closed nonempty set~$S$. For a~cube $Q=Q(x,r) \subset \mathbb{R}^{n}$
with $x \notin S$ call $\widetilde{Q}=\widetilde{Q}(\widetilde{x},r)$
a~\textit{reflected cube}, where~$\widetilde{x}$
is a~metric projection of~$x$ to~$S$.
\end{Def}

\begin{Remark}
\label{Rem2.5}
Clearly, the metric projection exits may not be unique.
We specify an~algorithm for choosing~$\widetilde{x}$ only when our constructions require that.
Otherwise, for a~given cube $Q(x,r)$ we fix one arbitrarily chosen point~%
$\widetilde{x}$ and the cube $\widetilde{Q}(\widetilde{x},r)$.
\end{Remark}

\begin{Lm}
\label{Lm2.4}
Take a~closed set~$S$ with a~Whitney decomposition $W_{S}=\{Q_{\varkappa}\}_{\varkappa \in I}$. Let $c \geq 1$.
Then every point $x \in \mathbb{R}^{n}$ belongs at most $C$ cubes $Q(\widetilde{x}_{\varkappa},cr_{\varkappa})$ with the same side length.
The constant $C > 0$ depends only on $c$ and $n$.
\end{Lm}

\textbf{Proof}. Suppose that $Q(\widetilde{x}_{\varkappa},c r_{\varkappa}) \cap Q(\widetilde{x}_{\varkappa'},c r_{\varkappa'}) \neq \emptyset$
for some $\varkappa,\varkappa' \in I$ and $2 r_{\varkappa} = \operatorname{diam}Q_{\varkappa}=\operatorname{diam}Q_{\varkappa'}=2 r_{\varkappa'}$.
In view of \eqref{eq2.4}, we have
$\operatorname{dist}(Q_{\varkappa}, \widetilde{x}_{\varkappa}) \le 4 \operatorname{diam}(Q_{\varkappa})$
and $\operatorname{dist}(Q_{\varkappa'}, \widetilde{x}_{\varkappa'}) \le 4 \operatorname{diam}(Q_{\varkappa'})$;
hence, $\operatorname{dist}(Q_{\varkappa},Q_{\varkappa'}) \le (8+c) \operatorname{diam}(Q_{\varkappa})$.
Clearly, if $\operatorname{dist}(Q_{\alpha},Q_{\alpha'}) < (8+c) \operatorname{diam}(Q_{\varkappa})$
then $Q_{\varkappa'} \subset (18+2c)Q_{\varkappa}$. Therefore, the number of \textit{Whitney cubes} of the same size as
$Q_{\varkappa}$ lying at a~distance of less than $ (8+c) \operatorname{diam}(Q_{\varkappa})$ (recall that Whitney cubes have mutually disjoint interiors)
is bounded above by the constant $C= \frac{\mathcal{L}_{n}((18+2c)Q_{\varkappa})}{\mathcal{L}_{n}(Q_{\varkappa})}=(18+2c)^{n}$.
This proves Lemma~\ref{Lm2.4}.

\begin{Lm}
\label{Lm2.5}
Let $S\subset \mathbb{R}^{n}$ be an arbitrary nonempty closed set. Let $\mathfrak{m}$ be a finite Borel measure with $\operatorname{supp}\mathfrak{m} \subset S$. Let $W_{S}=\{Q_{\varkappa}\}_{\varkappa \in I}$ be the Whitney decomposition of $\mathbb{R}^{n} \setminus S$. Then for every $c \geq 1$
\begin{equation}
\notag
\sum\limits_{\varkappa \in \mathcal{I}}\mathcal{L}_{n}(Q(\widetilde{x}_{\varkappa},r_{\varkappa}))\mathfrak{m}(Q(\widetilde{x}_{\varkappa},c)) \le C \mathfrak{m}(S),
\end{equation}
where the constant $C > 0$ depends only on $c$ and $n$.
\end{Lm}

\textbf{Proof.} Consider the family of cubes $\{Q(\widetilde{x}_{\varkappa},c)\}_{\varkappa \in \mathcal{I}}$. Using Theorem \ref{Th2.5} and Remark \ref{Rem2.7}, it is easy to find an index set $\widehat{\mathcal{I}} \subset \mathcal{I}$ such that all cubes from the
family $\{Q(\widetilde{x}_{\varkappa},c)\}_{\varkappa \in \widehat{\mathcal{I}}}$, are mutually disjoint and
\begin{equation}
\notag
\cup_{\varkappa \in \widehat{\mathcal{I}}}Q(\widetilde{x}_{\varkappa}, 5c) \supset \cup_{\varkappa \in \mathcal{I}}Q(\widetilde{x}_{\varkappa}, c).
\end{equation}

Note that if $Q(\widetilde{x}_{\varkappa'},r_{\varkappa'}) \cap Q(\widetilde{x}_{\varkappa},5c) \neq \emptyset$ for some $\varkappa, \varkappa' \in \mathcal{I}$, then
\begin{equation}
\label{incl}
Q(\widetilde{x}_{\varkappa'},r_{\varkappa'}) \subset Q(\widetilde{x}_{\varkappa},7c)
\end{equation}
because $c \geq 1$ and $r_{\varkappa'} \le 1$.

Using \eqref{eq2.4}, we conclude that the Whitney cube $Q_{\varkappa'}=Q(x_{\varkappa'},r_{\varkappa'}) \subset Q(x_{\varkappa}, 20c)$. Hence, using the fact that different Whitney cubes have disjoint interiors, we get
\begin{equation}
\begin{split}
\notag
&\sum\limits_{\substack{\varkappa' \in \mathcal{I} \\
Q(\widetilde{x}_{\varkappa'},r_{\varkappa'}) \cap Q(\widetilde{x}_{\varkappa},5c) \neq \emptyset}}\mathcal{L}_{n}(Q(\widetilde{x}_{\varkappa'},r_{\varkappa'})) \le \sum\limits_{\substack{\varkappa' \in \mathcal{I} \le Q_{\varkappa'} \subset Q(\widetilde{x}_{\varkappa},20c) \neq \emptyset}}\mathcal{L}_{n}(Q(\widetilde{x}_{\varkappa'},r_{\varkappa'}))\\
&\le \mathcal{L}_{n}(Q(x_{\varkappa},20c)) \le (20 c)^{n}.
\end{split}
\end{equation}

Using this fact, inclusion \eqref{incl} and Lemmas \ref{Lm2.8}, \ref{Lm2.9}, we obtain the estimate
\begin{equation}
\begin{split}
\notag
&\sum\limits_{\varkappa \in \mathcal{I}}\mathcal{L}_{n}(Q(\widetilde{x}_{\varkappa},r_{\varkappa}))\mathfrak{m}(Q(\widetilde{x}_{\varkappa},c)) \le
\sum\limits_{\varkappa \in \widehat{\mathcal{I}}}\sum\limits_{\substack{\varkappa' \in \mathcal{I} \\
Q(\widetilde{x}_{\varkappa'},r_{\varkappa'}) \cap Q(\widetilde{x}_{\varkappa},5c) \neq \emptyset}}\mathcal{L}_{n}(Q(\widetilde{x}_{\varkappa'},r_{\varkappa'}))\mathfrak{m}(Q(\widetilde{x}_{\varkappa'},c)) \\
& \le \sum\limits_{\varkappa \in \widehat{\mathcal{I}}}\mathfrak{m}(Q(\widetilde{x}_{\varkappa},7c))\sum\limits_{\substack{\varkappa' \in \mathcal{I} \\ Q(\widetilde{x}_{\varkappa'},r_{\varkappa})
\cap Q(\widetilde{x}_{\varkappa},5c) \neq \emptyset}}\mathcal{L}_{n}(Q(\widetilde{x}_{\varkappa'},r_{\varkappa'}))\\
& \le (20c)^{n} \sum\limits_{\varkappa \in
\widehat{\mathcal{I}}}\mathfrak{m}(Q(\widetilde{x}_{\varkappa},7c)) \le (20c)^{n} \mathfrak{m}(S).
\end{split}
\end{equation}

The lemma is proved.






Recall the notion of Ahlfors~$n$-regular sets (see Definition \ref{Def2.2}).

\begin{Lm}
\label{Lm2.6}
Let $S$ be a~closed Ahlfors~ $n$-regular set in $\mathbb{R}^{n}$. Let $\mathcal{W}_{S}=\{Q_{\varkappa}\}_{\varkappa \in \mathcal{I}}$ be the part of the Whitney decomposition of $\mathbb{R}^{n} \setminus S$ comprised of cubes of side length $\le 1$. Then there exists a~family
$\mathfrak{U}:=\{\mathcal{U}_{\varkappa} : \varkappa \in \mathcal{I}\}$
of Borel sets with the following properties:

{\rm (1)}
$\mathcal{U}_{\varkappa} \subset Q(\widetilde{x}_{\varkappa}, r_{\varkappa}) \subset (10Q_{\varkappa}) \cap S$
for all
$\varkappa \in \mathcal{I}$;

{\rm (2)}
$\mathcal{L}_{n}(Q_{\varkappa}) \le \kappa_{1}\mathcal{L}_{n}(\mathcal{U}_{\varkappa})$
for all
$\varkappa \in \mathcal{I}$;

{\rm (3)}
$\sum\limits_{\varkappa \in \mathcal{I}}\chi_{\mathcal{U}_{\varkappa}}(x) \le \kappa_{2}$
for
$x \in S$.

\noindent
Furthermore, the positive constants $\kappa_{1}$ and $\kappa_{2}$ depend only on $n$~and~the constants $C_{1},C_{2}$ from Definition \ref{Def2.2}.
\end{Lm}

\textbf{Proof}. Our arguments repeat almost verbatim the proof of Theorem 2.4 of~\cite{Shv1}.

\subsection{Covering Theorems}

\begin{Def}
\label{Def2.9}
A collection $\mathcal{F}$ of closed balls in $\mathbb{R}^{n}$ is a cover of a set $E \subset \mathbb{R}^{n}$ if
\begin{equation}
\notag
E \subset \bigcup\limits_{B \in \mathcal{F}}B.
\end{equation}
\end{Def}

The following theorems are classical. One can find the proofs in section 1.5 of \cite{Evans}.

\begin{Th}(Vitali's Covering Theorem)
\label{Th2.5}
Let $\mathcal{F}$ be any collection of closed nondegenerate balls in $\mathbb{R}^{n}$ with
\begin{equation}
\notag
\sup\{\operatorname{diam}B | B \in \mathcal{F}\} < \infty.
\end{equation}

Then there exists a countable family $\mathcal{G} \subset \mathcal{F}$ of disjoint balls such that

\begin{equation}
\label{eq2.6}
\bigcup\limits_{B \in \mathcal{F}} B \subset \bigcup\limits_{B \in \mathcal{G}} 5B.
\end{equation}
\end{Th}

\begin{Remark}
\label{Rem2.7}
Similarly we can define a cover of a given set $E \subset \mathbb{R}^{n}$ by closed cubes instead of closed balls. One can formulate and prove analog of the previous theorem using cubes instead of balls.
\end{Remark}

\begin{Th}(Besicovitch's Covering Theorem)
\label{Th2.6}
There exists a constant $N(n)$, depending only on $n$ with the following ptoperty:
If $\mathcal{F}$ is any collection of nondegenerate closed balls in $\mathbb{R}^{n}$ with
\begin{equation}
\notag
\sup\{\operatorname{diam}B | B \in \mathcal{F}\} < \infty
\end{equation}

and if $A$ is the set of centers of balls in $\mathcal{F}$, then there exist $\mathcal{G}_{1},...,\mathcal{G}_{N(n)} \subset \mathcal{F}$ such that
each $\mathcal{G}_{i}$ ($i \in \{1,...,N(n)\}$) is a countable collection of disjoint balls in $\mathcal{F}$ and

\begin{equation}
\notag
A \subset \bigcup\limits_{i=1}^{N(n)}\bigcup\limits_{B \in \mathcal{G}_{i}}B.
\end{equation}

\end{Th}

\begin{Def}
\label{Def2.10}
Let $E$ be a nonempty set in $\mathbb{R}^{n}$. Let $\varepsilon > 0$. Let $\{x_{j}\}_{j \in \mathcal{J}}$, $\mathcal{J} \subset \mathbb{N}$ be a subset of $E$ with
the following properties:

{\rm (i)} $\|x_{i}-x_{j}\|_{\infty} \geq \varepsilon$ for every $i,j \in \mathcal{J}$ and $i \neq j$;

{\rm (ii)} for every $x \in E \setminus \{x_{j}\}_{j \in \mathcal{J}}$ there is a point $x_{j}$ such that $\|x-x_{j}\|_{\infty} < \varepsilon$.

We call the set $\{x_{j}\}_{j \in \mathcal{J}}$ \textit{maximal $\varepsilon$-separated subset of} $E$.
\end{Def}

The following result is a direct corollary of Definition \ref{Def2.10}. Recall that all cubes are assumed to be closed.

\begin{Lm}
\label{Lm2.7}
Let $E$ be a nonempty  set in $\mathbb{R}^{n}$. Let $\varepsilon > 0$.
Let $\{x_{j}\}_{j \in \mathcal{J}}$ be a maximal $\varepsilon$-separated subset of $E$. Then

{\rm (1)} $E \supset \bigcup\limits_{j \in \mathcal{J}}Q(x_{j},\varepsilon)$;

{\rm (2)} the family $\{Q(x_{j},\frac{\varepsilon}{2})\}_{j \in \mathcal{J}}$ is pairwise disjoint;

{\rm (3)} every point $x \in E$ belongs at most $3^{n}$ cubes in $\{Q(x_{j},\frac{\varepsilon}{2})\}_{j \in \mathcal{J}}$.

\end{Lm}

\textbf{Proof.} We prove only item (3), because the other items are obvious. Fix an arbitrary point $x \in E$. If $x \in Q(x_{j},\varepsilon)$, then
$Q(x_{j}, \frac{\varepsilon}{2}) \subset Q(x,\frac{3\varepsilon}{2})$. It is clear that the cube $Q(x,\frac{3\varepsilon}{2})$ contains at most $3^{n}$ mutually disjoint cubes
with diameters $\varepsilon$. Hence, using item (2), we conclude.

\begin{Def}
\label{Def2.10}
Let $\mathcal{J}$ be an arbitrary finite or countably index set. Let $\{E_{j}\}_{j \in \mathcal{J}}$ be a family of Borel subsets of $\mathbb{R}^{n}$. We say that the multiplicity of overlapping of the sets $E_{j}$ is finite if and only if
there exists a number $N \in \mathbb{N}$ such that every point $x \in \cup_{j}E_{j}$ belongs at most than $N$ sets from the family $\{E_{j}\}_{j \in \mathcal{J}}$.
\end{Def}

The following simple lemmas will be often useful in what follows.

\begin{Lm}
\label{Lm2.8}
Let $r > 0$ and $c \geq 1$. Let $\{Q_{j}\}_{j \in \mathcal{J}}=\{Q(x_{j},r)\}_{j \in \mathcal{J}}$ be a family of mutually disjoint cubes with the same side length.
Then the multiplicity of overlapping of the cubes $cQ_{j}$ is finite and bounded above by a constant $C > 0$ depending only on $n$ and $c$.
\end{Lm}

\textbf{Proof.} Assume that $cQ_{j} \cap cQ_{j'} \neq \emptyset$. Then, using the fact that $\operatorname{diam}Q_{j}=\operatorname{diam}Q_{j'}$, we
conclude that $cQ_{j'} \subset 3c Q_{j}$. Hence the number of cubes $cQ_{j'}$ which have nonempty intersection with $cQ_{j}$ is bounded above by the number
of mutually disjoint cubes $Q_{j'}$ containing in $3c Q_{j}$. But the later is at most $(3c)^{n}$. This proves the lemma.

\begin{Lm}
\label{Lm2.9}
Let $\mathfrak{m}$ be a finite Borel measure on $\mathbb{R}^{n}$. Let $\{E_{j}\}_{j \in \mathcal{J}}$ be a family of Borel subsets of $\mathbb{R}^{n}$ such that
the multiplicity of overlapping of the sets $E_{j}$ is finite and bounded above by some constant $N \in \mathbb{N}$. Then
\begin{equation}
\label{eq2.9}
\sum\limits_{j \in \mathcal{J}}\mathfrak{m}(E_{j}) \le N \mathfrak{m}(\mathbb{R}^{n}).
\end{equation}
\end{Lm}

\textbf{Proof.} From the hypothesis of the theorem we see at once that

\begin{equation}
\notag
\sum\limits_{j \in \mathcal{J}}\chi_{E_{j}}(x) \le N, \quad x \in \mathbb{R}^{n}.
\end{equation}

Hence, we obtain the desirable estimate

\begin{equation}
\begin{split}
\notag
&\sum\limits_{j \in \mathcal{J}}\mathfrak{m}(E_{j}) = \sum\limits_{j \in \mathcal{J}}\int\limits_{\mathbb{R}^{n}}\chi_{E_{j}}(x)\,d\mathfrak{m}(x)\\
&=\int\limits_{\mathbb{R}^{n}}\sum\limits_{j \in \mathcal{J}}\chi_{E_{j}}(x)\,d\mathfrak{m}(x) \le N \mathfrak{m}(\mathbb{R}^{n}).
\end{split}
\end{equation}


\section{Main tools}

The aim of this section is to present key tools which are cornerstones for the proof of our main results. More precisely, we introduce $d$-\textit{regular systems of measures}, \textit{generalized Calderon-type maximal functions}, \textit{porous sets} and establish basic properties of these objects.

\subsection{$d$-regular system of measures}

The following result is a variant of the Frostman-type theorem adapted for our purposes (compare with Theorem 5.1.12 in \cite{A}). For convenience of the reader we will present the proof in Appendix.

\begin{Th}
\label{Th3.1}
Let $d \in [0,n]$. Let $S$ be a nonempty closed subset of $\mathbb{R}^{n}$. Then, there exists a system of  Borel measures $\{\mu_{k}\}_{k \in \mathbb{N}_{0}}$ such that for every $k \in \mathbb{N}_{0}$ the following properties hold:

{\rm (1)}
$$
\operatorname{supp}\mu_{k} \subset E;
$$

{\rm (2)}
\begin{equation}
\label{eq3.1}
\mu_{k}(B(x,r)) \le C r^{d}, \quad x \in \mathbb{R}^{n}, \quad r \in (0,2^{k}];
\end{equation}

{\rm (3)}
for every set $V$ of the form $V=\cup_{m \in \mathcal{A}} Q_{k,m}$, $\mathcal{A} \subset \mathbb{Z}^{n}$

\begin{equation}
\label{eq3.2}
\mu_{k}(V \cap S) \geq C \mathcal{H}^{d}_{\infty}(V \cap S),
\end{equation}

where the constant $C > 0$ depends only on $n$;

{\rm (4)}
for every Borel set $G \subset \mathbb{R}^{n}$

\begin{equation}
\label{eq3.3}
2^{d-n}\mu_{k}(G) \le \mu_{k-1}(G) \le \mu_{k}(G).
\end{equation}
\end{Th}

\begin{Def}
\label{Def3.1}
Let $d \in [0,n]$. Assume that $\operatorname{dim}_{H}S \geq d$. Let $\{\mu_{k}\}_{k \in \mathbb{N}_{0}}$ be a family of Borel measures on $\mathbb{R}^{n}$ with $\operatorname{supp}\mu_{k} \subset S$, $k \in \mathbb{N}_{0}$.
We say that $\{\mu_{k}\}_{k \in \mathbb{N}_{0}}$ is  a \textit{$d$-regular system of measures on} $S$ if and only if for some universal constants $C_{1}, C_{2}, C_{3} > 0$ the following properties hold for every $k \in \mathbb{N}_{0}$:

{\rm (1)}
\begin{equation}
\label{eq3.4}
\mu_{k}(B(x,r)) \le C_{1} r^{d} \text{ for every } x \in \mathbb{R}^{n}
\text{ and every } r \in (0,2^{-k}];
\end{equation}

{\rm (2)}
\begin{equation}
\label{eq3.5}
\mu_{k}(B(x,2^{-k})) \geq C_{2} 2^{-k} \text{ for every } x \in S;
\end{equation}

{\rm (3)}
\begin{equation}
\label{eq3.6}
2^{d-n}\mu_{k}(G) \le \mu_{k-1}(G) \le \mu_{k}(G) \text{ for every Borel set } G \subset S.
\end{equation}
\end{Def}

\begin{Remark}
\label{Rem3.1}
Note that in view of \eqref{eq3.6} a function $f \in L_{1}^{\text{loc}}(\mathbb{R}^{n},\mu_{k})$ for some fixed $k \in \mathbb{N}_{0}$ if and only if
$L_{1}^{\text{loc}}(\mathbb{R}^{n},\mu_{j})$ for every $j \in \mathbb{N}_{0}$. Furthermore, given a number $c \in [2^{j},2^{j+1})$, $j \in \mathbb{N}_{0}$, estimates \eqref{eq3.4}--\eqref{eq3.6} implies that for every $k \in \mathbb{N}_{0}$ and every $x \in S$
\begin{equation}
\begin{split}
\label{doubling}
&\mu_{k}(Q(x,\frac{2^{-k}}{c})) \geq 2^{(j+1)(d-n)}\mu_{k+j+1}(Q(x,\frac{2^{-k}}{c})) \geq 2^{(d-n)(j+1)}\mu_{k+j+1}(Q(x,\frac{2^{-k}}{2^{j+1}}))\\
&\geq \frac{C_{2}}{2^{n(j+1)}} 2^{-dk} \geq \frac{C_{2}}{C_{1}2^{n(j+1)}}\mu_{k}(Q(x,2^{-k})).
\end{split}
\end{equation}
Note that estimate \eqref{doubling} implies that the space $(S,\mu_{k}|_{S})$ is the space with doubling measure. But this fact does not imply that the measure $\mu_{k}$ is doubling on $\mathbb{R}^{n}$.
\end{Remark}

\begin{Ca}
\label{Ca3.1}
Let $d \in [0,n]$. Let $S$ be a closed $d$-thick set. Then there exists a \textit{$d$-regular system of measures on} $S$.
\end{Ca}

\textbf{Proof}. Apply Theorem \ref{Th3.1} to the set $S$. This gives the system of Borel measures $\{\mu_{k}\}_{k \in \mathbb{N}_{0}}$ on $S$. In order to prove our assertion it is sufficient to verify only \eqref{eq3.5}.

Fix some $x \in S$ and consider all dyadic cubes from $\mathcal{Q}_{k+2}$ which intersect $Q(x,2^{-k-2})$. It is clear that all such cubes are contained in $Q(x,2^{-k})$ and the union of these cubes contains $Q(x,2^{-k-2})$. Hence, using Definition \ref{Def2.3}, estimates \eqref{eq3.2}, \eqref{eq3.3}  and subadditivity of the Hausdorf content,  we obtain

\begin{equation}
\begin{split}
\notag
&2^{2(n-d)}\mu_{k}(Q(x,2^{-k})) \geq \mu_{k+2}(Q(x,2^{-k})) \geq \sum\limits_{\substack{m \in \mathbb{Z}^{n} \\ Q_{k+2,m} \cap Q(x,2^{-k-2}) \neq \emptyset}}\mu_{k+2}(Q_{k+2,m})\\
&\geq C \sum\limits_{\substack{m \in \mathbb{Z}^{n} \\ Q_{k+2,m} \cap Q(x,2^{-k-2}) \neq \emptyset}}\mathcal{H}^{d}_{\infty}(Q_{k+2,m})
\geq C \mathcal{H}^{d}_{\infty}\left(\bigcup\limits_{\substack{m \in \mathbb{Z}^{n} \\ Q_{k+2,m} \cap Q(x,2^{-k-2}) \neq \emptyset}}Q_{k+2,m}\right) \\
&\geq C \mathcal{H}^{d}_{\infty}(Q(x,2^{-k-2})) \geq C C'_{2} 2^{-2d} 2^{-kd}.
\end{split}
\end{equation}

This proves the claim.

\begin{Lm}
\label{Lm3.1}
Let $d \in [0,n]$. Let $S$ be a closed set in $\mathbb{R}^{n}$ with $\operatorname{dim}_{H}(S) \geq d$. Let $\{\mu_{k}\}_{k \in \mathbb{N}_{0}}$ be a $d$-regular system of measures on $S$. Let $E$ be a Borel subset of $S$. If
$\mathcal{H}^{d}(E) = 0$, then $\mu_{k}(E) = 0$ for every $k \in \mathbb{N}_{0}$.
\end{Lm}

\textbf{Proof.} Assume that $\mathcal{H}^{d}(E) = 0$. Then, for every $j \in \mathbb{N}$ there exist a number $\delta_{j}$ and a countable covering of $E$ by balls $\{B^{j}_{i}\}=\{B(x^{j}_{i},r^{j}_{i})\}_{i \in \mathbb{N}}$ with radii $\sup\limits_{i \in \mathbb{N}}r^{j}_{i} < \delta_{j}$ such that

\begin{equation}
\notag
\sum\limits_{i=1}^{\infty}(r^{j}_{i})^{d} < \frac{1}{j}.
\end{equation}

Using this and \eqref{eq3.4}, it is easy to see that for every $k \le j$

\begin{equation}
\label{eq3.7}
\mu_{k}(E) \le \mu_{k}(\cup_{i}B^{j}_{i}) \le \sum\limits_{i=1}^{\infty}\mu_{k}(B^{j}_{i}) \le \frac{C}{j}.
\end{equation}

Fix $k \in \mathbb{N}_{0}$ and letting $j \to \infty$ in \eqref{eq3.7} we conclude.

\begin{Remark}
\label{Rem3.2}
We will see in Example 6.3  that there exist a set $S$ with $\operatorname{dim}_{H}S \geq d$, $d \in [0,n]$, a $d$-regular system of measures $\{\mu_{k}\}_{k \in \mathbb{N}_{0}}$ on $S$ and set $E \subset S$ such that $\mathcal{H}^{d}(E) > 0$ but $\mu_{k}(E)=0$ for every $k \in \mathbb{N}_{0}$.
\end{Remark}

For every $r > 0$ by $k(r)$ we denote the unique integer number for which $r \in [2^{-k(r)},2^{-k(r)+1})$.

\begin{Lm}
\label{Lm3.2}
Let $d \in [0,n]$. Let $S$ be an arbitrary closed nonempty subset of $\mathbb{R}^{n}$ with $\operatorname{dim}_{H}S \geq d$. Let $\{\mu_{k}\}_{k \in \mathbb{N}_{0}}$ be a $d$-regular system of measures on $S$.
Then for every $r \in (0,1)$, $x \in S$ and every Borel set $G \subset Q(x,r) \cap S$
\begin{equation}
\label{eq3.8}
\frac{\mathcal{L}_{n}(G)}{\mathcal{L}_{n}(Q(x,r))} \le C \frac{\mu_{k(r)}(G)}{\mu_{k(r)}(Q(x,r) \cap S)}.
\end{equation}
The constant $C > 0$ depends only on $n$ and the constants $C_{1},C_{2}$ from Defintion \ref{Def3.1}.
\end{Lm}

\textbf{Proof.} Fix $x \in S$ and $r \in (0,1)$. Consider an arbitrary cube $Q(y,t) \subset Q(x,r)$ with $y \in S$. It is clear that

\begin{equation}
\label{eq3.9}
\frac{\mathcal{L}_{n}(Q(y,t) \cap S)}{\mathcal{L}_{n}(Q(x,r))} \le \frac{\mathcal{L}_{n}(Q(y,t) \cap S)}{\mathcal{L}_{n}(Q(y,t))}\frac{\mathcal{L}_{n}(Q(y,t))}{\mathcal{L}_{n}(Q(x,r))} \le 2^{n(k(r)-k(t)+1)}\frac{\mathcal{L}_{n}(Q(y,t) \cap S)}{\mathcal{L}_{n}(Q(y,t))}.
\end{equation}

On the other hand, using \eqref{eq3.4} -- \eqref{eq3.6} (we can use these estimates because $x,y \in S$), we have

\begin{equation}
\label{eq3.10}
\frac{\mu_{k(r)}(Q(y,t)\cap S)}{\mu_{k(r)}(Q(x,r) \cap S)} \geq 2^{(d-n)(k(t)-k(r))}\frac{\mu_{k(t)}(Q(y,t)\cap S)}{\mu_{k(r)}(Q(x,r) \cap S)} \geq C(C_{1},C_{2}) 2^{n(k(r)-k(t))}.
\end{equation}

Combining \eqref{eq3.9}, \eqref{eq3.10}, we obtain (the constants $C_{1}, C_{2} > 0$ are the same as in \eqref{eq3.4}, \eqref{eq3.5})

\begin{equation}
\begin{split}
\label{eq3.11}
&\frac{\mathcal{L}_{n}(Q(y,t) \cap S)}{\mathcal{L}_{n}(Q(x,r))} \le \widetilde{C} (C_{1},C_{2}) \frac{\mu_{k(r)}(Q(y,t)\cap S)}{\mu_{k(r)}(Q(x,r) \cap S)}\frac{\mathcal{L}_{n}(Q(y,t) \cap S)}{\mathcal{L}_{n}(Q(y,t))}\\
&\le \widetilde{C} (C_{1},C_{2}) \frac{\mu_{k(r)}(Q(y,t)\cap S)}{\mu_{k(r)}(Q(x,r) \cap S)}.
\end{split}
\end{equation}

Using $\sigma$-additivity of measures $\mathcal{L}_{n}$ and $\mu_{k(r)}$, we have (see Theorem 1, section 1 of \cite{Evans})

\begin{equation}
\label{eq3.12}
\mathcal{L}_{n}(G)=\lim\limits_{j \to \infty}\mathcal{L}_{n}(U_{j}), \quad \mu_{k(r)}(G) = \lim\limits_{j \to \infty}\mu_{k(r)}(U_{j}),
\end{equation}

where  $\{U_{j}\}$ is an arbitrary decreasing sequence of open sets such that $G=\bigcap\limits_{i=1}^{\infty}U_{i}$.

For every $j \in \mathbb{N}$ let $\{x^{j}_{i}\}$ be a maximal $\frac{r}{j}$ separated subset of $Q(x,r) \cap S$. Recall that $\{x^{j}_{i}\} \subset Q(x,r) \cap S$. Clearly $Q(x,r)\cap S \subset \cup_{i}\operatorname{int}Q(x^{j}_{i},\frac{2r}{j})$ and cubes $Q(x^{j}_{i},\frac{r}{2j})$ are pairewisely disjoint.

For every $j \in \mathbb{N}$ we consider the set

\begin{equation}
\notag
U_{j}:=\bigcup\limits_{i}\operatorname{int}Q(x^{j}_{i},\frac{2r}{j}).
\end{equation}

It is clear that $U_{j} \subset Q(x,3r)$ for every $j \in \mathbb{N}_{0}$.
Hence, using Lemma \ref{Lm2.8}  and \eqref{eq3.11}, we get

\begin{equation}
\label{eq3.13}
\begin{split}
&\mathcal{L}_{n}(U_{j}) \le \sum\limits_{i}\mathcal{L}_{n}(Q(x^{j}_{i},\frac{2r}{j})) \le C \frac{\mathcal{L}_{n}(Q(x,3r))}{\mu_{k(3r)}(Q(x,3r)\cap S)}\sum\limits_{i}\mu_{k(r)}(B(x^{j}_{i},\frac{2r}{j})\cap S) \\
&\le   C \frac{\mathcal{L}_{n}(Q(x,r))}{\mu_{k(r)}(Q(x,r)\cap S)}\mu_{k(r)}(U_{j}).
\end{split}
\end{equation}

Combining \eqref{eq3.12} and \eqref{eq3.13}, we complete the proof.

\begin{Ca}
\label{Ca3.2}
Let $d \in [0,n]$. Let $S$ be an arbitrary closed subset of $\mathbb{R}^{n}$ with $\operatorname{dim}_{H}S \geq d$. Let $\{\mu_{k}\}_{k \in \mathbb{N}_{0}}$
be a $d$-regular system of measures on $S$. Assume that a function $f \in L_{1}^{\operatorname{loc}}(S,\mu_{k})$ for every $k \in \mathbb{N}_{0}$. Then for every $x \in S$
and every $r \in (0,1)$
\begin{equation}
\label{eq3.14'}
\frac{1}{\mathcal{L}_{n}(Q(x,r))}\int\limits_{Q(x,r) \cap S}|f(y)|\,d\mathcal{L}_{n}(y) \le C \fint\limits_{Q(x,r) \cap S}|f(y)|\,d\mu_{k(r)}(y).
\end{equation}
The constant $C > 0$ does not depend on $x$, $r$, $f$.
\end{Ca}

\textbf{Proof.} For a simple function $f:S \to \mathbb{R}$  estimate \eqref{eq3.14'} clearly holds due to Lemma \ref{Lm3.2}. In general case we should construct increasing sequence of simple functions converging to $f$ and use monotone convergence theorem for integrals (see section 1.3 of \cite{Evans}).

\begin{Lm}
\label{Lm3.3}
Let $c \geq 1$ and $d \in [0,n]$.  Let $S$ be an arbitrary closed nonempty set in $\mathbb{R}^{n}$ with $\operatorname{dim}_{H}S \geq d$. Assume that there exists a $d$-regular system of measures $\{\mu_{k}\}_{k \in \mathbb{N}_{0}}$ on $S$. Assume that  $g \in
L_{1}^{\operatorname{loc}}(\mathbb{R}^{n},\mu_{k})$ for every $k \in \mathbb{N}_{0}$. Then for every two
cubes $Q=Q(x,r)$, $Q'=Q(x',cr)$ with $x \in S$, $x' \in \mathbb{R}^{n}$, $cr \in (0,1)$ and $Q \subset Q'$
\begin{equation}
\label{eq3.15}
\begin{split}
&\fint\limits_{Q(x,r)}|g(z)|\,d\mu_{k(r)}(z) \le C \fint\limits_{Q(x',cr)}|g(z)|\,d\mu_{k(r)}(z), \\
&\fint\limits_{Q(x,r)}|g(z)|\,d\mu_{k(r)}(z) \approx \fint\limits_{Q(x,r)}|g(z)|\,d\mu_{k(r)+1}(z),
\end{split}
\end{equation}

where the constant $C > 0$ does not depend on $x,x'$ and $r$.

\end{Lm}

\textbf{Proof}. Let us prove the first inequality in \eqref{eq3.15}. Let $N(c)$ be a number of all dyadic cubes with side length $2^{-k}$ which have nonempty
intersection with $Q(x,cr)$. Hence, using \eqref{eq3.4} we obtain

\begin{equation}
\label{eq3.16}
\mu_{k}(Q(x',cr)) \le C_{1}N(c)r^{d}.
\end{equation}

From \eqref{eq3.5}, \eqref{eq3.16}, using inclusion $Q(x,r) \subset Q(x',cr)$, we derive

\begin{equation}
\begin{split}
\notag
&\fint\limits_{Q(x,r)}|g(y)|\,d\mu_{k}(y) \le \frac{\mu_{k}(Q(x',cr))}{\mu_{k}(Q(x,r))}\fint\limits_{Q(x',cr)}|g(y)|\,d\mu_{k}(y) \le C \fint\limits_{Q(x',cr)}|g(y)|\,d\mu_{k}(y).
\end{split}
\end{equation}

The second inequality in \eqref{eq3.15} clearly follows from \eqref{eq3.6}.
The lemma is proved.

\subsection{Calderon-type maximal functions}

Now we introduce one of the main tool for this paper. In the sequel we will often use the following notation. Given $r > 0$, we denote by $k(r)$ the unique
integer number for which $r \in [2^{-k(r)},2^{-k(r)+1})$.

\begin{Def}
\label{Def3.2}
Let $d \in [0,n]$. Let $S$ be an arbitrary closed nonempty set in $\mathbb{R}^{n}$ with $\operatorname{dim}_{H}S \geq d$. Assume that there exists a $d$-regular system of measures $\{\mu_{k}\}=\{\mu_{k}\}_{k \in \mathbb{N}_{0}}$ on $S$. We define for every $t \in [0,1)$ \textit{generalized Calderon-type maximal function} as follows
\begin{equation}
\begin{split}
\notag
&f^{\sharp}_{\{\mu_{k}\}}(x,t):=\sup\limits_{r \in (t,1)} \frac{1}{r}\fint\limits_{Q(x,r)}\Bigl|f(y)-\fint\limits_{Q(x,r)}f(z)\,d\mu_{k(r)}(z)\Bigr|\,d\mu_{k(r)}(y)\\
&=\sup\limits_{r \in (t,1)} \frac{1}{r}\fint\limits_{Q(x,r) \cap S}\Bigl|f(y)-\fint\limits_{Q(x,r) \cap S}f(z)\,d\mu_{k(r)}(z)\Bigr|\,d\mu_{k(r)}(y).
\end{split}
\end{equation}
\end{Def}

\begin{Remark}
\label{Rem3.2}
In the case $t=0$ we will write $f^{\sharp}_{\{\mu_{k}\}}(x)$ instead of $f^{\sharp}_{\{\mu_{k}\}}(x,0)$. If the set $S$ is Ahlfors $n$-regular one can take $\mu_{k}=\mathcal{L}_{n}$ for every $k \in \mathbb{N}_{0}$. Hence, for such set $S$ our function $f^{\sharp}_{\{\mu_{k}\}}$ coincide with that introduced in \cite{Shv1} (it was denoted there by $f^{\sharp}_{S}$). In particular, if $S=\mathbb{R}^{n}$ we obtain \textit{the classical Calderon-type maximal function} \cite{Cal}.
\end{Remark}

\begin{Lm}
\label{Lm3.4}
Let $c \geq 1$, $d \in [0,n]$. Let $S$ be an arbitrary closed nonempty set in $\mathbb{R}^{n}$ with $\operatorname{dim}_{H}S \geq d$. Assume that there exists a $d$-regular system of measures $\{\mu_{k}\}=\{\mu_{k}\}_{k \in \mathbb{N}_{0}}$ on $S$. Assume that
$f \in L_{1}^{\rm loc}(\mathbb{R}^{n},\mu_{k})$ for every $k \in \mathbb{N}_{0}$.
Then for every  pair of cubes $Q=Q(x,r)$ and $Q'=Q(x',cr)$ such that $x,x' \in S$, $cr \in (0,1)$, $Q \subset Q'$
\begin{equation}
\label{eq3.17}
\begin{split}
&\Bigl|\fint\limits_{Q(x,r)}f(y)\,d\mu_{k(r)}(y)-\fint\limits_{Q(x',cr)}f(z)\,d\mu_{k(cr)}(z)\Bigr| \\
&\le C\fint\limits_{Q(x',cr)}\Bigl|f(y)-\fint\limits_{Q(x',cr)}f(z)\,d\mu_{k(cr)}(z)\Bigr|\,d\mu_{k(cr)}(y),
\end{split}
\end{equation}
where the constant $C > 0$ does not depend on $x,x'$ and $r$.
\end{Lm}

\textbf{Proof}. Using Lemma \ref{Lm3.3} and \eqref{eq3.6} we clearly have

\begin{equation}
\notag
\begin{split}
&\Bigl|\fint\limits_{Q(x,r)}f(y)\,d\mu_{k(r)}(y)-\fint\limits_{Q(x',cr)}f(z)\,d\mu_{k(cr)}(z)\Bigr| \\
& \le \fint\limits_{Q(x,r)}\Bigl|f(y)-\fint\limits_{Q(x',cr)}f(z)\,d\mu_{k(cr)}(z)\Bigr|\,d\mu_{k(r)}(y)\\
&\le C \fint\limits_{Q(x',cr)}\Bigl|f(y)-\fint\limits_{Q(x',cr)}f(z)\,d\mu_{k(cr)}(z)\Bigr|\,d\mu_{k(r)}(y)\\
&\le C\fint\limits_{Q(x',cr)}\Bigl|f(y)-\fint\limits_{Q(x',cr)}f(z)\,d\mu_{k(cr)}(z)\Bigr|\,d\mu_{k(cr)}(y).
\end{split}
\end{equation}

The lemma is proved.

\begin{Lm}
\label{Lm3.5}
Let $c \geq 1$, $d \in [0,n]$. Let $S$ be an arbitrary closed nonempty set in $\mathbb{R}^{n}$ with $\operatorname{dim}_{H}S \geq d$. Assume that there exists a $d$-regular system of measures $\{\mu_{k}\}=\{\mu_{k}\}_{k \in \mathbb{N}_{0}}$ on $S$. Assume that
$f \in L_{1}^{\rm loc}(\mathbb{R}^{n},\mu_{k})$ for every $k \in \mathbb{N}_{0}$. Let $Q(x,r) \subset Q(x',cr)$ for some $x,x' \in S$ and $r \in (0,1)$. Then
\begin{equation}
\notag
f^{\sharp}_{\{\mu_{k}\}}(x,r) \le C \Bigl(f^{\sharp}_{\{\mu_{k}\}}(x',r)+\fint\limits_{Q(x',c)}|f(y)|\,d\mu_{0}(y)\Bigr),
\end{equation}

where the constant $C > 0$ depends only on $c$ and $\{\mu_{k}\}_{k \in \mathbb{N}_{0}}$.
\end{Lm}

\textbf{Proof.} Assume that $rc < tc < 1$. Then, using Lemma \ref{Lm3.3} and inclusion $Q(x,t) \subset Q(x',ct)$ (which clearly follows from inclusion $Q(x,r) \subset Q(x',cr)$), we obtain

\begin{equation}
\begin{split}
\notag
&t^{-1}\fint\limits_{Q(x,t)}\Bigl|f(y)-\fint\limits_{Q(x,t)}f(z)\,d\mu_{k(t)}(z)\Bigr|\,d\mu_{k(t)}(y)
\\ &\le t^{-1}\fint\limits_{Q(x,t)}\fint\limits_{Q(x,t)}\Bigl|f(y)-f(z)\Bigr|\,d\mu_{k(t)}(z)\,d\mu_{k(t)}(y)  \\
&\le \frac{C}{ct}\fint\limits_{Q(x',ct)}\fint\limits_{Q(x',ct)}\Bigl|f(y)-f(z)\Bigr|\,d\mu_{k(t)}(z)\,d\mu_{k(t)}(y) \\
&\le \frac{C}{ct}\fint\limits_{Q(x',ct)}\Bigl|f(y)-\fint\limits_{Q(x',ct)}f(z)\,d\mu_{k(t)}(z)\Bigr|\,d\mu_{k(t)}(y).
\end{split}
\end{equation}

Hence, for every $r \in (0,\frac{1}{c})$

\begin{equation}
\label{eq3.18}
\sup\limits_{t \in (r,\frac{1}{c})}t^{-1}\fint\limits_{Q(x,t)}\Bigl|f(y)-\fint\limits_{Q(x,t)}f(z)\,d\mu_{k(t)}(z)\Bigr|\,d\mu_{k(t)}(y) \le C f^{\sharp}_{S}(x',cr) \le C f^{\sharp}_{S}(x',r).
\end{equation}

Assume now that $tc > rc \geq 1$. We use Lemma \ref{Lm3.3} and estimate \eqref{eq3.6}. This gives the estimate

\begin{equation}
\begin{split}
\label{eq3.18'}
&t^{-1}\fint\limits_{Q(x,t)}\Bigl|f(y)-\fint\limits_{Q(x,t)}f(z)\,d\mu_{k(t)}(z)\Bigr|\,d\mu_{k(t)}(y) \le C \fint\limits_{Q(x,t)}|f(y)|\,d\mu_{k(t)}(y)\\
&\le C \fint\limits_{Q(x',ct)}|f(y)|\,d\mu_{k(t)}(y) \le C \fint\limits_{Q(x',ct)}|f(y)|\,d\mu_{0}(y).
\end{split}
\end{equation}

Combining \eqref{eq3.18} and \eqref{eq3.18'}, we conclude.

\begin{Lm}
\label{Lm3.6}
Let $d \in [0,n]$. Let $S$ be a~closed set with $\operatorname{dim}_{H}S \geq d$. Assume that there exists a $d$-regular system of measures $\{\mu_{k}\}_{k \in \mathbb{N}_{0}}$ on $S$. Let $f \in L_{1}^{\rm loc}(S,\mu_{k})$ for every $k \in \mathbb{N}_{0}$. Given a point $x_{0} \in S$ and a number $r_{0} \in (0,1)$, for every $r \in (r_{0},1)$
\begin{equation}
\label{eq3.19}
\frac{1}{r}\Bigl|\fint\limits_{Q(x_{0},r_{0})}f(y)\,d\mu_{k(r_{0})}(y)-\fint\limits_{Q(x_{0},r)} f(z)\,d\mu_{k(r)}(z)\Bigr| \le C  f^{\sharp}_{\{\mu_{k}\}}(x_{0},r_{0}),
\end{equation}
If in addition
\begin{equation}
\label{eq3.20}
f(x_{0})=\lim\limits_{r \to 0}\fint\limits_{Q(x_{0},r)} f(y)\,d\mu_{k(r)}(z),
\end{equation}
then for every $r \in (0,1)$
\begin{equation}
\label{eq3.21}
\frac{1}{r}\Bigl|f(x_{0})-\fint\limits_{Q(x_{0},r)} f(y)\,d\mu_{k(r)}(z)\Bigr| \le C  f^{\sharp}_{\{\mu_{k}\}}(x_{0}).
\end{equation}
The constant $C > 0$ in \eqref{eq3.19} and \eqref{eq3.21} does not depend on $x_{0}$, $r$.
\end{Lm}

\textbf{Proof}. We prove only \eqref{eq3.21} because the proof of \eqref{eq3.19} is similar. Using \eqref{eq3.20} and Lemma \ref{Lm3.4}, we obtain

\begin{equation}
\label{eq3.22}
\begin{split}
&\Bigl|f(x_{0})-\fint\limits_{Q(x_{0},r)} f(y)\,d\mu_{k(r)}(z)\Bigr| \\
&\le \sum\limits_{j=0}^{\infty}\frac{r}{2^{j}}\frac{2^{j}}{r}\Bigl|\fint\limits_{Q(x_{0},\frac{r}{2^{j}})}f(z)d\mu_{k(\frac{r}{2^{j}})}(z)-\fint\limits_{Q(x_{0},\frac{r}{2^{j+1}})}f(z')d\mu_{k(\frac{r}{2^{j+1}})}(z')\Bigr|\\
&\le C \sum\limits_{j'=0}^{\infty}\frac{r}{2^{j'}}\sup\limits_{j \in \mathbb{N}_{0}}\frac{2^{j}}{r}\fint\limits_{Q(x_{0},2^{-j}r)}\Bigl|f(y)-\fint\limits_{Q(x_{0},2^{-j}r)}f(z)\,d\mu_{k(\frac{r}{2^{j}})}(z)\Bigr|\,d\mu_{k(\frac{r}{2^{j}})}(y)\\
&\le C r f^{\sharp}_{\{\mu_{k}\}}(x_{0}).
\end{split}
\end{equation}

The lemma is proved.

\subsection{Porous sets}

\begin{Def}
\label{Def3.3}
Let $S$ be a closed nonempty subset of $\mathbb{R}^{n}$ and $\lambda \in (0,1)$. For every $j \in \mathbb{N}_{0}$ define
$$
S_{j}(\lambda):=\{x \in S | \text{ there  exists }  y \in Q(x,2^{-j}) \text{ such that }  Q(y,\lambda 2^{-j}) \subset \mathbb{R}^{n}\setminus S\}.
$$
and call $S_{j}(\lambda)$ maximal \textit{$2^{-j}$-porous subset} of $S$. We say that $S$ \textit{is porous} if there exists a number $\lambda \in (0,1)$ such that $S_{j}(\lambda)=S$ for every $j \in \mathbb{N}_{0}$.
\end{Def}

\begin{Remark}
\label{Rem3.4}
Let us note useful facts about porous subsets. Fix an arbitrary $\lambda \in (0,1)$

{\rm (1)} It is easy to see that $S_{j}(\lambda)$ is closed for every $j \in \mathbb{N}_{0}$.

{\rm (2)} The observation that Ahlfors $d$-regular sets with $d \in [0,n)$ are porous was done in \cite{Jon5}.
See also Proposition 9.18 in \cite{Tr2} which gives this fact as a special case. Let us also mention that a
set $S\subset \mathbb{R}^{n}$ is porous if,and only if, its Assouad dimension is strictly less than $n$ \cite{Luu}.
\end{Remark}

\textbf{Example 3.1} Let $\beta: [0,+\infty) \to [0,+\infty)$ be continuous strictly increasing function such that $\beta(0)=0$ and $\beta(t) > 0$, $t > 0$. Consider the 
closed single cusp

\begin{equation}
\notag
G^{\beta}:=\{x=(x',x_{n}) \in \mathbb{R}^{n}| x_{n} \in [0,\infty), \|x'\| \le \beta(x_{n})\}.
\end{equation}

It is easy to see that $\partial G^{\beta}$ is porous.

Recall Lemma \ref{Lm2.1} and Definition \ref{Def2.7}. Recall also that by $k(\varkappa)$ we denoted the unique integer number such that $r_{\varkappa} = 2^{-k(\varkappa)}$.
Recall also that we measure distances in $\mathbb{R}^{n}$ in the uniform norm.

\begin{Lm}
\label{Lm3.7}
Let $S$ be a closed nonempty set in $\mathbb{R}^{n}$. Let $Q_{\varkappa}=Q(x_{\varkappa},r_{\varkappa})$ be a Whitney cube in $W_{S}$.
Then $\widetilde{x}_{\varkappa} \in S_{j}(\lambda)$ for every $j \geq k(\varkappa)$ and $\lambda \in (0,1)$.
Furthermore, $Q(\widetilde{x}_{\varkappa},\frac{r_{\varkappa}(c-1)}{c}) \cap S \subset S_{k(\varkappa)}(\lambda)$ for every $c >1$ and every $\lambda \in (0,\frac{1}{c})$
\end{Lm}

\textbf{Proof.} Consider the interval $(x_{\varkappa},\widetilde{x}_{\varkappa}):=\{x=x_{\varkappa}+t(\widetilde{x}_{\varkappa}-x_{\varkappa})| t \in (0,1)\}$.
It clear that $S \cap (x_{\varkappa},\widetilde{x}_{\varkappa}) = \emptyset$ because otherwise there exists a point $x' \in S$ such that $\|x_{\varkappa}-x'\| < \|x_{\varkappa}-\widetilde{x}_{\varkappa}\|=\operatorname{dist}(x_{\varkappa},S)$. For every $r \in (0,r_{\varkappa}]$ consider the point $y_{r}:=(x_{\varkappa},\widetilde{x}_{\varkappa}) \cap \partial Q(\widetilde{x}_{\varkappa},r)$. From Remark \ref{Rem2.5''} $\operatorname{dist}(y_{r},S) = r$. Hence for every $\lambda \in (0,1)$ the cube $Q(y_{r},\lambda r) \subset \mathbb{R}^{n} \setminus S$. This proves the first claim of the lemma.

Given a number $c >1$ we set $r_{c}:=\frac{r_{\varkappa}}{c}$. Then from Remark \ref{Rem2.5''} we conclude that $\operatorname{dist}(y_{r_{c}},S) = \frac{r_{\varkappa}}{c}$. On the other hand, it is clear that $y_{r_{c}} \in Q(x,r_{\varkappa})$ for every $x \in Q(\widetilde{x}_{\varkappa},\frac{c-1}{c}r_{\varkappa})$. This proves the second claim of the lemma.





\begin{Lm}
\label{Lm3.8}
Let $S$ be a closed nonempty set in $\mathbb{R}^{n}$. Let $W_{S}=\{Q_{\varkappa}\}_{\varkappa \in I}$ be the Whitney decomposition of $\mathbb{R}^{n} \setminus S$. Let $\lambda \in (0,1)$ and $k \in \mathbb{N}_{0}$. Let $x \in S_{k}(\lambda)$. Then there exists a point $y(x) \in Q(x,2^{-k})$ such that
\begin{equation}
\label{eq3.23}
\frac{\lambda 2^{-k}}{5} \le \operatorname{diam}Q_{\varkappa} \le 2^{-k}.
\end{equation}
for every Whitney cube $Q_{\varkappa} \ni y(x)$.
\end{Lm}

\textbf{Proof.} By Definition \ref{Def3.3} there exists a point $y \in Q(x,2^{-k})$ such that $Q(y,\lambda 2^{-k}) \subset \mathbb{R}^{n} \setminus S$.
We set $y(x):=y$. Now we prove \eqref{eq3.23}. Consider  an arbitrary Whitney cube $Q_{\varkappa} \ni y(x)$. From \eqref{eq2.4} we have

\begin{equation}
\notag
\operatorname{diam}Q_{\varkappa} \le \operatorname{dist}(Q_{\varkappa},S) \le \operatorname{dist}(S,y(x)) \le 2^{-k}.
\end{equation}

On the other hand, using \eqref{eq2.4} again, we have

\begin{equation}
\notag
\lambda 2^{-k} \le \operatorname{dist}(y(x),S) \le \operatorname{dist}(Q_{\varkappa},S)+\operatorname{diam}(Q_{\varkappa}) \le 5\operatorname{diam}(Q_{\varkappa}).
\end{equation}

Combining the estimates above, we conclude.

\smallskip

\section{Proof of the main results}

Recall that $\|x-y\|:=\|x-y\|_{\infty}:=\max\{|x_{i}-y_{i}|: i=1,...,n\}$ for $x,y \in \mathbb{R}^{n}$. Recall also that by $\mathcal{L}_{n}$ we denote the classical
Lebesgue $n$-dimensional measure on $\mathbb{R}^{n}$.

\subsection{Pointwise characterization of Sobolev functions}

The following theorem gives a~pointwise characterization of functions in the first-order Sobolev space $W_{p}^{1}(\mathbb{R}^{n})$.
It was proven for the first time in~\cite{Haj} in a~slightly different form. This theorem will help us to estimate the Sobolev norm of the extension.

\begin{Th}
\label{Th4.1}
Assume that $p \in (1,\infty]$ and take $F \in L_{p}(\mathbb{R}^{n})$.
Then $F \in W_{p}^{1}(\mathbb{R}^{n})$ if and only if there exist a~nonnegative function
$g \in L_{p}(\mathbb{R}^{n})$ a set $E_{F}$ with $\mathcal{L}_{n}(E_{F})=0$ and positive constant $\delta$ such that

\begin{equation}
\label{eq4.1}
|F(x)-F(y)| \le  \|x-y\|\Bigl(g(x)+g(y)\Bigr)
\end{equation}
for every $x,y \in \mathbb{R}^{n} \setminus E_{F}$  with $\|x-y\| < \delta$.

Furthermore

\begin{equation}
\label{eq4.2}
\|F|L^{1}_{p}(\mathbb{R}^{n})\| \le  C \|g|L_{p}(\mathbb{R}^{n})\|,
\end{equation}

where the constant $C > 0$ does not depend on $g$.
\end{Th}

\textbf{Proof.} The proof repeats that of Theorem~1 of \cite{Haj} with minor adjustments.

\subsection{Extension operator}

Recall Lemma \ref{Lm2.1} and Definition \ref{Def2.7}.
Let $S$ be a closed nonempty set $S$ in $\mathbb{R}^{n}$. Let $W_{S}=\{Q_{\varkappa}\}_{\varkappa \in I}$ be the Whitney decomposition of the set $\mathbb{R}^{n} \setminus S$. Recall that $\mathcal{I}\subset I$ denotes the index set labeling all Whitney cubes with side length $\le 1$. Recall that for every $\varkappa \in I$ the symbol $\widetilde{x}_{\varkappa}$ denotes a metric projection of the center $x_{\varkappa}$ of the cube $Q_{\varkappa}$ to the set $S$ and $\widetilde{Q}_{\varkappa}=Q(\widetilde{x}_{\varkappa}, r_{\varkappa})$.

In what follows we will frequently use the following notation. For every $\varkappa \in I$ let $k(\varkappa)$ be the unique integer
number for which  $r_{\varkappa} \in [2^{-k(\varkappa)},2^{-k(\varkappa)+1})$.

Now we are ready to present our construction of the extension operator.

\begin{Def}
\label{Def4.1}
Let $d \in [0,n]$. Let $S$ be a closed set with $\operatorname{dim}_{H}S \geq d$. Assume that there exists $d$-regular system of measures $\{\mu_{k}\}=\{\mu_{k}\}_{k \in \mathbb{N}_{0}}$ on $S$. Assume that $f \in L^{\text{loc}}_{1}(S,\mu_{k})$ for every $k \in \mathbb{N}_{0}$. For every $\varkappa \in \mathcal{I}$ we set
$$
f_{\varkappa}:=\frac{1}{\mu_{k(\varkappa)}(\widetilde{Q}_{\varkappa} \cap S)}\int\limits_{\widetilde{Q}_{\varkappa} \cap S}f(x)\,d\mu_{k(\varkappa)}(x).
$$

With the same family of functions $\{\varphi_{\varkappa}\}_{\varkappa \in I}$ as in Lemma~\ref{Lm2.3}, put
\begin{equation}
\label{eq4.3}
F(x):=\operatorname{Ext}[f](x):=\chi_{S}(x)f(x)+\sum\limits_{\varkappa \in \mathcal{I}}\varphi_{\varkappa}(x)f_{\varkappa}, \quad x \in \mathbb{R}^{n}.
\end{equation}
\end{Def}

\begin{Remark}
\label{Rem4.1}
Actually, \eqref{eq4.3} defines not just one extension operator, but a~whole family of operators.
The reason is that the choice of a~$d$-regular system of measures $\{\mu_{k}\}$ is not unique. Furthermore, the choice of metric projections $\widetilde{x}_{\varkappa}$ is also not unique.
\end{Remark}

\subsection{Poincare-type inequalities}

The aim of this subsection is to prove Poincare-type inequalities with a $d$-regular system of measures. This inequality will be the cornerstone in proving "direct" trace theorem.

Recall the classical Poincare-type inequlity.

\begin{Lm}
\label{Lm4.1}
Assume that $F \in W_{1}^{1,\operatorname{loc}}(\mathbb{R}^{n})$. Then for every cube $Q=Q(x,r) \subset \mathbb{R}^{n}$ with $r > 0$
\begin{equation}
\label{eq4.4}
\fint\limits_{Q}\Bigl|F(y)-\fint\limits_{Q}F(z)\,d\mathcal{L}_{n}(z)\Bigr|\,d\mathcal{L}_{n}(y) \le C(n) r \fint\limits_{Q}|\nabla F(y)|\,d\mathcal{L}_{n}(y).
\end{equation}
\end{Lm}

\textbf{Proof.} Using the density of smooth functions in the space $W^{1}_{1}(\operatorname{int}Q(x,r))$ for every $r > 0$ (see section 1.4 of \cite{Maz2}) it is sufficient to prove \eqref{eq4.4}
only for smooth functions. But the latter is a well known fact, see section 8.1 in \cite{Kos}.

The following lemma is standard. We present the proof to make our exposition complete.

\begin{Lm}
\label{Lm4.2}
Let $d \in [0,n]$. Let $F \in W^{1}_{p}(\mathbb{R}^{n})$ for some $p \in (1,\infty)$, $p > n-d$. Then there exists a representative $\widehat{F}$ of
the element $F$ such that for $\mathcal{H}^{d}$-a.e. point $x \in \mathbb{R}^{n}$ and for every cube $Q(y,r) \ni x$
\begin{equation}
\label{eq4.5}
\Bigl|\widehat{F}(x)-\fint\limits_{Q(y,r)}F(z)\,d\mathcal{L}_{n}(z)\Bigr| \le C \int\limits_{Q(y,r)}\frac{|\nabla F(z)|}{|x-z|^{n-1}}d\mathcal{L}_{n}(z),
\end{equation}
where the constant $C > 0$ does not depend on $F$,$x$ and $r$.
\end{Lm}

\textbf{Proof.} We give only a sketch of the proof because all steps are routine. For every $j \in \mathbb{N}$ we write the cube $Q=Q(y,r)$ as a union of $2^{jn}$ equal cubes with disjoint interiors and choose an
arbitrary such cube $Q_{j} \ni x$ with side length $\frac{r}{2^{j}}$. We set $Q_{0}=Q(y,r)$.

Using Theorem \ref{Th2.1} it is easy to show that there exists a representative $\widehat{F}$ such that

$$
\widehat{F}(x)=\lim\limits_{j \to \infty}\fint\limits_{Q_{j}}F(z)\,d\mathcal{L}_{n}(z).
$$

 Using this and triangle inequality we clearly have

\begin{equation}
\label{eq4.6}
\begin{split}
\Bigl|\widehat{F}(x)-\fint\limits_{Q_{0}}F(z)\,d\mathcal{L}_{n}(z)\Bigr| \le \sum\limits_{j=1}^{\infty}\Bigl|\fint\limits_{Q_{j-1}}F(z)\,d\mathcal{L}_{n}(z)-\fint\limits_{Q_{j}}F(z)\,d\mathcal{L}_{n}(z)\Bigr|.
\end{split}
\end{equation}

It is clear that $|x-y| \le \frac{r}{2^{j}}$ for every $y \in Q_{j}$. Hence, using \eqref{eq4.6} and Lemma \ref{Lm4.1}, we get for $\mathcal{H}^{d}$-a.e. point $x \in \mathbb{R}^{n}$ the estimate

\begin{equation}
\label{eq4.7}
\begin{split}
&\Bigl|\widehat{F}(x)-\fint\limits_{Q(y,r)}F(z)\,d\mathcal{L}_{n}(z)\Bigr| \le C \sum\limits_{j=1}^{\infty}\frac{r}{2^{j}}\fint\limits_{Q_{j}}|\nabla F(z)|\,d\mathcal{L}_{n}(z)\\
&\le C \sum\limits_{i=0}^{j}\Bigl(\frac{2^{i}}{r}\Bigr)^{n-1}\int\limits_{Q_{j}\setminus Q_{j+1}}|\nabla F(z)|\,d\mathcal{L}_{n}(z) \le C \int\limits_{Q(y,r)}\frac{|\nabla F(z)|}{|x-z|^{n-1}}\,d\mathcal{L}_{n}(z).
\end{split}
\end{equation}

The lemma is proved.

Let $\alpha \in [0,n)$. Given a function $g \in L^{\operatorname{loc}}_{1}(\mathbb{R}^{n})$, we set

\begin{equation}
\notag
\operatorname{I}_{\alpha}[g](x):=\int\limits_{\mathbb{R}^{n}}\frac{g(y)}{\|x-y\|^{\alpha}}\,d\mathcal{L}_{n}(y).
\end{equation}

Now we formulate a particular case of the result obtained in \cite{Ver}.

\begin{Th}
\label{Th4.2}
Let $\alpha \in (0,n)$ and $q \in (1,\infty)$, and let $\mu$ be a positive Borel measure on $\mathbb{R}^{n}$. The the following statements are equivalent:

{\rm (1)} the inequality
\begin{equation}
\label{eq4.8}
\int\limits_{\mathbb{R}^{n}}\operatorname{I}_{\alpha}[g](x)\,d\mu(x) \le C \|g|L_{q}(\mathbb{R}^{n})\|
\end{equation}
holds for every $g \in L_{q}(\mathbb{R}^{n})$ with the constant $C > 0$ independent on $g$;

{\rm (2)}
\begin{equation}
\notag
\int\limits_{\mathbb{R}^{n}}\int\limits_{0}^{+\infty}\left[\frac{\mu(B(x,r))}{r^{n-\alpha q}}\right]^{q'-1}\frac{dr}{r}\,d\mu(x) < +\infty.
\end{equation}

Moreover, the least possible constant $C$ in \eqref{eq4.8} satisfies the inequality

\begin{equation}
\label{eq4.9}
a\left(\int\limits_{\mathbb{R}^{n}}\int\limits_{0}^{+\infty}\left[\frac{\mu(B(x,r))}{r^{n-\alpha q}}\right]^{q'-1}\frac{dr}{r}\,d\mu(x)\right)^{\frac{1}{q'}} \le C
\le b\left(\int\limits_{\mathbb{R}^{n}}\int\limits_{0}^{+\infty}\left[\frac{\mu(B(x,r))}{r^{n-\alpha q}}\right]^{q'-1}\frac{dr}{r}\,d\mu(x)\right)^{\frac{1}{q'}},
\end{equation}

where the constants $a,b > 0$ do not depend on $\mu$.

\end{Th}

Now we are ready to formulate the main result of this subsection.

\begin{Th}
\label{Th4.3}
Let $d \in [0,n]$. Let $S$ be a~closed set with $\operatorname{dim}_{H}S \geq d$. Assume that there exists a $d$-regular system of measures $\{\mu_{k}\}_{k \in \mathbb{N}_{0}}$ on $S$.
Let $q \in (\max\{1,n-d\},\infty)$. Take $F \in W_{q}^{1}(\mathbb{R}^{n})$.
Then for every cube $Q=Q(x,r)$ with $x \in S$ and $r \in (0,1]$
\begin{equation}
\label{eq4.10}
\fint\limits_{Q \cap S}\Bigl|F|_{S}(y)-\fint\limits_{Q}F(z)\,d\mathcal{L}_{n}(z)\Bigr|\,d\mu_{k(r)}(y) \le C r
\left(\fint\limits_{Q}\sum\limits_{|\alpha|=1}|D^{\alpha}F(t)|^{q}\,d\mathcal{L}_{n}(t)\right)^{\frac{1}{q}},
\end{equation}
where the constant $C > 0$ is independent of~$F$.
\end{Th}

\textbf{Proof.} Let us fix a cube $Q=Q(x,r_{0})$. We set $g:=\chi_{Q}|\nabla F|$. We can rewrite \eqref{eq4.5} as follows. For $\mathcal{H}^{d}$-a.e. $y \in Q \cap S$

\begin{equation}
\label{eq4.11}
\Bigl|\widehat{F}(y)-\fint\limits_{Q(x,r_{0})}F(z)\,d\mathcal{L}_{n}(z)\Bigr| \le C \operatorname{I}_{n-1}[g](y).
\end{equation}

Now we consider the measure $\mu_{Q}:=\mu_{k(r)}\lfloor{Q \cap S}$. Apply Theorem \ref{Th4.2} with the measure $\mu_{Q}$ (instead of $\mu$) and with $\alpha=1$. It is clear that $Q \cap S \subset Q(x,2r_{0})$ for every $x \in Q \cap S$. Hence, using \eqref{eq3.4} and \eqref{eq3.5} we can easily estimate the least possible constant $C$ in \eqref{eq4.8} from above.
Direct computations give

\begin{equation}
\label{eq4.12}
\begin{split}
&\int\limits_{\mathbb{R}^{n}}\int\limits_{0}^{\infty}\left[\frac{\mu_{Q}(B(y,r))}{r^{n-q}}\right]^{q'-1}\frac{dr}{r}\,d\mu_{Q}(y)\\
&\le C_{1}(r_{0})^{d}\left(C_{2}\int\limits_{0}^{2r_{0}}r^{(q+d-n)(q'-1)-1}\,dr+C_{3}(r_{0})^{d(q'-1)}\int\limits_{2r_{0}}^{\infty}\frac{dr}{r^{(n-q)(q'-1)+1}}\right) \le C_{4}(r_{0})^{1+d-\frac{n}{q}}.
\end{split}
\end{equation}

The constants $C_{1},C_{2},C_{3},C_{4} > 0$ in \eqref{eq4.12} do not depend on $x,r_{0}$.

Combining \eqref{eq4.8}, \eqref{eq4.9}, \eqref{eq4.11}, \eqref{eq4.12} and using \eqref{eq3.5}, we obtain \eqref{eq4.10}.

The theorem is proved.

\begin{Ca}
\label{Ca4.1}
Let $d \in [0,n]$. Let $S$ be a~$d$-thick closed set. Assume that $q \in (\max\{1,n-d\},\infty)$  and take $F \in W_{q}^{1}(\mathbb{R}^{n})$. Set $f:=F|_{S}$. Then
for every $r \in [0,1)$ and every $x \in S$
\begin{equation}
\label{eq4.13}
f^{\sharp}_{\{\mu_{k}\}}(x,r) \le C \Bigl(\operatorname{M}^{<1}_{>r}[|\nabla F|^{q}](x)\Bigr)^{\frac{1}{q}}
\end{equation}
\end{Ca}

\textbf{Proof.} Using Theorem \ref{Th4.3}, we have for every $x \in S$

\begin{equation}
\label{eq4.14}
\begin{split}
&\sup\limits_{t \in (r,1)}\frac{1}{t}\fint\limits_{Q(x,t) \cap S}\Bigl|f(y)-\fint\limits_{Q(x,t) \cap S}f(z)\,d\mu_{k(t)}(z)\Bigr|\,d\mu_{k(t)}(y)\\
&\le \sup\limits_{t \in (r,1)}\frac{2}{t}\fint\limits_{Q(x,t) \cap S}\Bigl|f(y)-\fint\limits_{Q(x,t)}f(z)\,d\mathcal{L}_{n}(y)\Bigr|\,d\mu_{k(t)}(y)\\
&\le C \sup\limits_{t \in (r,1)}\left(\fint\limits_{Q(x,t)}|\nabla F(y)|^{q}\,d\mathcal{L}_{n}(y)\right)^{\frac{1}{q}} \le C \Bigl(\operatorname{M}^{<1}_{>r}[|\nabla F|^{q}](x)\Bigr)^{\frac{1}{q}}.
\end{split}
\end{equation}

This proves the claim.

\subsection{Pointwise estimates of the extension}

In this section we prove the lemma which will be the cornerstone in proving "inverse" trace theorem. Recall that given $r > 0$ we denoted by $k(r)$ the unique integer number for which
$r \in [2^{-k(r)},2^{-k(r)+1})$. For every $\varkappa \in I$ we set $k(\varkappa):=k(r_{\varkappa})$.

\begin{Lm}
\label{Lm4.3}
Let $d \in [0,n]$. Let $S$ be a~closed set in~$\mathbb{R}^{n}$ with $\operatorname{dim}_{H}S \geq d$. Assume that there exists  a $d$-regular system of measures $\{\mu_{k}\}=\{\mu_{k}\}_{k \in \mathbb{N}_{0}}$  on $S$.
Let $f \in L_{1}^{\operatorname{loc}}(S,\mu_{k})$ for every $k \in \mathbb{N}_{0}$. Suppose that
\begin{equation}
\label{eq4.15}
\lim\limits_{k \to \infty}\fint\limits_{Q(x,2^{-k})}|f(x)-f(y)|\,d\mu_{k}(y)=0
\end{equation}
for $\mathcal{L}_{n}$-a.e. point $x \in S$.
Then there exist numbers $\delta \in (0,1)$ and $C > 0$ such that the function
$F:=\operatorname{Ext}[f]: \mathbb{R}^{n} \to \mathbb{R}$ defined in~\eqref{eq4.3} satisfies
\begin{equation}
\label{eq4.16}
|F(x)-F(y)| \le C \|x-y\|\Bigl(g(x)+g(y)\Bigr)
\end{equation}
for $x,y \in \mathbb{R}^{n}$ with $\|x-y\| < \delta$, where
\begin{equation}
\begin{split}
\label{eq4.17}
&g(x)=\chi_{S}(x)f^{\sharp}_{\{\mu_{k}\}}(x)\\
&+\sum\limits_{\varkappa \in \mathcal{I}}\chi_{Q_{\varkappa}}(x)\sum\limits_{\varkappa' \in b(x)}\Bigl(f^{\sharp}_{\{\mu_{k}\}}(\widetilde{x}_{\varkappa'},r_{\kappa'})+\fint\limits_{\widetilde{Q}_{\varkappa'}\cap
S}|f(z)|\,d\mu_{k(\varkappa')}(z)\Bigr), \quad x \in \mathbb{R}^{n}.
\end{split}
\end{equation}
\end{Lm}

\textbf{Proof}. Let us verify that \eqref{eq4.17} holds for all $\delta \in (0,\frac{1}{50})$.
It is obvious that we should consider five cases:

(1)
$x,y \in S$ with $\|x-y\| < \delta$;

(2)
$x \in S$ and $y \in \mathbb{R}^{n} \setminus S$ with
$\|x-y\| < \delta$;

(3)
$y \in S$ and $x \in \mathbb{R}^{n} \setminus S$
with $\|x-y\| < \delta$;

(4)
$x,y \in \mathbb{R}^{n} \setminus S$
with $\|x-y\| < \delta$ and
$x,y \in U_{\frac{1}{25}}(S)$;

(5)
$x,y \in \mathbb{R}^{n} \setminus S$
with $\|x-y\| < \delta$ and either
$x \notin U_{\frac{1}{25}}(S)$ or $y \notin U_{\frac{1}{25}}(S)$.

Clearly, in the last case $x,y \in \mathbb{R}^{n} \setminus U_{\frac{1}{50}}S$ because
$\delta \in (0,\frac{1}{50})$. In~addition, by the symmetry of the left-hand side of \eqref{eq4.16}
with respect to $x$~and~$y$, we can combine cases 2~and~3.

\textit{Case 1.} Take $f : S \to \mathbb{R}$. Assume that $x,y \in S$ and $\|x-y\| < \delta$  with $\delta \in (0,\frac{1}{50})$. Let $k$ be the unique natural number such that $\|x-y\| \in [2^{-k},2^{-k+1})$.
We have
\begin{equation}
\begin{split}
\label{eq4.18}
&|F(x)-F(y)|=|f(x)-f(y)| \\
& \le \Bigl|f(x)-\fint\limits_{Q(x,\|x-y\|)}f(z)d\mu_{k}(z)\Bigr|+\Bigl|f(y)-\fint\limits_{Q(x,\|x-y\|)}f(z)d\mu_{k}(z)\Bigr|.
\end{split}
\end{equation}
It is clear that $Q(x,\|x-y\|) \subset Q(y,2\|x-y\|)$. Hence, using \eqref{eq4.15}, Lemma \ref{Lm3.5} and Lemma \ref{Lm3.6}, for  $\mathcal{L}_{n}$-a.e. $y \in S$ we obtain
\begin{equation}
\label{eq4.19}
\begin{split}
&\Bigl|f(y)-\fint\limits_{Q(x,\|x-y\|)}f(z)d\mu_{k}(z)\Bigr| \le \Bigl|\fint\limits_{Q(y,2\|x-y\|)}f(z')d\mu_{k}(z')-\fint\limits_{Q(x,\|x-y\|)}f(z)d\mu_{k}(z)\Bigr|
\\
&+ \Bigl|f(y)-\fint\limits_{Q(y,2\|x-y\|)}f(z')d\mu_{k}(z')\Bigr| \le C\|x-y\|f^{\sharp}_{\{\mu_{k}\}}(y).
\end{split}
\end{equation}
Similar arguments for $\mathcal{L}_{n}$-almost every $x \in S$ yield
\begin{equation}
\label{eq4.20}
\begin{split}
&\Bigl|f(x)-\fint\limits_{Q(x,\|x-y\|)}f(z)d\mu_{k}(z)\Bigr| \le  C\|x-y\|f^{\sharp}_{\{\mu_{k}\}}(x).
\end{split}
\end{equation}
As a result, from \eqref{eq4.18}, \eqref{eq4.19}, \eqref{eq4.20} for $\mathcal{L}_{n}$-a.e.  $x \in S$ and $\mathcal{L}_{n}$-a.e. $y \in S$ we obtain

\begin{equation}
\label{eq4.21}
\begin{split}
|F(x)-F(y)|=|f(x)-f(y)| \le C\|x-y\|\bigl(f^{\sharp}_{\{\mu_{k}\}}(x) + f^{\sharp}_{\{\mu_{k}\}}(y)\bigr).
\end{split}
\end{equation}

\textit{Case 2.} Consider the case that $x \in S$ and $y \in U_{\delta}(S) \setminus S$.
As we noted at the begining of the proof, the case when $y \in S$ and $x \in U_{\delta}(S) \setminus S$ is similar.
Assume also that $\|x-y\| \le \delta$ with $\delta \in (0,\frac{1}{50})$.

If $\delta \in (0,1)$ then by \eqref{eq2.4} each point $y \in U_{\delta}(S)$
lies in some cube $Q_{\varkappa}$ with $r_{\varkappa} < 1$. Hence, $\varkappa \in \mathcal{I}$.
If $\delta \in (0,\frac{1}{4})$ then by \eqref{eq2.5} we may assume in addition that
\begin{equation}
\label{eq4.22}
\sum\limits_{\varkappa \in \mathcal{I}}\varphi_{\varkappa}(y) = 1, \quad \text{ for every } y \in U_{\delta}(S).
\end{equation}
Using claim~2 of Lemma~\ref{Lm2.3}, observe that $b(y):=\{\varkappa \in I: \varphi_{\varkappa}(y) \neq 0\}$.
Therefore, \eqref{eq4.3} and \eqref{eq4.22} yield
\begin{equation}
\label{eq4.23}
|F(x)-F(y)|=|f(x)-F(y)| \le \sum\limits_{\varkappa \in b(y)}\varphi_{\varkappa}(y)\Bigl|f(x)-\fint\limits_{\widetilde{Q}_{\varkappa} \cap S}f(z)\,d\mu_{k(\varkappa)}(z)\Bigr|.
\end{equation}
Fix $\varkappa \in b(y)$ and consider the cube $Q=Q(x,r)$ with
$$
r = 4\max\{\|x-\widetilde{x}_{\varkappa}\|, \operatorname{diam}\widetilde{Q}_{\varkappa}\}.
$$

Note that $\widetilde{Q}_{\varkappa}:=Q(\widetilde{x}_{\varkappa},r_{\varkappa}) \subset Q(\widetilde{x}_{\varkappa},\frac{r}{4}) \subset Q(x,r)$ and $r < 1$ for $\delta \in (0,\frac{1}{50})$.
Observe also that \eqref{eq2.4} implies the estimates

\begin{equation}
\begin{split}
\label{eq4.24}
&\|x-y\| \geq \operatorname{dist}(x,Q_{\varkappa})- \frac{1}{16}\operatorname{diam}Q_{\varkappa} \geq \frac{1}{2}\operatorname{diam}\widetilde{Q}_{\varkappa},\\
& \|x-\widetilde{x}_{\varkappa}\| \le \|x-y\|+\|y-x_{\varkappa}\|+\|x_{\varkappa}-\widetilde{x}_{\varkappa}\| \le \|x-y\| + 5\operatorname{diam}\widetilde{Q}_{\varkappa}.
\end{split}
\end{equation}

From \eqref{eq4.24} we derive with the help of elementary computations

\begin{equation}
\label{eq4.25}
r < C\|x-y\|
\end{equation}

with the constant $C > 0$ independent of $x$~and~$y$.

Hence, using \eqref{eq4.15}, \eqref{eq4.25} and Lemmas \ref{Lm3.4}, \ref{Lm3.6}, we obtain for every $\varkappa \in b(y)$ and
for $\mathcal{L}_{n}$-a.e. $x \in S$
\begin{equation}
\label{eq4.26}
\begin{split}
&\Bigl|f(x)-\fint\limits_{\widetilde{Q}_{\varkappa} \cap S}f(z)\,d\mu_{k(\varkappa)}(z)\Bigr| \le \Bigl|f(x)-\fint\limits_{Q(x,r) \cap S}f(y)\,d\mu_{k(r)}(y)\Bigr|\\
&+\Bigl|\fint\limits_{Q(\widetilde{x}_{\varkappa},\frac{r}{4}) \cap S}f(y)\,d\mu_{k(r)}(y)-\fint\limits_{Q(\widetilde{x}_{\varkappa},r_{\varkappa}) \cap S}f(z)\,d\mu_{k(\varkappa)}(z)\Bigr| \\
&+ \Bigl|\fint\limits_{Q(\widetilde{x}_{\varkappa},\frac{r}{4}) \cap S}f(y)\,d\mu_{k(r)}(y)- \fint\limits_{Q(x,r) \cap S}f(y)\,d\mu_{k(r)}(y)\Bigr| \le C \|x-y\|\Bigl(f^{\sharp}_{\{\mu_{k}\}}(x) + f^{\sharp}_{\{\mu_{k}\}}(\widetilde{x}_{\varkappa}, r_{\varkappa})\Bigr).
\end{split}
\end{equation}
As a result,
\begin{equation}
\label{eq4.27}
\begin{split}
&|F(x)-F(y)| \le C\|x-y\|\Bigl(f^{\sharp}_{\{\mu_{k}\}}(x) + \sum\limits_{\varkappa \in b(y)} f^{\sharp}_{\{\mu_{k}\}}(\widetilde{x}_{\varkappa}, r_{\varkappa})\Bigr) \\
&\le C\|x-y\|\Bigl(f^{\sharp}_{\{\mu_{k}\}}(x) + \sum\limits_{\varkappa \in \mathcal{I}}\chi_{Q_{\varkappa}}(y) f^{\sharp}_{\{\mu_{k}\}}(\widetilde{x}_{\varkappa}, r_{\varkappa})\Bigr)=C\|x-y\|\Bigl(g(x)+g(y)\Bigr).
\end{split}
\end{equation}

\textit{Case 3.} Fix $\delta \in (0,\frac{1}{50})$. Take $x,y \in \mathbb{R}^{n} \setminus S$ with $\|x-y\| \le \delta$ and $x,y \in U_{\frac{1}{25}}(S)$.

Suppose that $x \in Q_{\varkappa_{0}}$ and $y \in Q_{\varkappa_{1}}$ for some $\varkappa_{0}, \varkappa_{1} \in I$.
Observe that by \eqref{eq2.4} for $\delta \in (0,\frac{1}{25})$we have $r(Q_{\varkappa_{0}}), r(Q_{\varkappa_{1}}) < 1$.
Hence, we may assume that $\varkappa_{0}, \varkappa_{1} \in \mathcal{I}$.

There are two subcases here. In the first one there exist cubes $Q_{\varkappa_{0}} \ni x$
and $Q_{\varkappa_{1}} \ni y$ with $Q_{\varkappa_{0}} \cap Q_{\varkappa_{1}} = \emptyset$,
and in the second one $Q_{\varkappa_{0}} \cap Q_{\varkappa_{1}} \neq \emptyset$ for all cubes
$Q_{\varkappa_{0}}$ and $Q_{\varkappa_{1}}$ containing $x$~and~$y$ respectively.

Consider the first subcase. Arguing as in \eqref{eq4.23}, we see that
\begin{equation}
\label{eq4.28}
\begin{split}
|F(x)-F(y)| \le \sum\limits_{\varkappa \in b(y)} \sum\limits_ {\varkappa' \in b(x)}\Bigl|\fint\limits_{\widetilde{Q}_{\varkappa} \cap
S}f(z)\,d\mu_{k(\varkappa)}(z)-\fint\limits_{\widetilde{Q}_{\varkappa'} \cap S}f(z)\,\mu_{k(\varkappa')}(z)\Bigr|.
\end{split}
\end{equation}

For fixed $\varkappa \in b(x)$ and $\varkappa' \in b(y)$ consider the cube $Q=Q(\widetilde{x}_{\varkappa},r)$ with

$$
r:=4\max\{\|\widetilde{x}_{\varkappa}-\widetilde{x}_{\varkappa'}\|, \operatorname{diam}(\widetilde{Q}_{\varkappa}), \operatorname{diam}(\widetilde{Q}_{\varkappa'})\}.
$$

From \eqref{eq2.4} is clear that $r < 1$ for $x,y \in U_{\frac{1}{25}}(S)$ and $\|x-y\| \le \delta$ because $ \delta \in \frac{1}{50}$.

Using the condition $Q_{\varkappa_{0}} \cap Q_{\varkappa_{1}} = \emptyset$, we get

\begin{equation}
\label{eq4.29}
\|x-y\| \geq \frac{1}{2}\bigl(\operatorname{diam}(\widetilde{Q}_{\varkappa}) + \operatorname{diam}(\widetilde{Q}_{\varkappa'})\bigr).
\end{equation}

On the other hand, using \eqref{eq2.4} we have

\begin{equation}
\begin{split}
\label{eq4.30}
&\|\widetilde{x}_{\varkappa}-\widetilde{x}_{\varkappa'}\| \le \|\widetilde{x}_{\varkappa}-x_{\varkappa}\|+\|\widetilde{x}_{\varkappa'}-x_{\varkappa'}\|+\|x-y\|+\|x-x_{\varkappa}\|+\|y-x_{\varkappa'}\|\\ &\le \|x-y\|+5\operatorname{diam}(\widetilde{Q}_{\varkappa})+5\operatorname{diam}(\widetilde{Q}_{\varkappa'}).
\end{split}
\end{equation}

Combining \eqref{eq4.29} and \eqref{eq4.30}, we easily derive

\begin{equation}
\label{eq4.31}
r < C\|x-y\|
\end{equation}

with the constant $C > 0$ independent of $x$~and~$y$, as well as $\varkappa$~and~$\varkappa'$.

It is clear also that $\widetilde{Q}_{\varkappa}, \widetilde{Q}_{\varkappa'} \subset Q(\widetilde{x}_{\varkappa},\frac{r}{4})$ and $Q(\widetilde{x}_{\varkappa'},\frac{r}{4}) \subset Q:=Q(\widetilde{x}_{\varkappa},r)$.
Consequently, using this, \eqref{eq4.31} and Lemmas \ref{Lm3.4}, \ref{Lm3.6}, we obtain
\begin{equation}
\label{eq4.32}
\begin{split}
&\Bigl|\fint\limits_{\widetilde{Q}_{\varkappa} \cap S}f(z)\,d\mu_{k(\varkappa)}(z)-\fint\limits_{\widetilde{Q}_{\varkappa'} \cap S}f(z)\,d\mu_{k(\varkappa')}(z)\Bigr|
\\
&\le \Bigl|\fint\limits_{\widetilde{Q}_{\varkappa} \cap S}f(z)\,d\mu_{k(\varkappa)}(z)-\fint\limits_{Q(\widetilde{x}_{\varkappa},r) \cap S}f(z')\,d\mu_{k(r)}(z')\Bigr|\\
&+\Bigl|\fint\limits_{Q(\widetilde{x}_{\varkappa'},\frac{r}{4}) \cap S}f(z')\,d\mu_{k(\frac{r}{4})}(z')-\fint\limits_{\widetilde{Q}_{\varkappa'} \cap S}f(z)\,d\mu_{k(\varkappa')}(z)\Bigr|
\\
&+\Bigl|\fint\limits_{Q(\widetilde{x}_{\varkappa'},\frac{r}{4}) \cap S}f(z')\,d\mu_{k(\frac{r}{4})}(z')-\fint\limits_{Q(\widetilde{x}_{\varkappa},r) \cap S}f(z')\,d\mu_{k(r)}(z')\Bigr|\\
&\le C\|x-y\|\Bigl(f^{\sharp}_{\{\mu_{k}\}}(\widetilde{x}_{\varkappa}, r_{\kappa})+f^{\sharp}_{\{\mu_{k}\}}(\widetilde{x}_{\varkappa'},r_{\kappa'})\Bigr).
\end{split}
\end{equation}
From \eqref{eq4.28} and \eqref{eq4.32} we conclude that
\begin{equation}
\begin{split}
\label{eq4.33}
&|F(x)-F(y)| \le C \|x-y\|\Bigl(\sum\limits_{\varkappa \in \mathcal{I}}\chi_{Q_{\varkappa}}(x)f^{\sharp}_{\{\mu_{k}\}}(\widetilde{x}_{\varkappa},r_{\kappa}) + \sum\limits_{\varkappa' \in \mathcal{I}}\chi_{Q_{\varkappa'}}(y)f^{\sharp}_{\{\mu_{k}\}}(\widetilde{x}_{\varkappa'},r_{\kappa'})\Bigr) \\
&\le C\|x-y\|\bigl(g(x)+g(y)\bigr).
\end{split}
\end{equation}

Consider now the second subcase. Since $F \in C^{\infty}(\mathbb{R}^{n} \setminus S)$,
the mean value inequality applies. By claim~4 of Lemma~\ref{Lm2.3},
\begin{equation}
\begin{split}
\label{eq4.34}
&\frac{1}{\|x-y\|}|F(x)-F(y)| \le C(n)\max\limits_{t \in [0,1]}|\nabla F(x+t(y-x))| \\
&\le \frac{C}{r(Q_{\varkappa_{0}})} \sum\limits_{\varkappa \in b(\varkappa_{0})}\Bigl|\fint\limits_{\widetilde{Q}_{\varkappa_{0}} \cap S}f(z')\,d\mu_{k(\varkappa_{0})}(z')-
\fint\limits_{\widetilde{Q}_{\varkappa} \cap S}f(z)\,d\mu_{k(\varkappa)}(z)\Bigr|.
\end{split}
\end{equation}
Using \eqref{eq2.5} and Lemma \ref{Lm3.4}, we obtain
\begin{equation}
\begin{split}
\label{eq4.35}
&\frac{C}{r(Q_{\varkappa_{0}})} \sum\limits_{\varkappa \in b(\varkappa_{0})}\Bigl|\fint\limits_{\widetilde{Q}_{\varkappa_{0}}}f(z')\,d\mu_{k(\varkappa_{0})}(z')-
\fint\limits_{\widetilde{Q}_{\varkappa}}f(z)\,d\mu_{k(\varkappa)}(z)\Bigr| \\
&\le C \sum\limits_{\varkappa \in b(\varkappa_{0})}f^{\sharp}_{\{\mu_{k}\}}(\widetilde{x}_{\varkappa}, r_{\varkappa}) \le C\bigl(g(x)+g(y)\bigr).
\end{split}
\end{equation}
Combining \eqref{eq4.33}, \eqref{eq4.34}, and~\eqref{eq4.35}, we handle case~3.

\textit{Case 4.} Fix $\delta \in (0,\frac{1}{50})$ and $x,y \in \mathbb{R}^{n}$  such that $\|x-y\| < \delta$ and at least one of this points
does not lie in $\mathbb{R}^{n} \setminus U_{\frac{1}{25}}(S)$. Then $x,y \in \mathbb{R}^{n} \setminus U_{\frac{1}{50}}(S)$.
By~\eqref{eq2.4}, this implies that for every $Q_{\varkappa_{0}} \ni x$, $Q_{\varkappa_{1}} \ni y$
\begin{equation}
\label{eq4.36}
\begin{split}
&\frac{1}{50} \le \operatorname{dist}(x,S) \le \operatorname{diam}Q_{\varkappa_{0}}+\operatorname{dist}(Q_{\varkappa_{0}},S) \le 5\operatorname{diam}Q_{\varkappa_{0}},
\\
&\frac{1}{50} \le \operatorname{dist}(y,S) \le \operatorname{diam}Q_{\varkappa_{1}}+\operatorname{dist}(Q_{\varkappa_{1}},S) \le 5\operatorname{diam}Q_{\varkappa_{1}}.
\end{split}
\end{equation}

Consider two subcases by analogy with case~3.

In the first subcase there exist disjoint cubes $Q_{\varkappa_{0}} \ni x$ and $Q_{\varkappa_{1}} \ni y$.
Then \eqref{eq4.36} yields
\begin{equation}
\label{eq4.37}
\|x-y\| \geq \min\{\operatorname{diam}Q_{\varkappa_{0}}, \operatorname{diam}Q_{\varkappa_{1}}\} \geq \frac{1}{250}.
\end{equation}
By \eqref{eq4.3} and \eqref{eq4.17} this implies that
\begin{equation}
\begin{split}
\label{eq4.38}
&\frac{1}{\|x-y\|}|F(x)-F(y)| \\
&\le C \Bigl(\sum\limits_{\varkappa \in \mathcal{I}}\chi_{Q_{\varkappa}}(x)\int\limits_{\widetilde{Q}_{\varkappa}}|f(z)|\,\mu_{k(\varkappa)}(z) + \sum\limits_{\varkappa \in
\mathcal{I}}\chi_{Q_{\varkappa'}}(y)\fint\limits_{\widetilde{Q}_{\varkappa'}}|f(z')|\,d\mu_{k(\varkappa')}(z'))\Bigr) \\
&\le C\bigl(g(x)+g(y)\bigr).
\end{split}
\end{equation}

In the second subcase every cube $Q_{\varkappa_{0}} \ni x$ has nonempty intersection with every cube
$Q_{\varkappa_{1}} \ni y$ meet. By claim~4 of Lemma~\ref{Lm2.3} together with \eqref{eq4.3} and \eqref{eq4.17},
we obtain
\begin{equation}
\begin{split}
\label{eq4.39}
&\frac{1}{\|x-y\|}|F(x)-F(y)| \le \max\limits_{t \in [0,1]}|\nabla F(x+t(y-x))| \\
& \le C \sum\limits_{\varkappa \in b(\varkappa_{0}) \cup b(\varkappa_{1})}\fint\limits_{\widetilde{Q}_{\varkappa} \bigcap S}|f(z)|\,d\mu_{k(\varkappa)}(z) \le C \bigl(g(x)+g(y)\bigr).
\end{split}
\end{equation}
Combining \eqref{eq4.38} and \eqref{eq4.39}, we handle case~4.

The proof of Lemma~\ref{Lm4.3} is complete.
\smallskip

\subsection{Trace norm}

In this section, given a closed nonempty set $S \subset \mathbb{R}^{n}$, we introduce the functional $\mathcal{N}_{S,p,\varkappa}$ and show that this functional is bounded on the trace space $W^{1}_{p}(\mathbb{R}^{n})|_{S}$.

Recall Defintions \ref{Def3.1} and \ref{Def3.3}.

\begin{Def}
\label{Def4.2}
Let $d \in [0,n]$ and $\lambda \in (0,1)$. Let $S$ be a closed nonempty set in $\mathbb{R}^{n}$ with $\operatorname{dim}_{H}S \geq d$. Assume that there exists a $d$-regular system of measures  $\{\mu_{k}\}=\{\mu_{k}\}_{k \in \mathbb{N}_{0}}$ on $S$. Assume that $f \in L^{\operatorname{loc}}_{1}(S,\mu_{k})$ for every $k \in \mathbb{N}_{0}$.  For every $p \in (1,\infty)$ we set
\begin{equation}
\label{eq4.40}
\begin{split}
&\mathcal{SN}_{S,p}[f]:=\left(\int\limits_{S}\bigl(f^{\sharp}_{\{\mu_{k}\}}(x)\bigr)^{p}\,d\mathcal{L}_{n}(x)\right)^{\frac{1}{p}};\\
&\mathcal{BN}_{S,p,\lambda}[f]:=\left(\int\limits_{S}|f(x)|^{p}\,d\mu_{0}(x)\right)^{\frac{1}{p}}+
\left(\sum\limits_{k=1}^{\infty}\int\limits_{S_{k}(\lambda)}\bigl(f^{\sharp}_{\{\mu_{k}\}}(x,2^{-k})\bigr)^{p}\,d\mu_{k}(x)\right)^{\frac{1}{p}};\\
&\mathcal{N}_{S,p,\lambda}[f]:=\mathcal{SN}_{S,p}[f]+\mathcal{BN}_{S,p,\lambda}[f].
\end{split}
\end{equation}
\end{Def}

\begin{Remark}
\label{Rem4.2}
The functionals $\mathcal{BN}_{S,p,\lambda}$, $\mathcal{SN}_{S,p}$, and $\mathcal{N}_{S,p,\lambda}$ have values in $[0,+\infty]$. Lemma \ref{Lm3.1} implies that the values of functionals $\mathcal{SN}_{S,p}[f]$, $\mathcal{BN}_{S,p,\lambda}[f]$, $\mathcal{N}_{S,p,\lambda}[f]$ will remain the same after changing of the function $f$ on set of $\mathcal{H}^{d}$-measure zero.
Below we establish that these functionals are bounded on $W_{p}^{1}(\mathbb{R}^{n})|_{S}$.
\end{Remark}

\begin{Remark}
\label{Rem4.3}
Our notation  $\mathcal{SN}_{S,p}$ and  $\mathcal{BN}_{S,p,\lambda}$ is not picked at random. Informally speaking,
the functional $\mathcal{SN}_{p}$ is the "Sobolev part" of the norm of the function~$f$
on the set~$S$, while we may regard the functional $\mathcal{BN}_{S,p,\lambda}$ as the "Besov part" of the norm of~$f$ on~$S$. We clarify this in Example 6.1 and Example 6.2 respectively.
\end{Remark}

\begin{Lm}
\label{Lm4.4'}
Let $d \in [0,n]$, $p \in (1,\infty)$. Let $S$ be a closed set in $\mathbb{R}^{n}$ with $\operatorname{dim}_{H}S \geq d$. Assume that there exists a $d$-regular system of measures  $\{\mu_{k}\}=\{\mu_{k}\}_{k \in \mathbb{N}_{0}}$ on $S$. Assume that $f \in L^{\operatorname{loc}}_{1}(S,\mu_{k})$ for every $k \in \mathbb{N}_{0}$. Let $\{Q_{\varkappa}\}_{\varkappa \in \mathcal{I}}=\{Q(x_{\varkappa},r_{\varkappa})\}_{\varkappa \in \mathcal{I}}$ be the family of all Whitney cubes with $r_{\varkappa} \le 1$, $\varkappa \in \mathcal{I}$. Then for every $\lambda \in (0,1)$

\begin{equation}
\begin{split}
\label{eq4.41'}
&\sum\limits_{k=0}^{\infty}2^{k(d-n)}\int\limits_{S_{k}(\lambda)}\Bigl(f^{\sharp}_{\{\mu_{k}\}}(x,2^{-k})\Bigr)^{p}\,d\mu_{k}(x)\\
&\le C_{1}\sum\limits_{\varkappa \in \mathcal{I}}\mathcal{L}_{n}(Q_{\varkappa})\Bigl(f^{\sharp}_{\{\mu_{k}\}}(\widetilde{x}_{\varkappa}, r_{\kappa})\Bigr)^{p}+C_{1}\int\limits_{S}|f(y)|^{p}\,d\mu_{0}(y),
\end{split}
\end{equation}
\begin{equation}
\begin{split}
\label{eq4.42'}
&\sum\limits_{\varkappa \in \mathcal{I}}\mathcal{L}_{n}(Q_{\varkappa})\Bigl(f^{\sharp}_{\{\mu_{k}\}}(\widetilde{x}_{\varkappa}, r_{\kappa})\Bigr)^{p} \\
&\le C_{2}\sum\limits_{k=0}^{\infty}2^{k(d-n)}\int\limits_{S_{k}(\lambda)}\Bigl(f^{\sharp}_{\{\mu_{k}\}}(x,2^{-k})\Bigr)^{p}\,d\mu_{k}(x) + C_{2}\int\limits_{S}|f(x)|^{p}\,d\mu_{0}(x).
\end{split}
\end{equation}

The constants $C_{1},C_{2} > 0$ do not depend on $f$.
\end{Lm}

\textbf{Proof.} Firstly we prove the estimate \eqref{eq4.41'}. Fix $k \in \mathbb{N}_{0}$.
Let $\{x_{k,j}\}_{j \in \mathcal{J}_{k}}$ be a maximal $2^{-k}$ separated subset of $S_{k}(\lambda)$. Using item (1) of Lemma \ref{Lm2.7}, we have
\begin{equation}
\label{eq4.43'}
\int\limits_{S_{k}(\lambda)}\bigl(f^{\sharp}_{\{\mu_{k}\}}(x,2^{-k})\bigr)^{p}\,d\mu_{k}(x) \le \sum\limits_{j \in \mathcal{J}_{k}}\int\limits_{S_{k}(\lambda) \cap Q(x_{k,j},2^{-k})}\Bigl(f^{\sharp}_{\{\mu_{k}\}}(x,2^{-k})\Bigr)^{p}\,d\mu_{k}(x).
\end{equation}

Now we use Lemma \ref{Lm3.8}. Choose arbitrary point $y(x_{k,j})$ and index $\varkappa(k,j) \in \mathcal{I}$ such that $Q_{\varkappa(k,j)} \ni y(x_{k,j})$ and \eqref{eq3.23} holds. Define the function

\begin{equation}
\label{eq4.44'}
\Theta(k,j):=\varkappa(k,j), \quad \text{ for every } k \in \mathbb{N}_{0} \text{ and } j \in \mathcal{J}_{k}.
\end{equation}

It follows from \eqref{eq3.23} that if $Q_{\varkappa} \cap Q(x_{k,j},2^{-k}) \neq \emptyset$ then $Q_{\varkappa} \supset 3Q(x_{k,j},2^{-k})$. Using this, item (2) of Lemma \ref{Lm2.7} and Lemma \ref{Lm2.8} with $c=6$, it is easy to see that there exists a constant $C(n) > 0$ such that

\begin{equation}
\notag
\operatorname{card}\{\Theta^{-1}(\varkappa)\} \le C(n), \quad \varkappa \in \mathcal{I}.
\end{equation}
Hence, using this and \eqref{eq3.23}, we obtain for every $\varkappa \in \mathcal{I}$
\begin{equation}
\label{eq4.45'}
\sum\limits_{(k,j) \in \Theta^{-1}(\varkappa)}2^{-kn} \le C \mathcal{L}_{n}(Q_{\varkappa}).
\end{equation}
If $\varkappa=\Theta(k,j)$ for some $k \in \mathbb{N}_{0}$, $j \in \mathcal{J}_{k}$, then from \eqref{eq2.4} and \eqref{eq3.23} it follows that
\begin{equation}
\notag
\operatorname{dist}(x_{k,j},\widetilde{x}_{\varkappa}) \le 2^{-k}+\operatorname{dist}(Q_{\varkappa},S)+\operatorname{diam}(Q_{\varkappa}) \le 2^{-k} + 5 \operatorname{diam}(Q_{\varkappa}) \le \frac{11}{2^{k}}.
\end{equation}
This gives inclusion

\begin{equation}
\label{eq4.46'}
Q\bigl(\widetilde{x}_{\Theta(k,j)},\frac{25}{2^{k}}\bigr) \supset Q(x,2^{-k}), \quad x \in Q(x_{k,j},2^{-k}).
\end{equation}
Using \eqref{eq3.4}, \eqref{eq4.46'} and Lemma \ref{Lm3.5} with $c=25$, we obtain

\begin{equation}
\begin{split}
\label{eq4.47'}
&\int\limits_{S_{k}(\lambda) \cap Q(x_{k,j},2^{-k})}\Bigl(f^{\sharp}_{\{\mu_{k}\}}(x,2^{-k})\Bigr)^{p}\,d\mu_{k}(x) \\
&\le C 2^{-kd} \Bigl(f^{\sharp}_{\{\mu_{k}\}}\bigl(\widetilde{x}_{\Theta(k,j)},2^{-k}\bigr)+\fint\limits_{Q(\widetilde{x}_{\Theta(k,j)},25)}|f(y)|\,d\mu_{0}(y)\Bigr)^{p}.
\end{split}
\end{equation}

From \eqref{eq3.23} it is clear that $2^{-k} \geq r_{\varkappa} \geq \frac{\lambda}{5}2^{-k}$ for $\varkappa = \Theta(k,j)$. Hence, using this, monotonicity of $f^{\sharp}_{\{\mu_{k}\}}(x,t)$ with respect to $t$, and, combining the estimates \eqref{eq4.43'}, \eqref{eq4.45'}, \eqref{eq4.47'}, we derive

\begin{equation}
\begin{split}
\label{eq4.48'}
&\sum\limits_{k=0}^{\infty}2^{k(d-n)}\int\limits_{S_{k}(\lambda)}\Bigl(f^{\sharp}_{\{\mu_{k}\}}(x,2^{-k})\Bigr)^{p}\,d\mu_{k}(x)\\
&\le C \sum\limits_{k=0}^{\infty}\sum\limits_{j \in \mathcal{J}_{k}}2^{-kn}\Bigl(f^{\sharp}_{\{\mu_{k}\}}\bigl(\widetilde{x}_{\Theta(k,j)},2^{-k}\bigr)+\fint\limits_{Q(\widetilde{x}_{\Theta(k,j)},25)}|f(y)|\,d\mu_{0}(y)\Bigr)^{p} \\
& \le C \sum\limits_{\varkappa \in \mathcal{I}} \sum\limits_{(k,j) \in \Theta^{-1}(\varkappa)}2^{-kn}\Bigl(f^{\sharp}_{\{\mu_{k}\}}\bigl(\widetilde{x}_{\Theta(k,j)},2^{-k}\bigr)+\fint\limits_{Q(\widetilde{x}_{\Theta(k,j)},25)}|f(y)|\,d\mu_{0}(y)\Bigr)^{p}\\
&\le C \sum\limits_{\varkappa \in \mathcal{I}}\mathcal{L}_{n}(Q_{\varkappa})\Bigl(f^{\sharp}_{\{\mu_{k}\}}\bigl(\widetilde{x}_{\varkappa}, r_{\varkappa}\bigr)\Bigr)^{p}+C \sum\limits_{\varkappa \in \mathcal{I}}\mathcal{L}_{n}(Q_{\varkappa})\Bigl(\fint\limits_{Q(\widetilde{x}_{\varkappa},25)}|f(y)|\,d\mu_{0}(y)\Bigr)^{p}.
\end{split}
\end{equation}

Using H\"older inequality, \eqref{eq3.5} and Lemma \ref{Lm2.5} with $d\mathfrak{m}(x)=|f(x)|^{p}d\mu_{0}(x)$ and $c=25$, we have

\begin{equation}
\begin{split}
\label{eq4.48''}
&\sum\limits_{\varkappa \in \mathcal{I}}\mathcal{L}_{n}(Q_{\varkappa})\Bigl(\fint\limits_{Q(\widetilde{x}_{\varkappa},25)}|f(y)|\,d\mu_{0}(y)\Bigr)^{p} \le \sum\limits_{\varkappa \in \mathcal{I}}\mathcal{L}_{n}(Q_{\varkappa})\int\limits_{Q(\widetilde{x}_{\varkappa},25)}|f(y)|^{p}\,d\mu_{0}(y)\\
&\le C \int\limits_{S}|f(y)|^{p}\,d\mu_{0}(y).
\end{split}
\end{equation}

Combining \eqref{eq4.48'} and \eqref{eq4.48''} we obtain \eqref{eq4.41'}.

Now we prove \eqref{eq4.42'}. Using \eqref{eq3.5}, Lemma \ref{Lm3.5}, and simple inclusions $Q(\widetilde{x}_{\varkappa},3r) \supset Q(x,2r)  \supset Q(\widetilde{x}_{\varkappa},r)$ which hold for every $x \in Q(\widetilde{x}_{\varkappa},r_{\varkappa})$ and for every $r
\in [r_{\varkappa},1]$, we obtain

\begin{equation}
\begin{split}
\label{eq4.54'}
&f^{\sharp}_{\{\mu_{k}\}}(\widetilde{x}_{\varkappa},r_{\kappa}) \le C \inf\limits_{x \in Q(\widetilde{x}_{\varkappa},r_{\varkappa}) \cap S}f^{\sharp}_{\{\mu_{k}\}}(x,r_{\varkappa})+C \inf\limits_{x \in Q(\widetilde{x}_{\varkappa},r_{\varkappa}) \cap S} \fint\limits_{Q(x,2)}|f(y)|\,d\mu_{0}(y)\\
&\le C \inf\limits_{x \in Q(\widetilde{x}_{\varkappa},r_{\varkappa}) \cap S}f^{\sharp}_{\{\mu_{k}\}}(x,r_{\varkappa}) + C \int\limits_{Q(\widetilde{x}_{\varkappa},3)}|f(y)|\,d\mu_{0}(y).
\end{split}
\end{equation}
Using \eqref{doubling} with $\frac{1}{c}=1-\lambda$, \eqref{eq3.5} and H\"older inequality, we derive from \eqref{eq4.54'} (we use also identity $\mu_{k(\varkappa)}(Q(\widetilde{x}_{\varkappa},(1-\lambda)r_{\varkappa}) \cap S)=\mu_{k(\varkappa)}(Q(\widetilde{x}_{\varkappa},(1-\lambda)r_{\varkappa}))$ because $\operatorname{supp}\mu_{k} \subset S$)
\begin{equation}
\begin{split}
\label{eq4.55'}
&\mathcal{L}_{n}(Q_{\varkappa})\Bigl(f^{\sharp}_{\{\mu_{k}\}}(\widetilde{x}_{\varkappa},r_{\kappa})\Bigr)^{p} = \mathcal{L}_{n}(Q_{\varkappa})\frac{\mu_{k(\varkappa)}(Q(\widetilde{x}_{\varkappa},(1-\lambda)r_{\varkappa}))}{\mu_{k(\varkappa)}(Q(\widetilde{x}_{\varkappa},(1-\lambda)r_{\varkappa}))}
\Bigl(f^{\sharp}_{\{\mu_{k}\}}(\widetilde{x}_{\varkappa},r_{\kappa})\Bigr)^{p}\\
&\le C (r_{\varkappa})^{n-d}\int\limits_{Q(\widetilde{x}_{\varkappa},(1-\lambda)r_{\varkappa})\cap S}\Bigl(f^{\sharp}_{\{\mu_{k}\}}(x,r_{\varkappa})\Bigr)^{p}\,d\mu_{k(\varkappa)}(x)\\
&+C \mathcal{L}_{n}(Q_{\varkappa})\int\limits_{Q(\widetilde{x}_{\varkappa},3)}|f(y)|^{p}\,d\mu_{0}(y).
\end{split}
\end{equation}

From Lemma \ref{Lm2.4} it follows that the multiplicity of overlapping of the cubes $Q(\widetilde{x}_{\varkappa},r_{\varkappa})$ with the same side length is finite
and bounded above by a constant $C=C(n)$. Hence, it is easy to see that for every $j \in \mathbb{N}_{0}$
the multiplicity of overlapping of the sets $Q(\widetilde{x}_{\varkappa},(1-\lambda)r_{\varkappa}) \cap S$ with $r_{\varkappa}=2^{-j}$ is bounded above by the same constant $C(n)$.
Furthermore, from Lemma \ref{Lm3.7} it follows that $Q(\widetilde{x}_{\varkappa},(1-\lambda)r_{\varkappa}) \cap S \subset S_{j}(\lambda)$ for $r_{\varkappa}=2^{-j}$.
We combine these facts and apply for every fixed $j\in\mathbb{N}_{0}$ Lemma \ref{Lm2.9} with $d\mathfrak{m}(x)=\Bigl(f^{\sharp}_{\{\mu_{k}\}}(x,r_{\varkappa})\Bigr)^{p}\,d\mu_{k(\varkappa)}(x)$.
This gives
\begin{equation}
\begin{split}
\notag
&\sum\limits_{r_{\varkappa} = 2^{-j}}(r_{\varkappa})^{n-d}\int\limits_{Q(\widetilde{x}_{\varkappa},(1-\lambda)r_{\varkappa})\cap
S}\Bigl(f^{\sharp}_{\{\mu_{k}\}}(x,r_{\varkappa})\Bigr)^{p}\,d\mu_{k(\varkappa)}(x)\\
&\le C \sum\limits_{j \in \mathbb{N}_{0}}2^{j(d-n)}\int\limits_{S_{j}(\lambda)}\Bigl(f^{\sharp}_{\{\mu_{k}\}}(x,2^{-j})\Bigr)^{p}\,d\mu_{j}(x).
\end{split}
\end{equation}
Hence, we derive
\begin{equation}
\begin{split}
\label{eq4.56'}
&\sum\limits_{\varkappa \in \mathcal{I}} (r_{\varkappa})^{n-d}\int\limits_{Q(\widetilde{x}_{\varkappa},(1-\lambda)r_{\varkappa})\cap S}\Bigl(f^{\sharp}_{\{\mu_{k}\}}(x,r_{\varkappa})\Bigr)^{p}\,d\mu_{k(\varkappa)}(x)\\
&\le C \sum\limits_{j \in \mathbb{N}_{0}}2^{j(d-n)}\int\limits_{S_{j}(\lambda)}\Bigl(f^{\sharp}_{\{\mu_{k}\}}(x,2^{-j})\Bigr)^{p}\,d\mu_{j}(x).
\end{split}
\end{equation}

Using Lemma \ref{Lm2.5} with $d\mathfrak{m}(x)=|f(x)|^{p}d\mu_{0}(x)$ and $c=3$, we obtain

\begin{equation}
\begin{split}
\label{eq4.57'}
\sum\limits_{\varkappa \in \mathcal{I}}\mathcal{L}_{n}(Q_{\varkappa})\int\limits_{Q(\widetilde{x}_{\varkappa},3)}|f(y)|^{p}\,d\mu_{0}(y) \le C \int\limits_{S}|f(x)|^{p}\,d\mu_{0}(x).
\end{split}
\end{equation}

Combining  \eqref{eq4.55'}, \eqref{eq4.56'}, \eqref{eq4.57'}, we conclude.

Now we are ready to prove the main result of this subsection.

\begin{Th}
\label{Th4.4}
Let $d \in [0,n]$, $p \in (1,\infty)$ and $\lambda \in (0,1)$. Let $S$ be a closed set in $\mathbb{R}^{n}$ with $\operatorname{dim}_{H}S \geq d$. Assume that there exists a $d$-regular system of measures $\{\mu_{k}\}_{k \in \mathbb{N}_{0}}$ on $S$. Then the functional $\mathcal{N}_{S,p,\lambda}$ is bounded on the space $W_{p}^{1}(\mathbb{R}^{n})|_{S}$.
\end{Th}

\textbf{Proof.} Recall Definition \ref{Def2.5}. It is sufficient to verify that there exists a (universal) constant $C > 0$
such that the inequality
\begin{equation}
\label{eq4.41}
\mathcal{N}_{S,p,\lambda}[f] \le C \|F|W_{p}^{1}(\mathbb{R}^{n})\|
\end{equation}
holds for  every $F \in W_{p}^{1}(\mathbb{R}^{n})$ such that $F|_{S}=f$.

\textit{Step 1.} First of all we estimate $\mathcal{SN}_{S,p}[f]$ from above. Fix some $q \in (\max\{1,n-d\},p)$ and apply Corollary \ref{Ca4.1}. Using  (2.3), the fact that $p > q$ and Theorem \ref{Th2.2} with the exponent $\frac{p}{q}$ instead of $p$, we obtain

\begin{equation}
\label{eq4.42}
\mathcal{SN}_{S,p}[f] \le C(p,q,d,n) \|F|L_{p}^{1}(\mathbb{R}^{n})\|.
\end{equation}

\textit{Step 2.} Now we are going to estimate the second term in the right-hand side of \eqref{eq4.41'}. Fix some $q \in (\max\{1,n-d\},p)$ and apply Corollary \ref{Ca4.1}.
We also use Remark 2.3, Theorem \ref{Th2.2} with the exponent $\frac{p}{q}$ instead of $p$ and the fact that interiors of Whitney cubes are mutually disjoint.
This gives the estimate

\begin{equation}
\begin{split}
\label{eq4.45}
&\sum\limits_{\varkappa \in \mathcal{I}}\mathcal{L}_{n}(Q_{\varkappa})\Bigl(f^{\sharp}_{\{\mu_{k}\}}(\widetilde{x}_{\varkappa}, r_{\kappa})\Bigr)^{p} \le C \sum\limits_{\varkappa \in \mathcal{I}}\mathcal{L}_{n}(Q_{\varkappa})\Bigl(\operatorname{M}_{>r_{\varkappa}}[|\nabla F|^{q}](\widetilde{x}_{\varkappa})\Bigr)^{\frac{p}{q}}\\
&\le C \sum\limits_{\varkappa \in \mathcal{I}}\mathcal{L}_{n}(Q_{\varkappa})\inf\limits_{x \in Q_{\varkappa}}\Bigl(\operatorname{M}_{>r_{\varkappa}}[|\nabla F|^{q}](x)\Bigr)^{\frac{p}{q}} \le C \int\limits_{\mathbb{R}^{n}}\Bigl(\operatorname{M}_{>r_{\varkappa}}[|\nabla F|^{q}](x)\Bigr)^{\frac{p}{q}}\,d\mathcal{L}_{n}(x)\\
& \le C \|F|L_{p}^{1}(\mathbb{R}^{n})\|^{p}.
\end{split}
\end{equation}

\textit{Step 3.} It remains to estimate $\|f|L_{p}(S,\mu_{0})\|$ from above. Let $\{x_{j}\}_{j \in \mathcal{J}}$ be a maximal $1$-separated subset of $S$. Consider the family of cubes $\{Q_{j}\}_{j \in \mathcal{J}}:=\{Q(x_{j},1)\}_{j \in \mathcal{J}}$. Using item (1) of Lemma \ref{Lm2.7}, we derive the estimate

\begin{equation}
\label{eq4.47}
\begin{split}
&\int\limits_{S}|f(x)|^{p}\,d\mu_{0}(x) \\
&\le \sum\limits_{j \in \mathcal{J}}\int\limits_{Q_{j}}\Bigl|f(x)-\fint\limits_{Q_{j}}F(y)\,d\mathcal{L}_{n}(y)\Bigr|^{p}\,d\mu_{0}(x) + \sum\limits_{j \in \mathcal{J}}\int\limits_{Q_{j}}\Bigl(\fint\limits_{Q_{j}}|F(y)|\,d\mathcal{L}_{n}(y)\Bigr)^{p}\,d\mu_{0}(x).
\end{split}
\end{equation}

Using H\"older inequality, estimate \eqref{eq3.4}, item (3) of Lemma \ref{Lm2.7} and Lemma \ref{Lm2.9} with $d\mathfrak{m}(x)=|F(x)|^{p}\,d\mathcal{L}_{n}(x)$, we easily obtain

\begin{equation}
\label{eq4.48}
\sum\limits_{j \in \mathcal{J}}\int\limits_{Q_{j}}\Bigl(\fint\limits_{Q_{j}}|F(y)|\,d\mathcal{L}_{n}(y)\Bigr)^{p}\,d\mu_{0}(x) \le C \sum\limits_{j \in \mathcal{J}}\int\limits_{Q_{j}}|F(y)|^{p}\,d\mathcal{L}_{n}(y) \le C\int\limits_{\mathbb{R}^{n}}|F(y)|^{p}\,d\mathcal{L}_{n}(y).
\end{equation}

Fix $q \in (\max\{1,n-d\},p)$ and $Q_{j}=Q(x_{j},1)$. It is clear that $Q_{j} \subset Q(x,3)$ for every $x \in S \cap Q_{j}$. Recall Theorem \ref{Th2.4} and Lemma \ref{Lm3.1}. This gives for $\mu_{0}$-a.e. $x \in S \cap Q_{j}$ and for every sufficiently small $\delta \in (0,1)$

\begin{equation}
\notag
\begin{split}
&\Bigl|f(x)-\fint\limits_{Q_{j}}F(y)\,d\mathcal{L}_{n}(y)\Bigr| \le \Bigl|\fint\limits_{Q(x,1)}F(y)\,d\mathcal{L}_{n}(y)-\fint\limits_{Q_{j}}F(y')\,d\mathcal{L}_{n}(y')\Bigr|\\
&+\sum\limits_{i=1}^{\infty}\frac{2^{i\delta}}{2^{i\delta}}\Bigl|\fint\limits_{Q(x,2^{-i+1})}F(y)\,d\mathcal{L}_{n}(y)-\fint\limits_{Q(x,2^{-i})}F(y')\,d\mathcal{L}_{n}(y)\Bigr|\\
& \le C \fint\limits_{Q(x,3)}\Bigl|F(y)-\fint\limits_{Q(x,3)}F(y')\,d\mathcal{L}_{n}(y')\Bigr|\,d\mathcal{L}_{n}(y)\\
&+ C\sup\limits_{r \in (0,3)}r^{-\delta}\fint\limits_{Q(x,r)}\Bigl|F(z)-\fint\limits_{Q(x,r)}F(z')\,d\mathcal{L}_{n}(z')\Bigr|\,d\mathcal{L}_{n}(z).
\end{split}
\end{equation}

Using Lemma \ref{Lm4.1} and H\"older inequality, we continue the previous estimate and get

\begin{equation}
\label{eq4.49}
\begin{split}
&\Bigl|f(x)-\fint\limits_{Q_{j}}F(y)\,d\mathcal{L}_{n}(y)\Bigr| \\
&\le C \sup\limits_{r \in (0,3)}\frac{r}{r^{\delta}}\fint\limits_{Q(x,r)}|\nabla F(y)|\,d\mathcal{L}_{n}(y) \le  C \sup\limits_{r \in (0,3)}r^{1-\delta}\left(\fint\limits_{Q(x,r)}|\nabla F(y)|^{q}\,d\mathcal{L}_{n}(y)\right)^{\frac{1}{q}}.
\end{split}
\end{equation}

Recall that $d > n-p$. Hence we can choose $\delta$ and $q$ such that $d > n-q(1-\delta)$. Now we use item (3) of Lemma \ref{Lm2.7}, then apply Lemma \ref{Lm2.9} with $d\mathfrak{m}(x)=\Bigl(\operatorname{M}^{<3}[|\nabla F|^{q},q(1-\delta)](x)\Bigr)^{\frac{p}{q}}\,d\mu_{0}(x)$ and finally apply Theorem \ref{Th2.2} with $\alpha = q(1-\delta)$.
As a result,  we derive from \eqref{eq4.49}

\begin{equation}
\begin{split}
\label{eq4.50}
&\sum\limits_{j \in \mathcal{J}}\int\limits_{Q_{j}}\Bigl|f(x)-\fint\limits_{Q_{j}}F(y)\,d\mathcal{L}_{n}(y)\Bigr|^{p}\,d\mu_{0}(x) \le C \sum\limits_{j \in \mathcal{J}}\int\limits_{Q_{j}}\Bigl(\operatorname{M}^{<3}[|\nabla F|^{q},q(1-\delta)](x)\Bigr)^{\frac{p}{q}}d\mu_{0}(x) \\
&\le C \int\limits_{\mathbb{R}^{n}}\Bigl(\operatorname{M}^{<3}[|\nabla F|^{q},q(1-\delta)](x)\Bigr)^{\frac{p}{q}}d\mu_{0}(x) \le C \int\limits_{\mathbb{R}^{n}}|\nabla F(y)|\,d\mathcal{L}_{n}(y).
\end{split}
\end{equation}

Combining estimates \eqref{eq4.47}, \eqref{eq4.48} and \eqref{eq4.50}, we get

\begin{equation}
\label{eq4.51}
\|f|L_{p}(S,\mu_{0})\| \le C \|F|W_{p}^{1}(\mathbb{R}^{n})\|.
\end{equation}

Combining estimates \eqref{eq4.42}, \eqref{eq4.45}, \eqref{eq4.51} and \eqref{eq4.41'}, we conclude.

\subsection{Proof of the main result}

Let $g$ be the function defined in \eqref{eq4.17}. Recall Definition \ref{Def4.2}.

\begin{Lm}
\label{Lm4.5}
Let $d \in [0,n]$. Let $S$ be a closed set in $\mathbb{R}^{n}$ with $\operatorname{dim}_{H}S \geq d$. Let $\{\mu_{k}\}=\{\mu_{k}\}_{k\in\mathbb{N}_{0}}$ be a $d$-regular system of measures on $S$. Let $p \in (1,\infty)$ and $\lambda \in (0,1)$. Then
\begin{equation}
\label{eq4.52}
\|g|L_{p}(\mathbb{R}^{n})\| \le C \mathcal{N}_{S,p,\lambda}[f].
\end{equation}
\end{Lm}

\textbf{Proof.} It is sufficient to estimate $\|g|L_{p}(\mathbb{R}^{n} \setminus S)\|$. Let $\mathcal{W}_{S}=\{Q_{\varkappa}\}_{\varkappa \in \mathcal{I}}=\{Q(x_{\varkappa},r_{\varkappa})\}_{\varkappa \in \mathcal{I}}$ be the set of all Whitney cubes with $r_{\varkappa} \le 1$. Set for every $x \in \mathbb{R}^{n} \setminus S$

\begin{equation}
\begin{split}
\notag
&g_{1}(x):=\sum\limits_{\varkappa \in \mathcal{I}}\chi_{Q_{\varkappa}}(x)\sum\limits_{\varkappa' \in b(x)}f^{\sharp}_{\{\mu_{k}\}}(\widetilde{x}_{\varkappa'},r_{\kappa'}),\\
&g_{2}(x):=\sum\limits_{\varkappa \in \mathcal{I}}\chi_{Q_{\varkappa}}(x)\sum\limits_{\varkappa' \in b(x)}\fint\limits_{Q(\widetilde{x}_{\varkappa'},r_{\kappa'})}|f(y)|\,d\mu_{k(\varkappa')}(y).
\end{split}
\end{equation}
It is clear that $g(x)=g_{1}(x)+g_{2}(x)$ for every $x \in \mathbb{R}^{n} \setminus S$.

Recall that all cubes $Q_{\varkappa}$ have disjoint interiors. Furthermore, every index $\varkappa'$ belongs to sets $b(\varkappa)$ only for a finite, and independent of $\varkappa'$, number of indices $\varkappa \in I$.  Recall also that the fact $\varkappa' \in b(\varkappa)$ implies that side length of $Q_{\varkappa}$ and $Q_{\varkappa'}$ are comparable (see \eqref{eq2.5}). Hence, direct computations gives

\begin{equation}
\begin{split}
\label{eq4.53}
&\|g_{1}|L_{p}(\mathbb{R}^{n} \setminus S)\| \le C \sum\limits_{\varkappa \in \mathcal{I}}\mathcal{L}_{n}(Q_{\varkappa})\Bigl(f^{\sharp}_{\{\mu_{k}\}}(\widetilde{x}_{\varkappa},r_{\kappa})\Bigr)^{p}, \\
&\|g_{2}|L_{p}(\mathbb{R}^{n} \setminus S)\| \le C \sum\limits_{\varkappa \in
\mathcal{I}}\mathcal{L}_{n}(Q_{\varkappa})\Bigl(\fint\limits_{Q(\widetilde{x}_{\varkappa},r_{\kappa})}|f(y)|\,d\mu_{k(\varkappa)}(y)\Bigr)^{p}.
\end{split}
\end{equation}

Using Lemma \ref{Lm3.6}, for every $\varkappa \in \mathcal{I}$ we clearly have
\begin{equation}
\label{eq4.57}
\begin{split}
&\Bigl(\fint\limits_{Q(\widetilde{x}_{\varkappa},r_{\kappa})}|f(y)|\,d\mu_{k(\varkappa)}(y)\Bigr)^{p} \le
\Bigl(\fint\limits_{Q(\widetilde{x}_{\varkappa},r_{\kappa})}\Bigl|f(y)-\fint\limits_{Q(\widetilde{x}_{\varkappa},1)}f(z)\,d\mu_{0}(z)\Bigr|\,d\mu_{k(\varkappa)}(y)\Bigr)^{p} \\
&+\Bigl(\fint\limits_{Q(\widetilde{x}_{\varkappa},1)}|f(z)|\,d\mu_{0}(z)\Bigr)^{p} \le C\Bigl(f^{\sharp}_{\{\mu_{k}\}}(\widetilde{x}_{\varkappa},r_{\varkappa})\Bigr)^{p}+C\Bigl(\fint\limits_{Q(\widetilde{x}_{\varkappa},1)}|f(z)|\,d\mu_{0}(z)\Bigr)^{p}.
\end{split}
\end{equation}

Now we use \eqref{eq3.5}, H\"older inequality and Lemma \ref{Lm2.5} with $d\mathfrak{m}(x)=|f(x)|^{p}d\mu_{0}(x)$. This gives

\begin{equation}
\label{eq4.58}
\begin{split}
&\sum\limits_{\varkappa \in \mathcal{I}}\mathcal{L}_{n}(Q(\widetilde{x}_{\varkappa},r_{\varkappa}))\Bigl(\fint\limits_{Q(\widetilde{x}_{\varkappa},1)}|f(z)|\,d\mu_{0}(z)\Bigr)^{p}\\ &\le C \sum\limits_{\varkappa \in \mathcal{I}}\mathcal{L}_{n}(Q(\widetilde{x}_{\varkappa},r_{\varkappa}))\int\limits_{Q(\widetilde{x}_{\varkappa},1)}|f(x)|^{p}\,d\mu_{0}(x) \le C\int\limits_{S}|f(x)|^{p}\,d\mu_{0}(x).
\end{split}
\end{equation}

Combining estimates \eqref{eq4.42'}, \eqref{eq4.53}, \eqref{eq4.57} and \eqref{eq4.58} we conclude the proof.

\begin{Lm}
\label{Lm4.5'}
Let $d \in [0,n]$, $p \in (1,\infty)$ and $\lambda \in (0,1)$. Let $S$ be a closed set in $\mathbb{R}^{n}$ with $\operatorname{dim}_{H}S \geq d$. Let $\{\mu_{k}\}=\{\mu_{k}\}_{k\in\mathbb{N}_{0}}$ be a $d$-regular system of measures on $S$. Let $f:S \to \mathbb{R}$ be a Borel function such that $\mathcal{N}_{S,p,\lambda}[f] < \infty$. Let $F$ be the function constructed in \eqref{eq4.3}. Then
\begin{equation}
\label{norma}
\|F|L_{p}(\mathbb{R}^{n})\| \le C \mathcal{N}_{S,p,\lambda}[f].
\end{equation}
The constant $C > 0$ in \eqref{norma} does not depend on $f$.
\end{Lm}

\textbf{Proof.} It is clear that $|F(x)| \le \chi_{S}(x)f(x)+g_{2}(x)$, where $g_{2}$ is the same as in the proof of Lemma \ref{Lm4.5}. Thus, in order to establish \eqref{norma} it remains to show that
\begin{equation}
\notag
\|f|L_{p}(S,\mathcal{L}_{n})\| \le C \|f|L_{p}(S,\mu_{0})\|,
\end{equation}
where the constant $C > 0$ does not depend on $f$.
But this estimate clearly follows from \eqref{eq3.5} and Corollary \ref{Ca3.2}.

Now we can state the main result of this section. Combining this result with Remark \ref{Rem1.1}, we clearly obtain Theorem \ref{Th1.1}.

Recall that for every $r \in (0,\infty)$ we denoted by $k(r)$ the unique integer number such that $r \in [2^{-k(r)},2^{-k(r)+1})$.

\begin{Th}
\label{Th4.5}
Let $p \in (1,\infty)$, $\lambda \in (0,1)$, $d \in [0,n]$ and $d > n-p$. Let $S$ be a closed subset of $\mathbb{R}^{n}$ with $\operatorname{dim}_{H}S \geq d$. Assume that there
exists a $d$-regular system of measures $\{\mu_{k}\}_{k \in \mathbb{N}_{0}}$ on $S$. Then $f$ belongs to the $d$-trace space $W_{p}^{1}(\mathbb{R}^{n})|_{S}$
if and only if there exists a set $S' \subset S$ with $\mathcal{H}^{d}(S \setminus S')=0$ such that
\begin{equation}
\label{eq4.59}
\lim\limits_{r \to 0}\int\limits_{Q(x,r) \cap S}|f(x)-f(z)|\,d\mu_{k(r)}(z)=0, \quad \forall x \in S',
\end{equation}
and
\begin{equation}
\notag
\mathcal{N}_{S,p,\lambda}[f]< \infty.
\end{equation}
Furthermore,
\begin{equation}
\label{eq4.60}
\mathcal{N}_{S,p,\lambda}[f] \sim \|f|W_{p}^{1}(\mathbb{R}^{n})|_{S}\|,
\end{equation}
the operator $\operatorname{Ext}$ constructed in \eqref{eq4.3}
is a~bounded linear operator from $W_{p}^{1}(\mathbb{R}^{n})|_{S}$
to $W_{p}^{1}(\mathbb{R}^{n})$.

The constants of equivalence in \eqref{eq4.60} and the operator norm of\/ $\operatorname{Ext}$ depend only on parameters $p$, $n$, $\lambda$, $d$ and the constants $C_{1}, C_{2}$ in \eqref{eq3.4}, \eqref{eq3.5}.
\end{Th}

\textbf{Proof}. Split the proof in two parts.

\textit{Necessity}. Take $F \in W_{p}^{1}(\mathbb{R}^{n})$. Using Theorem \ref{Th2.4} we choose a representative $\widehat{F}$ which has Lebesgue points $\mathcal{H}^{d}$-almost everywhere in $\mathbb{R}^{n}$. Set $f:=F|_{S}$. Let $S' \subset S$ be the intersection of $S$ with the set of Lebesgue points of the function $\widehat{F}$. It is clear that $\mathcal{H}^{d}(S \setminus S')=0$.

For every $x \in S'$ we have
\begin{equation}
\begin{split}
\label{eq4.61}
&\fint\limits_{Q(x,r) \cap S}|f(x)-f(z)|\,d\mu_{k(r)}(z)\\
&\le \Bigl|f(x)-\fint\limits_{Q(x,r)}F(y)\,d\mathcal{L}_{n}(y)\Bigr|+ \int\limits_{Q(x,r) \cap S}\Bigl|f(z)-\fint\limits_{Q(x,r)}F(y)\,d\mathcal{L}_{n}(y)\Bigr|\,d\mu_{k(r)}(z).
\end{split}
\end{equation}
The first term in the right-hand side vanishes as $r \to 0$ by the definition of~$S'$.

Choose $q \in (\max\{1,n-d\},p)$ and apply Theorem \ref{Th4.3} to estimate the second term in the right-hand side of \eqref{eq4.61}. Combining this with Theorem \ref{Th2.3} we see that the second term in the right-hand side of \eqref{eq4.61} tends to zero when $r \to 0$ for $\mathcal{H}^{d}$-a.e. point $x \in S$.

Theorem \ref{Th4.4} in combination with \eqref{eq4.61} complete the proof of the necessity.

\textit{Sufficiency}. Assume that $\mathcal{N}_{S,p,\lambda}[f] < \infty$.
Then it is obvious that $f_{\varkappa} < \infty$ for all $\varkappa \in \mathcal{I}$.
Consequently, \eqref{eq4.3} yields a well-defined function $F:=\operatorname{Ext}[f] \in C^{\infty}(\mathbb{R}^{n} \setminus S)$
whose pointwise restriction to~$S$ coincides with the original function~$f$.

Show that $\mathcal{H}^{d}$-almost every point $x \in S'$  is a~Lebesgue point of the function~$F$.
Fix a~cube $Q(x,r)$ with $r \in (0,\frac{1}{100})$.
Using \eqref{eq2.5}, we see that $Q^{\ast}_{\varkappa} \cap Q(x,r) \neq \emptyset$ implies $\varkappa \in \mathcal{I}$.
Thus, using item (3) of Lemma \ref{Lm2.3}, it is easy to derive from \eqref{eq4.3} the estimate
\begin{equation}
\begin{split}
\label{eq4.63}
&\fint\limits_{Q(x,r)}|F(x)-F(z)|\,d\mathcal{L}_{n}(z) \le \frac{1}{\mathcal{L}_{n}(Q(x,r))}\int\limits_{Q(x,r) \cap S}|f(x)-f(z)|\,d\mathcal{L}_{n}(z) 
\\
&+\frac{C}{\mathcal{L}_{n}(Q(x,r))}\sum\limits_{\substack{\varkappa \in \mathcal{I} \\ Q^{\ast}_{\varkappa} \cap Q(x,r) \neq \emptyset}}\mathcal{L}_{n}(Q_{\varkappa})\fint\limits_{\widetilde{Q}_{\varkappa}\cap
S}|f(x)-f(z)|\,d\mu_{k(\varkappa)}(z).
\end{split}
\end{equation}
In \eqref{eq4.63} we also used the fact that every cube $Q^{\ast}_{\varkappa}$, $\varkappa \in I$ has nonempty intersection with at most $C=C(n)$ cubes $Q^{\ast}_{\varkappa'}$, $\varkappa' \in I$.

Using Corollary \ref{Ca3.2} and \eqref{eq4.59}, we have for $\mathcal{H}^{d}$-a.e. point $x \in S$
\begin{equation}
\begin{split}
\label{eq4.64}
&\frac{1}{\mathcal{L}_{n}(Q(x,r))}\int\limits_{Q(x,r) \cap S}|f(x)-f(z)|\,d\mathcal{L}_{n}(z) \\
&\le C \fint\limits_{Q(x,r)\cap S}|f(x)-f(z)|\,\mu_{k(r)}(z) \to 0, \quad r \to 0.
\end{split}
\end{equation}
Observe that \eqref{eq2.4} implies the inclusions
\begin{equation}
\label{eq4.71'''}
\begin{split}
&Q_{\varkappa}=Q(x_{\varkappa},r_{\varkappa}) \subset Q(x_{\varkappa},r) \subset Q(x,5r),\\
&\widetilde{Q}_{\varkappa}=Q(\widetilde{x}_{\varkappa},r_{\varkappa}) \subset Q(x,10r) \subset Q(\widetilde{x}_{\varkappa},15r)
\end{split}
\end{equation}
for all $\varkappa \in \mathcal{I}$ with $Q_{\varkappa} \cap Q(x,r) \neq \emptyset$.
Hence, we have
\begin{equation}
\label{eq4.65}
\sum\limits_{\substack{\varkappa \in \mathcal{I} \\ Q^{\ast}_{\varkappa} \cap Q(x,r) \neq \emptyset}}\mathcal{L}_{n}(Q_{\varkappa}) \le \mathcal{L}_{n}(Q(x,5r)) \le 5^{n} \mathcal{L}_{n}(Q(x,r)).
\end{equation}
From Lemma \ref{Lm3.4} and \eqref{eq4.71'''}  it is clear that (we also use the fact that $15r < 1$)
\begin{equation}
\begin{split}
\notag
&\fint\limits_{Q(\widetilde{x}_{\varkappa},r)\cap
S}\Bigl|f(z)-\fint\limits_{Q(x,10r)\cap S}f(z')\,d\mu_{k(10r)}(z')\Bigr|\,d\mu_{k(r)}(z)\\
&\le C \fint\limits_{Q(\widetilde{x}_{\varkappa},15r)\cap
S}\Bigl|f(z)-\fint\limits_{Q(\widetilde{x}_{\varkappa},15r)\cap S}f(z')\,d\mu_{k(15r)}(z')\Bigr|\,d\mu_{k(15r)}(z).
\end{split}
\end{equation}
Thus, using monotonicity of $f^{\sharp}_{\{\mu_{k}\}}(\cdot,r)$ with respect to $r$ and Lemma \ref{Lm3.6}, we have
\begin{equation}
\label{eq4.73'''}
\begin{split}
&\frac{1}{r}\fint\limits_{Q(\widetilde{x}_{\varkappa},r_{\varkappa})\cap
S}\Bigl|f(z)-\fint\limits_{Q(x,10r)\cap S}f(z')\,d\mu_{k(10r)}(z')\Bigr|\,d\mu_{k(r)}(z)\\
&\le \frac{C}{r}\fint\limits_{Q(\widetilde{x}_{\varkappa},15r)\cap
S}\Bigl|f(z)-\fint\limits_{Q(\widetilde{x}_{\varkappa},15r)\cap S}f(z')\,d\mu_{k(15r)}(z')\Bigr|\,d\mu_{k(15r)}(z)\\
&+\frac{C}{r}\fint\limits_{Q(\widetilde{x}_{\varkappa},r_{\varkappa})\cap
S}\Bigl|f(z)-\fint\limits_{Q(\widetilde{x}_{\varkappa},r) \cap S}f(z')\,d\mu_{k(r)}(z')\Bigr|\,d\mu_{k(\varkappa)}(z)\\
&\le C \Bigl(f^{\sharp}_{\{\mu_{k}\}}(\widetilde{x}_{\varkappa},15r)+f^{\sharp}_{\{\mu_{k}\}}(\widetilde{x}_{\varkappa},r_{\varkappa})\Bigr) \le Cf^{\sharp}_{\{\mu_{k}\}}(\widetilde{x}_{\varkappa},r_{\varkappa}).\\
\end{split}
\end{equation}
Using \eqref{eq4.65} and \eqref{eq4.73'''}, we obtain
\begin{equation}
\label{eq4.66}
\begin{split}
&\frac{1}{\mathcal{L}_{n}(Q(x,r))}\sum\limits_{\substack{\varkappa \in \mathcal{I} \\ Q^{\ast}_{\varkappa} \cap Q(x,r) \neq \emptyset}}\mathcal{L}_{n}(Q_{\varkappa})\fint\limits_{\widetilde{Q}_{\varkappa}\cap
S}|f(x)-f(z)|\,d\mu_{k(\varkappa)}(z)\\
& \le  \frac{1}{\mathcal{L}_{n}(Q(x,r))}\sum\limits_{\substack{\varkappa \in \mathcal{I} \\ Q^{\ast}_{\varkappa} \cap Q(x,r) \neq \emptyset}}\mathcal{L}_{n}(Q_{\varkappa}) \Bigl|f(x)-\fint\limits_{Q(x,10r)}f(z')\,d\mu_{k(10r)}(z')\Bigr|\\
&+\frac{1}{\mathcal{L}_{n}(Q(x,r))}\sum\limits_{\substack{\varkappa \in \mathcal{I} \\ Q^{\ast}_{\varkappa} \cap Q(x,r) \neq \emptyset}}\mathcal{L}_{n}(Q_{\varkappa}) \fint\limits_{\widetilde{Q}_{\varkappa}\cap
S}\Bigl|f(z)-\fint\limits_{Q(x,10r)}f(z')\,d\mu_{k(10r)}(z')\Bigr|\,d\mu_{k(\varkappa)}(z)\\
&\le C \fint\limits_{Q(x,10r) \cap S}|f(x)-f(z')|\,d\mu_{k(10r)}(z')\\
&+C\frac{r}{\mathcal{L}_{n}(Q(x,r))}\sum\limits_{\substack{\varkappa \in \mathcal{I} \\ Q^{\ast}_{\varkappa} \cap Q(x,r) \neq
\emptyset}}\mathcal{L}_{n}(Q_{\varkappa})f^{\sharp}_{\{\mu_{k}\}}(\widetilde{x}_{\varkappa},r_{\varkappa}).
\end{split}
\end{equation}
Using \eqref{eq4.65} once again, by H\"older's inequality for sums with exponents $p$~and~$p'$ we see that
\begin{equation}
\begin{split}
\label{eq4.67}
&\frac{r}{\mathcal{L}_{n}(Q(x,r))}\sum\limits_{\substack{\varkappa \in \mathcal{I} \\ Q^{\ast}_{\varkappa} \cap Q(x,r) \neq \emptyset}}\mathcal{L}_{n}(Q_{\varkappa})f^{\sharp}_{\{\mu_{k}\}}(\widetilde{x}_{\varkappa},r_{\varkappa})\\
&\le \left(\frac{r^{p}}{\mathcal{L}_{n}(Q(x,r))}\sum\limits_{\substack{\varkappa \in \mathcal{I} \\ Q^{\ast}_{\varkappa} \cap Q(x,r) \neq
\emptyset}}\mathcal{L}_{n}(Q_{\varkappa})\Bigl(f^{\sharp}_{\{\mu_{k}\}}(\widetilde{x}_{\varkappa},r_{\varkappa})\Bigr)^{p}\right)^{\frac{1}{p}}=:\bigl(J(x,r)\bigr)^{\frac{1}{p}}.
\end{split}
\end{equation}

Verify that there exists a set $S'' \subset S' \subset S$ such that $\mathcal{H}^{d}(S \setminus S'')=0$,
 and $J(x,r) \to 0$ as $r \to 0$ for all points $x \in S''$.

Consider the case $p \geq n$. In this case we can take $S''=S'$, because $\mathcal{BN}_{S,p,\lambda}[f] < \infty$ and \eqref{eq4.42'} holds.

Assume now that $1 < p < n$. Fix a constant $c > 0$ and consider the set

$$
S_{c}:=\{x \in S: \limsup\limits_{r \to 0} J(x,r) > c\}.
$$
Verify that $\mathcal{H}^{d}(S_{c})=0$. Assume on the contrary that $\mathcal{H}^{d}(S_{c}) > 0$. Then there are $\varepsilon > 0$, $\delta_{0} \in (0,1)$
such that for all $\delta \in (0,\delta_{0})$ we have
\begin{equation}
\label{eq4.68}
\mathcal{H}^{d}_{\delta}(S_{c}) \geq \varepsilon > 0.
\end{equation}

Fix an arbitrarily number $\delta \in (0,\delta_{0})$. For each point $x \in S_{c}$
find $\delta_{x} \in (0,\frac{\delta}{50})$ with $J(x,\delta_{x}) > c$. The family $\{10Q(x,\delta_{x})\}_{x \in S_{c}}$ of cubes covers $S_{c}$.
Using Theorem \ref{Th2.5} and Remark 2.8,
we find a~sequence $\{10Q_{k}\}=\{10Q(x_{k},\delta_{x_{k}})\}$
of disjoint cubes such that
$$
S_{c} \subset \bigcup_{k=1}^{\infty}50Q_{k}.
$$
Consequently, the definition of $\mathcal{H}^{d}_{\delta}(S_{c})$ and \eqref{eq4.68} yield
\begin{equation}
\label{eq4.69}
\sum\limits_{k=0}^{\infty}50^{d}\bigl(\operatorname{diam}Q_{k}\bigr)^{d} \geq \varepsilon > 0.
\end{equation}
Recall that $d > n-p$ and $\operatorname{diam}Q_{k} < 1$. Hence, $\bigl(\operatorname{diam}Q_{k}\bigr)^{d} < \bigl(\operatorname{diam}Q_{k}\bigr)^{n-p}$. This fact, \eqref{eq4.67} and the definition of $S_{c}$ yield
\begin{equation*}
\sum\limits_{k=0}^{\infty}\sum\limits_{\substack{\varkappa \in \mathcal{I} \\ Q^{\ast}_{\varkappa} \cap Q_{k} \neq \emptyset}}\mathcal{L}_{n}(Q_{\varkappa})\Bigl(f^{\sharp}_{S}(\widetilde{x}_{\varkappa},r_{\varkappa})\Bigr)^{p}
\geq 50^{-d}c^{p}\varepsilon.
\end{equation*}
However, if the cubes $10Q_{k}$ are disjoint then by Lemma~\ref{Lm2.2}
for each $\varkappa \in \mathcal{I}$ the cube
$Q^{\ast}_{\varkappa}$ can have a~nonempty intersection with at most one cube $Q_{k}$.
Then
\begin{equation}
\label{eq4.70}
\sum\limits_{\substack{\varkappa \in \mathcal{I} \\ \operatorname{diam}Q_{\varkappa} < 2\delta}}\mathcal{L}_{n}(Q_{\varkappa})\Bigl(f^{\sharp}_{S}(\widetilde{x}_{\varkappa},r_{\varkappa})\Bigr)^{p} \geq 50^{-d}c\varepsilon.
\end{equation}
On the other hand, taking into account the fact that $\mathcal{BN}_{S,p,\lambda}[f] < \infty$ and \eqref{eq4.42'} holds, it follows that
\begin{equation}
\label{eq4.71}
\sum\limits_{\substack{\varkappa \in \mathcal{I} \\ \operatorname{diam}Q_{\varkappa} < 2\delta}}\mathcal{L}_{n}(Q_{\varkappa})\Bigl(f^{\sharp}_{S}(\widetilde{x}_{\varkappa},r_{\varkappa})\Bigr)^{p} \to 0, \quad \delta \to 0.
\end{equation}
Clearly, \eqref{eq4.71} contradicts \eqref{eq4.70}.

Thus, we see that $\mathcal{H}^{d}(S_{c}) = 0$
for every $c > 0$. Consequently, $\mathcal{H}^{d}(\bigcup\limits_{n=1}^{\infty}S_{\frac{1}{n}})=0$.
However, if $x \in S' \setminus \bigcup\limits_{n=1}^{\infty}S_{\frac{1}{n}}$
then obviously $J(x,r) \to 0$ as $r \to 0$. Combined with \eqref{eq4.59}, \eqref{eq4.63}, \eqref{eq4.64}, and \eqref{eq4.66},
this implies that every point of the set $S' \setminus \bigcup\limits_{n=1}^{\infty}S_{\frac{1}{n}}$
is a~Lebesgue point of the function~%
$F$.

Finally, to conclude the sufficiency part of the theorem we should establish the estimate
\begin{equation}
\label{eq4.80}
\|F|W_{p}^{1}(\mathbb{R}^{n})\| \le C\mathcal{N}_{S,p,\lambda}[f],
\end{equation}
where the constant $C > 0$ does not depend on $f$.
But this estimate follows from Theorem \ref{Th4.1} and Lemmas \ref{Lm4.5}, \ref{Lm4.5'}.

The proof of Theorem~\ref{Th4.5} is complete.

\begin{Remark}
\label{Rem4.4}
Presently the authors do not know whether it is possible to obtain \eqref{eq4.59} from the condition $\mathcal{N}_{S,p,\lambda}[f] < \infty$.
This question is not as simple as it may seem at first. Indeed, from the proof of Theorem~\ref{Th4.5} it is clear that given a function $f:S \to \mathbb{R}$, condition $\mathcal{N}_{S,p,\lambda}[f] < \infty$
implies only that $F=\operatorname{Ext}[f] \in W_{p}^{1}(\mathbb{R}^{n})$. Then, we obtain from Theorem \ref{Th2.4} that for $\mathcal{H}^{d}$-a.e. point $x \in S$ there exists

$$
\widehat{F}(x):=\lim\limits_{r \to 0}\fint\limits_{Q(x,r)}F(y)\,d\mathcal{L}_{n}(y).
$$

But is not obvious why $f(x)=\widehat{F}(x)$ for $\mathcal{H}^{d}$-almost every $x \in S$?
\end{Remark}

\begin{Remark}
\label{Rem4.5}
As we noted above, while constructing the extension operator, we make a~choice of a $d$-regular system of measures.
It is remarkable, however, that both the statement and proof of Theorem \ref{Th4.5}, as well as the constants occurring in the proof,
depend only on the constants $C_{1}$, $C_{2}$ from Definition \ref{Def3.1} but independent of the concrete choice of a $d$-regular system of measures.
\end{Remark}

\section{Simplified criterion for sets with porous boundary}

In this section we are going to prove Theorem \ref{Th1.2} which is a simplified version of Theorem \ref{Th1.1} in the case of sets with porous boundary.
Recall Definition \ref{Def3.3}.

Let $S$ be a closed set in $\mathbb{R}^{n}$ with porous boundary. Given $\lambda > 0$, for every $k \in \mathbb{N}_{0}$ consider the sets

\begin{equation}
\label{eq5.1}
\begin{split}
&\partial S^{+}_{k}(\lambda):=\{x \in \partial S | \operatorname{ there }  \operatorname{ exists }  y \in Q(x,2^{-k}) \operatorname{ for}  \operatorname{ which } \quad Q(y,\lambda 2^{-k}) \subset \mathbb{R}^{n}\setminus S\}\\
&\partial S^{-}_{k}(\lambda):=\{x \in \partial S | \operatorname{ there }  \operatorname{ exists }  y' \in Q(x,2^{-k}) \operatorname{ for}  \operatorname{ which } \quad Q(y',\lambda 2^{-k}) \subset S \setminus \partial S\}.
\end{split}
\end{equation}

From Definition \ref{Def3.3} it is clear that if $\partial S$ is porous then there exists a number $\lambda > 0$ such that

\begin{equation}
\label{eq5.2}
\partial S= \partial S^{+}_{k}(\lambda)\cup \partial S^{-}_{k}(\lambda)
\end{equation}

for every $k \in \mathbb{N}_{0}$.

\begin{Def}
\label{Def5.1}
Let $S$ be a closed set in $\mathbb{R}^{n}$. Let $\mathfrak{m}$ be an arbitrary Radon measure with $\operatorname{supp}\mathfrak{m} = S$. Let $Q=Q(x,r)$ be a cube with $x \in S$ and $r > 0$. Given a function $f \in L_{1}^{\operatorname{loc}}(S,\mathfrak{m})$ define
\textit{the normalized with respect to the measure $\mathfrak{m}$ best approximation of $f$ by constants  on $Q$}
\begin{equation}
\label{eq5.3}
\mathcal{E}_{\mathfrak{m}}(f,Q):=\inf\limits_{c \in \mathbb{R}}\fint\limits_{Q \cap S}|f(y)-c|\,d\mathfrak{m}(y).
\end{equation}
\end{Def}


\begin{Remark}
\label{Rem5.1}
It is easy to see that
\begin{equation}
\label{eq5.4}
\mathcal{E}_{\mathfrak{m}}(f,Q) \le  \widetilde{\mathcal{E}}_{\mathfrak{m}}(f,Q) \le 2\mathcal{E}_{\mathfrak{m}}(f,Q),
\end{equation}

where

$$
\widetilde{\mathcal{E}}_{\mathfrak{m}}(f,Q):=\fint\limits_{Q \cap S}\Bigl|f(y)-\fint\limits_{Q \cap S}f(z)\,d\mathfrak{m}(z)\Bigr|\,d\mathfrak{m}(y)
$$

The exact value of $\mathcal{E}_{\mathfrak{m}}(f,Q)$ will not be important for us in the sequel. Hence, in practice we can work without loss of generality with
$\widetilde{\mathcal{E}}_{\mathfrak{m}}(f,Q)$ instead of $\mathcal{E}_{\mathfrak{m}}(f,Q)$ which is easier to compute.
\end{Remark}

\begin{Def}
\label{Def5.2}
Let $S$ be an arbitrary closed nonempty subset of $\mathbb{R}^{n}$. For every $k \in \mathbb{N}_{0}$ consider the set

\begin{equation}
\label{eq5.5}
\Sigma_{k}:=\Sigma_{k}(S):=\{x \in S| \operatorname{dist}(x,\partial S) \le 2^{-k}\}.
\end{equation}
\end{Def}

\begin{Remark}
\label{Rem5.2}
It is clear that for every $k \in \mathbb{N}_{0}$ the set $\Sigma_{k}$ is closed.
\end{Remark}


\begin{Lm}
\label{Lm5.1}
Let $d \in [0,n]$, $p \in (1,\infty)$ and $p > n-d$. Let $S$ be a closed set in $\mathbb{R}^{n}$ with $\operatorname{dim}_{H}S \geq d$. Assume that there exists a $d$-regular system of measures $\{\mu_{k}\}_{k \in \mathbb{N}_{0}}$ on $S$. Assume that $\partial S$ is porous. Then for every $F \in W_{p}^{1}(\mathbb{R}^{n})$
\begin{equation}
\label{eq5.6}
\sum\limits_{k=0}^{\infty}2^{kp(1-\frac{n-d}{p})}\int\limits_{\Sigma_{k}}\Bigl(\mathcal{E}_{\mu_{k}}(F|_{S},Q(x,2^{-k}))\Bigr)^{p}\,d\mu_{k}(x) \le C \|F|W_{p}^{1}(\mathbb{R}^{n})\|.
\end{equation}
The constant $C > 0$ in \eqref{eq5.6} does not depend on $F$.
\end{Lm}

\textbf{Proof.} Fix $k \in \mathbb{N}_{0}$ and consider an arbitrary maximal $2^{-k}$ separated subset $\{x_{k,j}\}_{j \in \mathcal{J}_{k}}$  of $\Sigma_{k}$. We set $Q_{k,j}:=Q(x_{k,j},2^{-k})$ for every $j \in \mathcal{J}_{k}$.

It is clear from the construction that for every $x \in \Sigma_{k}$ there exists index $j \in \mathcal{J}_{k}$ such that $3Q_{k,j} \supset Q(x,2^{-k})$. Using this and \eqref{eq3.5}, \eqref{eq3.6}, we get

\begin{equation}
\label{eq5.7}
\widetilde{\mathcal{E}}_{\mu_{k}}(F|_{S},Q(x,2^{-k})) \le C \widetilde{\mathcal{E}}_{\mu_{k}}(F|_{S},3Q_{k,j})
\end{equation}

Using item (1) of Lemma \ref{Lm2.7}, \eqref{eq3.4}, \eqref{eq5.7} and Remark \ref{Rem5.1}, we obtain

\begin{equation}
\label{eq5.8}
\begin{split}
&\int\limits_{S}\Bigl(\mathcal{E}_{\mu_{k}}(F|_{S},Q(x,2^{-k}))\Bigr)^{p}\,d\mu_{k}(x) \le \sum\limits_{j \in \mathcal{J}_{k}}\int\limits_{Q_{k,j}}\Bigl(\mathcal{E}_{\mu_{k}}(F|_{S},Q(x,2^{-k}))\Bigr)^{p}\,d\mu_{k}(x) \\
&\le C 2^{-kd}\sum\limits_{j \in \mathcal{J}_{k}}\Bigl(\mathcal{E}_{\mu_{k}}(F|_{S},3Q_{k,j}\Bigr)^{p}.
\end{split}
\end{equation}

Fix some $q \in (\max\{1,n-d\},p)$. Using \eqref{eq4.10} and \eqref{eq5.8},  we derive the following estimate

\begin{equation}
\label{eq5.9}
\begin{split}
&\sum\limits_{k=0}^{\infty}2^{kp(1-\frac{n-d}{p})}\int\limits_{S}\Bigl(\mathcal{E}_{\mu_{k}}(F|_{S},Q(x,2^{-k}))\Bigr)^{p}\,d\mu_{k}(x)\\
&\le C \sum\limits_{k=0}^{\infty}2^{kp(1-\frac{n}{p})}\sum\limits_{j \in \mathcal{J}_{k}}\Bigl(\mathcal{E}_{\mu_{k}}(F|_{S},3Q_{k,j})\Bigr)^{p}\\
&\le C \sum\limits_{k=0}^{\infty} \sum\limits_{j \in \mathcal{J}_{k}} \mathcal{L}_{n} (Q_{k,j}) \Bigl(\fint\limits_{3Q_{k,j}}|\nabla F(y)|^{q}\,d\mathcal{L}_{n}(y)\Bigr)^{\frac{p}{q}}.
\end{split}
\end{equation}

Fix some $\lambda > 0$ such that \eqref{eq5.2} holds. Let $\mathcal{J}^{1}_{k}$ be the set of all $j \in \mathcal{J}_{k}$ for each of which $Q_{k,j} \cap \partial S^{+}_{k}(\lambda) \neq \emptyset$. Let $\mathcal{J}^{2}_{k}$ be the set of all $j \in \mathcal{J}_{k}$ for each of which $Q_{k,j} \cap \partial S^{-}_{k}(\lambda) \neq \emptyset$. It is clear that $\mathcal{J}_{k}=\mathcal{J}^{1}_{k} \cup \mathcal{J}^{2}_{k}$ for every $k \in \mathbb{N}_{0}$.

Fix a Whitney decomposition $W^{1}$ of the set $\mathbb{R}^{n} \setminus S$ and a Whitney decomposition $W^{2}$ of the set $\operatorname{int} S$. Let $\mathcal{I}^{1}$ and $\mathcal{I}^{2}$ be the sets of all indices corresponding to cubes with side length $\le 1$ from the family $W^{1}$ and $W^{2}$ respectively.

It is clear that for every $j \in \mathcal{J}^{1}_{k}$ there is a point $x'_{k,j} \in Q_{k,j} \cap \partial S^{+}_{k}(\lambda)$. Similarly for every $j \in \mathcal{J}^{2}_{k}$ there is a point $x'_{k,j} \in Q_{k,j} \cap \partial S^{-}_{k}(\lambda)$. Using Lemma \ref{Lm3.8}, we find for every $j \in \mathcal{J}^{1}_{k}$ a point $y(x'_{k,j}) \in \mathbb{R}^{n} \setminus S$ such that for every cube $Q_{\varkappa} \ni y(x'_{k,j})$, $\varkappa \in \mathcal{I}^{1}$
\begin{equation}
\label{eq5.10}
\frac{\lambda}{5}2^{-k} \le 2r_{\varkappa} \le 2^{-k}.
\end{equation}
Similarly, for every  $j \in \mathcal{J}^{2}_{k}$ we  find a point $z(x'_{k,j}) \in \operatorname{int} S$ such that for every cube $Q_{\varkappa'} \ni z(x'_{k,j})$, $\varkappa' \in \mathcal{I}^{2}$
\begin{equation}
\label{eq5.11}
\frac{\lambda}{5}2^{-k} \le 2r_{\varkappa'} \le 2^{-k}.
\end{equation}

Consider the map $\Theta^{1}$ which takes a pair $(k,j)$ with $j \in \mathcal{J}^{1}_{k}$ and gives back an arbitrary chosen $\varkappa=\Theta^{1}(k,j) \in \mathcal{I}^{1}$
such that \eqref{eq5.10} holds. Similarly, we built the map $\Theta^{2}$ which takes a pair $(k,j)$ with $j \in \mathcal{J}^{2}_{k}$ and gives back an arbitrary chosen $\varkappa=\Theta^{2}(k,j) \in \mathcal{I}^{2}$ such that \eqref{eq5.11} holds.

Using item (3) of Lemma \ref{Lm2.7}, it is
easy to conclude that for every $\varkappa \in \mathcal{I}^{1}$ Whitney cube $Q_{\varkappa}$ has nonempty intersection with at most $C(n)$ cubes from the family
$\{Q_{k,j}\}_{j \in \mathcal{J}_{k}}$. The similar fact clearly holds for every $\varkappa' \in \mathcal{I}^{2}$. Then, there exists a constant $C(n)$ such that
\begin{equation}
\label{eq5.12}
\operatorname{card}\{(\Theta^{1})^{-1}(\varkappa)\} \le C(n), \quad \operatorname{card}\{(\Theta^{2})^{-1}(\varkappa)\} \le C(n).
\end{equation}

Let $k \in \mathbb{N}_{0}$ and $j \in \mathcal{J}^{1}_{k}$. Let $\varkappa  = \Theta^{1}(k,j)$. Then from \eqref{eq5.10} it follows that $7Q_{k,j} \supset Q_{\varkappa}$. Hence,
from \eqref{eq2.4'} we derive
\begin{equation}
\begin{split}
\label{eq5.13}
&\fint\limits_{3Q_{k,j}}|\nabla F(y)| \,d\mathcal{L}_{n}(y) \le C \fint\limits_{7Q_{k,j}}|\nabla F(y)| \,d\mathcal{L}_{n}(y) \\
&\le C \inf\limits_{x \in Q_{\varkappa}}\operatorname{M}_{> \frac{7}{2^{k}}}[|\nabla F|](x) \le C \inf\limits_{x \in Q_{\varkappa}}\operatorname{M}_{> 2^{-k}}[|\nabla F|](x).
\end{split}
\end{equation}
Similarly,  if $k \in \mathbb{N}_{0}$,   $j \in \mathcal{J}^{2}_{k}$ and $\varkappa'  = \Theta^{2}(k,j)$, then
\begin{equation}
\label{eq5.14}
\fint\limits_{3Q_{k,j}}|\nabla F(y)| \,d\mathcal{L}_{n}(y) \le C \inf\limits_{x \in Q_{\varkappa'}}\operatorname{M}_{> 2^{-k}}[|\nabla F|](x).
\end{equation}

Combining \eqref{eq5.10}, \eqref{eq5.12} and \eqref{eq5.13}, we have
\begin{equation}
\label{eq5.15}
\begin{split}
&\sum\limits_{k \in \mathbb{N}_{0}} \sum\limits_{j \in \mathcal{J}^{1}_{k}}\mathcal{L}_{n} (Q_{k,j}) \Bigl(\fint\limits_{3Q_{k,j}}|\nabla F(y)|^{q}\,d\mathcal{L}_{n}(y)\Bigr)^{\frac{p}{q}}\\
&\le C \sum\limits_{\varkappa \in \mathcal{I}^{1}}\sum\limits_{(k,j) \in (\Theta^{1})^{-1}(\varkappa)}\mathcal{L}_{n}(Q_{k,j})\inf\limits_{x \in Q_{\varkappa}}\Bigl(\operatorname{M}_{> 2^{-k}}[|\nabla F|^{q}](x)\Bigr)^{\frac{p}{q}} \\
&\le C \sum\limits_{\varkappa \in \mathcal{I}^{1}}\mathcal{L}_{n}(Q_{\varkappa})\inf\limits_{x \in Q_{\varkappa}}\Bigl(\operatorname{M}_{> 2^{-k}}[|\nabla F|^{q}](x)\Bigr)^{\frac{p}{q}} \le C \int\limits_{\mathbb{R}^{n} \setminus S}\Bigl(\operatorname{M}_{> 2^{-k}}[|\nabla F|^{q}](x)\Bigr)^{\frac{p}{q}}\,d\mathcal{L}_{n}(x).
\end{split}
\end{equation}
Similarly, from \eqref{eq5.11}, \eqref{eq5.12} and \eqref{eq5.14} we obtain
\begin{equation}
\label{eq5.16}
\begin{split}
&\sum\limits_{k \in \mathbb{N}_{0}} \sum\limits_{j \in \mathcal{J}^{2}_{k}}\mathcal{L}_{n} (Q_{k,j}) \Bigl(\fint\limits_{3Q_{k,j}}|\nabla F(y)|^{q}\,d\mathcal{L}_{n}(y)\Bigr)^{\frac{p}{q}} \le C \int\limits_{\operatorname{int} S}\Bigl(\operatorname{M}_{> 2^{-k}}[|\nabla F|^{q}](x)\Bigr)^{\frac{p}{q}}\,d\mathcal{L}_{n}(x).
\end{split}
\end{equation}

Using Theorem \ref{Th2.1}, and combining estimates \eqref{eq5.9}, \eqref{eq5.15}, \eqref{eq5.16}, we derive

\begin{equation}
\label{eq5.17}
\begin{split}
&\sum\limits_{k=0}^{\infty}2^{kp(1-\frac{n-d}{p})}\int\limits_{S}\Bigl(\mathcal{E}_{\mu_{k}}(F|_{S},Q(x,2^{-k}))\Bigr)^{p}\,d\mu_{k}(x) \le C \int\limits_{\mathbb{R}^{n}}\Bigl(\operatorname{M}[|\nabla F|^{q}](x)\Bigr)^{\frac{p}{q}}\,d\mathcal{L}_{n}(x)\\
&\le C \|F|L^{1}_{p}(\mathbb{R}^{n})\|^{p}.
\end{split}
\end{equation}

The lemma is proved.

\begin{Lm}
\label{Lm5.2}
Let $d \in [0,n]$, $p \in (1,\infty)$ and $p > n-d$. Let $S$ be a closed set in $\mathbb{R}^{n}$ with $\operatorname{dim}_{H}S \geq d$. Assume that there exists a $d$-regular system of measures $\{\mu_{k}\}=\{\mu_{k}\}_{k \in \mathbb{N}_{0}}$ on $S$. Assume that $\partial S$ is porous. Assume that $f \in L_{1}^{\operatorname{loc}}(\mathbb{R}^{n},\mu_{k})$ for every $k \in \mathbb{N}_{0}$. Then for every $k \in \mathbb{N}$, $k \geq 4$
\begin{equation}
\begin{split}
\label{eq5.18'}
&\sum\limits_{\substack{\varkappa \in \mathcal{I} \\ r_{\varkappa} \le 2^{-k}}}\mathcal{L}_{n}(Q_{\varkappa})\left(\mathcal{E}_{\mu_{k}}(f,Q(\widetilde{x}_{\varkappa},2^{-k}))\right)^{p}\\
&\le C 2^{(k-4)(d-n)}\int\limits_{\Sigma_{k-4}}\Bigl(\mathcal{E}_{\mu_{k}}(f,Q(x,2^{-(k-4)}))\Bigr)^{p}\,d\mu_{k-4}(x).
\end{split}
\end{equation}
The constant $C > 0$ does not depend on $f$ and $k$.
\end{Lm}

\textbf{Proof.} Fix an arbitrary natural number $k \geq 4$. We set $\mathcal{I}_{k}:=\{\varkappa \in \mathcal{I}| r_{\varkappa} \le 2^{-k}\}$. Using Theorem \ref{Th2.5} and Remark \ref{Rem2.7}, we find an index set $\widehat{\mathcal{I}}_{k} \subset \mathcal{I}_{k}$ such that

\begin{equation}
\label{eq5.19'}
\bigcup\limits_{\varkappa \in \mathcal{I}_{k}}Q\bigl(\widetilde{x}_{\varkappa},\frac{1}{2^{k}}\bigr) \subset \bigcup\limits_{\varkappa \in \widehat{\mathcal{I}}_{k}}Q\bigl(\widetilde{x}_{\varkappa},\frac{5}{2^{k}}\bigr).
\end{equation}

Note that if $Q\bigl(\widetilde{x}_{\varkappa'},\frac{1}{2^{k}}\bigr) \cap Q\bigl(\widetilde{x}_{\varkappa},\frac{5}{2^{k}}\bigr) \neq \emptyset$ for some $\varkappa' \in \mathcal{I}_{k}$, $\varkappa \in \widehat{\mathcal{I}}_{k}$, then $Q\bigl(\widetilde{x}_{\varkappa'},\frac{1}{2^{k}}\bigr) \subset Q\bigl(\widetilde{x}_{\varkappa},\frac{7}{2^{k}}\bigr)$ for such $\varkappa,\varkappa'$. Using Remark \ref{Rem5.1}, and arguing as in the proof of Lemma \ref{Lm3.5}, it is easy to show that for such $\varkappa,\varkappa'$ we have the estimate

\begin{equation}
\label{eq5.20'}
\mathcal{E}_{\mu_{k}}(f,Q(\widetilde{x}_{\varkappa'},2^{-k})) \le C \mathcal{E}_{\mu_{k}}\Bigl(f,Q\bigl(\widetilde{x}_{\varkappa},\frac{7}{2^{k}}\bigr)\Bigr).
\end{equation}

Using the same arguments as in the proof of Lemma \ref{Lm2.5}, we find that for every $\varkappa \in \widehat{\mathcal{J}}_{k}$

\begin{equation}
\label{eq5.21'}
\sum\limits_{\substack{\varkappa' \in \mathcal{I}_{k}\\ Q\bigl(\widetilde{x}_{\varkappa'},\frac{1}{2^{k}}\bigr) \cap Q\bigl(\widetilde{x}_{\varkappa},\frac{5}{2^{k}}\bigr)}}\mathcal{L}_{n}(Q_{\varkappa'}) \le C 2^{-kn}.
\end{equation}

Combining \eqref{eq5.19'}, \eqref{eq5.20'} and \eqref{eq5.21'}, we obtain

\begin{equation}
\label{eq5.22'}
\sum\limits_{\varkappa \in \mathcal{I}_{k}}\mathcal{L}_{n}(Q_{\varkappa})\left(\mathcal{E}_{\mu_{k}}(f,Q(\widetilde{x}_{\varkappa},2^{-k}))\right)^{p} \le C 2^{-kn}\sum\limits_{\varkappa \in \widehat{\mathcal{I}}_{k}}\left(\mathcal{E}_{\mu_{k}}\Bigl(f,Q\bigl(\widetilde{x}_{\varkappa},\frac{7}{2^{k}}\bigr)\Bigr)\right)^{p}.
\end{equation}

It is clear that $Q\bigl(\widetilde{x}_{\varkappa},\frac{7}{2^{k}}\bigr) \subset Q(x,\frac{15}{2^{k}})$ for every $x \in S \cap Q\bigl(\widetilde{x}_{\varkappa},\frac{7}{2^{k}}\bigr)$. Hence, using Remark \ref{Rem5.1}, \eqref{eq3.6}, and arguing as in the proof of Lemma \ref{Lm3.5}, we have

\begin{equation}
\begin{split}
\label{eq5.23'}
&\left(\mathcal{E}_{\mu_{k}}\Bigl(f,Q\bigl(\widetilde{x}_{\varkappa},\frac{7}{2^{k}}\bigr)\Bigr)\right)^{p} \le C\inf\limits_{x \in Q\bigl(\widetilde{x}_{\varkappa},\frac{7}{2^{k}}\bigr) \cap S}\left(\mathcal{E}_{\mu_{k}}\Bigl(f,Q\bigl(x,\frac{2^{4}}{2^{k}}\bigr)\Bigr)\right)^{p} \\
&\le C \inf\limits_{x \in Q\bigl(\widetilde{x}_{\varkappa},\frac{7}{2^{k}}\bigr) \cap S}\left(\mathcal{E}_{\mu_{k-4}}\Bigl(f,Q\bigl(x,2^{-(k-4)}\bigr)\Bigr)\right)^{p}\\
&\le C 2^{(k-4)d}\int\limits_{Q\bigl(\widetilde{x}_{\varkappa},\frac{7}{2^{k}}\bigr) \cap S}\left(\mathcal{E}_{\mu_{k-4}}\Bigl(f,Q\bigl(x,2^{-(k-4)}\bigr)\Bigr)\right)^{p}\,d\mu_{k-4}(x).
\end{split}
\end{equation}

It is clear that $Q\bigl(\widetilde{x}_{\varkappa},\frac{7}{2^{k}}\bigr) \cap S \subset \Sigma_{k-4}$. Furthermore the overlapping multiplicity of the sets
$Q\bigl(\widetilde{x}_{\varkappa},\frac{7}{2^{k}}\bigr) \cap S$ is finite and independent on $k$.
Hence, substituting \eqref{eq5.23'} into \eqref{eq5.22'} and, using Lemma \ref{Lm2.9}, we obtain \eqref{eq5.18'}.

The lemma is proved.

\begin{Lm}
\label{Lm5.3}
Let $\lambda \in (0,1)$, $d \in [0,n]$, $p \in (1,\infty)$ and $p > n-d$. Let $S$ be a closed set in $\mathbb{R}^{n}$ with $\operatorname{dim}_{H}S \geq d$. Assume that there exists a $d$-regular system of measures $\{\mu_{k}\}=\{\mu_{k}\}_{k \in \mathbb{N}_{0}}$ on $S$. Assume that $\partial S$ is porous. Assume that $f \in L_{1}^{\operatorname{loc}}(\mathbb{R}^{n},\mu_{k})$ for every $k \in \mathbb{N}_{0}$. Then
\begin{equation}
\begin{split}
\label{eq5.24'}
&\sum\limits_{k=0}^{\infty}2^{k(d-n)}\int\limits_{S_{k}(\lambda)}\Bigl(f^{\sharp}_{\{\mu_{k}\}}(x,2^{-k})\Bigr)^{p}\,d\mu_{k}(x)\\
&\le C \sum\limits_{k=0}^{\infty}2^{kp(1-\frac{n-d}{p})}\int\limits_{\Sigma_{k}}\Bigl(\mathcal{E}_{\mu_{k}}(f,Q(x,2^{-k}))\Bigr)^{p}\,d\mu_{k}(x) + C \int\limits_{S}|f(x)|^{p}\,d\mu_{0}(x).
\end{split}
\end{equation}
The constant $C > 0$ in \eqref{eq5.24'} does not depend on $f$.
\end{Lm}

\textbf{Proof.} It is clear that for every $\varkappa \in \mathcal{I}$ one can choose $k_{\varkappa} \in \mathbb{N}_{0}$ such that

\begin{equation}
\label{eq5.25'}
f^{\sharp}_{\{\mu_{k}\}}(\widetilde{x}_{\varkappa},r_{\varkappa}) \le  \frac{2^{n}}{2^{-k_{\varkappa}}}\mathcal{E}_{\mu_{k_{\varkappa}}}(f,Q(\widetilde{x}_{\varkappa},2^{-k_{\varkappa}})).
\end{equation}

Given $k \in \mathbb{N}_{0}$, let $\mathcal{I}_{k}:=\{\varkappa \in \mathcal{I}| k_{\varkappa} = k\}$. It is clear that

\begin{equation}
\label{eq5.26'}
\mathcal{I} \subset \bigcup\limits_{k \in \mathbb{N}_{0}}\mathcal{I}_{k}.
\end{equation}

Assume that $r_{\varkappa} \geq 2^{-4}$. Then, using H\"older inequality and \eqref{eq3.5}, \eqref{eq3.6}, we easily get

\begin{equation}
\label{eq5.27'}
\Bigl(f^{\sharp}_{\{\mu_{k}\}}(\widetilde{x}_{\varkappa},r_{\varkappa})\Bigr)^{p} \le C \int\limits_{Q(\widetilde{x}_{\varkappa},1)}|f(y)|^{p}\,d\mu_{0}(y).
\end{equation}

Hence, using Lemma \ref{Lm2.5} with $c=1$ and $d\mathfrak{m}(x)=|f(x)|^{p}d\mu_{0}(x)$, it is easy to find that

\begin{equation}
\begin{split}
\label{eq5.28'}
&\sum\limits_{\substack{\varkappa \in \mathcal{I}\\ r_{\varkappa} \geq 2^{-4}}}\mathcal{L}_{n}(Q_{\varkappa})\Bigl(f^{\sharp}_{\{\mu_{k}\}}(\widetilde{x}_{\varkappa},r_{\varkappa})\Bigr)^{p} \\
&\le C \sum\limits_{\substack{\varkappa \in \mathcal{I}\\ r_{\varkappa} \geq 2^{-4}}}\mathcal{L}_{n}(Q_{\varkappa})\int\limits_{Q(\widetilde{x}_{\varkappa},1)}|f(y)|^{p}\,d\mu_{0}(y) \le C \int\limits_{S}|f(y)|^{p}\,d\mu_{0}(y).
\end{split}
\end{equation}

Combining \eqref{eq5.18'} and \eqref{eq5.28'}, we obtain the estimate

\begin{equation}
\label{eq5.29'}
\begin{split}
&\sum\limits_{\varkappa \in \mathcal{I}}\mathcal{L}_{n}(Q_{\varkappa})\Bigl(f^{\sharp}_{\{\mu_{k}\}}(\widetilde{x}_{\varkappa},r_{\varkappa})\Bigr)^{p} \le 2^{n} \sum\limits_{k=0}^{\infty}\sum\limits_{\varkappa \in \mathcal{I}_{k}}\mathcal{L}_{n}(Q_{\varkappa})2^{kp}\mathcal{E}_{\mu_{k}}(f,Q(\widetilde{x}_{\varkappa},2^{-k}))\\
& \le C \sum\limits_{k=0}^{\infty}2^{k(p-n+d)}\int\limits_{\Sigma_{k}}\Bigl(\mathcal{E}_{\mu_{k}}(f,Q(x,2^{-k}))\Bigr)^{p}\,d\mu_{k}(x) + C \int\limits_{S}|f(x)|^{p}\,d\mu_{0}(x).
\end{split}
\end{equation}

Now the lemma follows from \eqref{eq4.41'} and \eqref{eq5.29'}.

Now we are ready to formulate the main result of this subsection. Combining the following theorem with Remark \ref{Rem1.1}, we get Theorem \ref{Th1.2}.

\begin{Th}
\label{Th5.1}
Let $\lambda \in (0,1)$, $d \in [0,n]$, $p \in (1,\infty)$ and $p > n-d$. Let $S$ be a closed set in $\mathbb{R}^{n}$ with $\operatorname{dim}_{H}S \geq d$. Assume that there exists a $d$-regular system of measures $\{\mu_{k}\}=\{\mu_{k}\}_{k \in \mathbb{N}_{0}}$ on $S$. Assume that $\partial S$ is porous. Then $f$ belongs to the $d$-trace space $W_{p}^{1}(\mathbb{R}^{n})|_{S}$ if and only if there exists a set $S' \subset S$ with $\mathcal{H}^{d}(S \setminus S') = 0$ such that
\begin{equation}
\notag
\lim\limits_{r \to 0}\fint\limits_{Q(x,r) \cap S} |f(x)-f(z)|\,d\mu_{k(r)}(z) = 0, \quad \forall x \in S',
\end{equation}
and
\begin{equation}
\notag
\begin{split}
&\mathcal{R}_{S,p}[f]:=\left(\int\limits_{S}|f(x)|^{p}\,d\mu_{0}(x)\right)^{\frac{1}{p}}+\left(\int\limits_{S}\bigl(f^{\sharp}_{\{\mu_{k}\}}(x)\bigr)^{p}\,d\mathcal{L}_{n}(x)\right)^{\frac{1}{p}}\\
&+\left(\sum\limits_{k=0}^{\infty}2^{kp(1-\frac{n-d}{p})}\int\limits_{\Sigma_{k}}\Bigl(\mathcal{E}_{\mu_{k}}(f,Q(x,2^{-k}))\Bigr)^{p}\,d\mu_{k}(x)\right)^{\frac{1}{p}} <
\infty.
\end{split}
\end{equation}
Furthermore
\begin{equation}
\notag
\|f|W_{p}^{1}(\mathbb{R}^{n})|_{S}\|  \sim \mathcal{R}_{S,p}[f],
\end{equation}
and there exists a bounded linear operator $\operatorname{Ext}:W_{p}^{1}(\mathbb{R}^{n})|_{S} \to  W_{p}^{1}(\mathbb{R}^{n})$ such that $\operatorname{Tr}|_{S} \circ \operatorname{Ext} = \operatorname{Id}$ on
$W^{1}_{p}(\mathbb{R}^{n})|_{S}$.
\end{Th}

\textbf{Proof.} The theorem follows from Theorem \ref{Th4.5} and Lemmas \ref{Lm5.1}, \ref{Lm5.3}.

\section{Examples}

In this section we present several useful examples which illustrate our main results.

\textbf{Example 6.1.} Let $S$ be a closed Ahlfors $n$-regular subset of $\mathbb{R}^{n}$ and $p > 1$. In this case one can take $\mu_{k} = \mathcal{L}_{n}\lfloor S$ for every $k \in \mathbb{N}_{0}$ to obtain an $n$-regular system of measures on $S$.

We show that $\mathcal{BN}_{S,p,\lambda}[f] \le C (\mathcal{SN}_{S,p}[f]+\|f|L_{p}(S)\|)$ with a constant $C > 0$ independent on $f$. We use estimate \eqref{eq4.41'}, Lemmas \ref{Lm2.5}, \ref{Lm2.6} and \ref{Lm3.5}. We also use simple inclusions $Q(\widetilde{x}_{\varkappa},3r_{\varkappa}) \supset Q(x,2r_{\varkappa})\supset Q(\widetilde{x}_{\varkappa},r_{\varkappa})$ for every $x \in Q(\widetilde{x}_{\varkappa},r_{\varkappa}) \cap S$. This gives
\begin{equation}
\begin{split}
\notag
&\Bigl(\mathcal{BN}_{S,p,\lambda}[f]\Bigr)^{p} \le C \sum\limits_{\varkappa \in \mathcal{I}}\mathcal{L}_{n}(\mathcal{U}_{\varkappa})\left(f^{\sharp}_{\{\mu_{k}\}}(\widetilde{x}_{\varkappa},r_{\varkappa})\right)^{p} + C \int\limits_{S}|f(x)|^{p}\,d\mathcal{L}_{n}(x)\\
&\le C \sum\limits_{\varkappa \in \mathcal{I}}\mathcal{L}_{n}(\mathcal{U}_{\varkappa})\inf\limits_{x \in \mathcal{U}_{\varkappa}}\left(f^{\sharp}_{\{\mu_{k}\}}(x,r_{\varkappa})\right)^{p}+C \sum\limits_{\varkappa \in \mathcal{I}}\mathcal{L}_{n}(Q_{\varkappa})\fint\limits_{Q(\widetilde{x}_{\varkappa},3)}|f(y)|^{p}\,d\mathcal{L}_{n}(y)\\
&+C \int\limits_{S}|f(x)|^{p}\,d\mathcal{L}_{n}(x) \le C \int\limits_{S}\bigl(f^{\sharp}_{\{\mu_{k}\}}(x)\bigr)^{p}\,d\mathcal{L}_{n}(x)+C \int\limits_{S}|f(x)|^{p}\,d\mathcal{L}_{n}(x)\\
&=C (\mathcal{SN}_{S,p}[f])^{p}+C\|f|L_{p}(S)\|^{p}.
\end{split}
\end{equation}

In this case we have
\begin{equation}
\notag
f^{\sharp}_{\{\mu_{k}\}}(x)=\sup\limits_{r \in (0,1)}\frac{1}{r}\fint\limits_{Q(x,r) \cap S}\Bigl|f(y) - \fint\limits_{Q(x,r) \cap S}f(z)\,d\mathcal{L}_{n}(z)\Bigr|\,d\mathcal{L}_{n}(y).
\end{equation}
To simplify notation we set $f^{\sharp}_{S}:=f^{\sharp}_{\{\mu_{k}\}}$. Such notation was used in \cite{Shv1}.

The estimate above together with the fact that $\mathcal{L}_{n}$-a.e. point $x \in S$ is a Lebesgue point of a function $f \in L_{p}(S)$ allows us to obtain simple trace criterion.
\textit{Namely, a function $f:S \to \mathbb{R}$ belongs to the $n$-trace space $W_{p}^{1}(\mathbb{R}^{n})|_{S}$ if and only if
\begin{equation}
\|f|L_{p}(S,\mathcal{L}_{n})\|+\|f^{\sharp}_{S}|L_{p}(S,\mathcal{L}_{n})\| < +\infty.
\end{equation}
Moreover,
\begin{equation}
\begin{split}
\notag
&\|f|L_{p}(S,\mathcal{L}_{n})\|+\|f^{\sharp}_{S}|L_{p}(S,\mathcal{L}_{n})\| \sim \|f|W_{p}^{1}(\mathbb{R}^{n})|_{S}\|
\end{split}
\end{equation}
and operator $\operatorname{Ext}$ constructed in \eqref{eq4.3} is a bounded linear extension operator $\operatorname{Ext}:W_{p}^{1}(\mathbb{R}^{n})|_{S} \to W_{p}^{1}(\mathbb{R}^{n})$ such that $\operatorname{Ext} \circ \operatorname{Tr}|_{S} = \operatorname{Id}$ on $W_{p}^{1}(\mathbb{R}^{n})|_{S}$.}

This result coincide with that of obtained in \cite{Shv1} in the case of first order Sobolev spaces. In the case $S=\mathbb{R}^{n}$ we have $\|f^{\sharp}|L_{p}(\mathbb{R}^{n})\| \sim \|L_{p}^{1}(\mathbb{R}^{n})\|$. Such equivalence motivated us to call $\mathcal{SN}_{S,p}[f]$ the "Sobolev part" of the norm (see Remark 4.3).

\textbf{Example 6.2.} Let $d \in [0,n)$ and $p \in (\max\{1,n-d\},\infty)$. Let $S$ be a closed Ahlfors $d$-regular subset of $\mathbb{R}^{n}$. In this case there exists a simple $d$-regular system of measures on $S$. More precisely, we set $\mu_{k}=\mathcal{H}^{d}\lfloor S$ for every $k \in \mathbb{N}_{0}$. Clearly $\mathcal{L}_{n}(S)=0$, because $d < n$. Furthermore, $\operatorname{int}S=\emptyset$ and $\partial S$ is porous (see Remark \ref{Rem3.4}).

Note that the measure $\mathcal{H}^{d}\lfloor S$ is Radon. Hence, by Theorem 1 in section 1.7.1 of \cite{Evans} we conclude that if $f \in L^{\operatorname{loc}}_{1}(S,\mathcal{H}^{d}\lfloor S)$ then

\begin{equation}
\notag
\fint\limits_{Q(x,r) \cap S}|f(x)-f(y)|\,d\mathcal{H}^{d}(y) = 0
\end{equation}
for $\mathcal{H}^{d}$-almost every $x \in S$.

Now we apply Theorem \ref{Th1.2} and take into account all facts mentioned above. This gives us very simple criterion.

\textit{Namely, a function $f : S \to \mathbb{R}$ belongs to the $d$-trace space $W_{p}^{1}(\mathbb{R}^{n})|_{S}$ if and only if
\begin{equation}
\begin{split}
\notag
&\|f|L_{p}(S,\mathcal{H}^{d}\lfloor S)\|+\left(\sum\limits_{k=0}^{\infty}2^{kp(1-\frac{n-d}{p})}\int\limits_{S}\mathcal{E}_{\mathcal{H}^{d}}(f,Q(x,2^{-k}))\,d\mathcal{H}^{d}(x)\right)^{\frac{1}{p}} < \infty.
\end{split}
\end{equation}
Moreover,
\begin{equation}
\begin{split}
\notag
&\|f|L_{p}(S,\mathcal{H}^{d}\lfloor S)\|\\
&+\left(\sum\limits_{k=0}^{\infty}2^{kp(1-\frac{n-d}{p})}\int\limits_{S}\mathcal{E}_{\mathcal{H}^{d}}(f,Q(x,2^{-k}))\,d\mathcal{H}^{d}(x)\right)^{\frac{1}{p}} \sim \|f|W_{p}^{1}(\mathbb{R}^{n})|_{S}\|
\end{split}
\end{equation}
and operator $\operatorname{Ext}$ constructed in \eqref{eq4.3} is a bounded linear extension operator $\operatorname{Ext}:W_{p}^{1}(\mathbb{R}^{n})|_{S} \to W_{p}^{1}(\mathbb{R}^{n})$ such that $\operatorname{Ext} \circ \operatorname{Tr}|_{S} = \operatorname{Id}$ on $W_{p}^{1}(\mathbb{R}^{n})|_{S}$.}

Note that this result coincide with that of obtained \cite{Ihn} in the case of first order Sobolev spaces.

In the simplest case $S = \mathbb{R}^{d}$, $d=1,..,n-1$ one can recognize the classical result. Namely, the trace space of the Sobolev first order space to the plane $\mathbb{R}^{d}$ is the classical Besov space $B^{1-\frac{n-d}{p}}_{p,p}(\mathbb{R}^{d})$. This fact together with Theorem 1.1 implies that
$\mathcal{BN}_{\mathbb{R}^{d},p,\lambda}[f] \sim \|f|B^{1-\frac{n-d}{p}}_{p,p}(\mathbb{R}^{d})\|$ for every $\lambda \in (0,1)$. Such equivalence motivated us to call $\mathcal{BN}_{S,p,\lambda}[f]$ the "Besov part" of the norm (see Remark 4.3).

\textbf{Example 6.3.} Let $\beta: [0,+\infty) \to [0,+\infty)$ be an arbitrary continuous strictly increasing function such that $\beta(0)=0$ and $\beta(t) > 0$ for every $t > 0$.
By $\beta^{-1}$ we denote the inverse function, i.e. $\beta^{-1}\circ \beta = \operatorname{id}$ on $[0,+\infty)$.
Consider the closed single cusp in $\mathbb{R}^{n}$
\begin{equation}
\label{eq5.30}
G^{\beta}:=\{x=(x',x_{n})| \max\limits_{i=1,..,,n-1}|x_{i}| \le \beta(x_{n})\}.
\end{equation}
Given a number $k \in \mathbb{N}_{0}$, we consider also the sets
\begin{equation}
\begin{split}
\label{eq5.30'}
&G^{\beta}_{k}:=\{x=(x',x_{n})| \max\limits_{i=1,..,,n-1}|x_{i}| \le \beta(x_{n}), 0 \le x_{n} \le \beta^{-1}(2^{-k})\} \cup\\
&\cup \{x=(x',x_{n})| \beta(x_{n}) \geq \max\limits_{i=1,..,,n-1}|x_{i}| > \beta(x_{n})-2^{-k},  x_{n} > \beta^{-1}(2^{-k})\}.
\end{split}
\end{equation}
Recall Definition \ref{Def5.2}. It is clear that $G^{\beta}_{k}$ coincides with $\Sigma_{k}(G^{\beta})$.

Given a number $k \in \mathbb{N}_{0}$, consider the measure $d\mu_{k}(x)=w^{\beta}_{k}(x)\,d\mathcal{L}_{n}(x)$, where
\begin{equation}
\label{eq5.31}
w^{\beta}_{k}(x):=w^{\beta}_{k}(x',x_{n}):=
\begin{cases}
&(\beta(x_{n}))^{1-n}, \quad x_{n} \in [0,\beta^{-1}(2^{-k})];\\
&2^{k(n-1)}, \quad x_{n} \geq \beta^{-1}(2^{-k});\\
&0, \quad x \notin G^{\beta}.
\end{cases}
\end{equation}

Using elementary geometrical observations, it is easy to see  that for every $x=(x',x_{n}) \in G^{\beta}$ and $r \in (0,1)$
\begin{equation}
\label{eq6.4'''}
\mu_{k}(Q((x',0),r)) \geq \mu_{k}(Q((x',x_{n}),r)).
\end{equation}
On the other hand, using \eqref{eq5.31} and monotonicity and continuity properties of $\beta$, it is easy to show that for every $x=(x',x_{n}) \in G^{\beta}$
\begin{equation}
\label{eq6.5''}
\mu_{k}(Q((x',0),2^{-k})) \le C(\beta)\mu_{k}(Q(x,2^{-k})).
\end{equation}
But direct computations give for every $x=(x',x_{n}) \in G^{\beta}$
\begin{equation}
\label{eq6.6''}
\mu_{k}(Q((0,x_{n}),r) = c(n)\int\limits_{x_{n}-r}^{x_{n}+r}(\beta(t))^{n-1}\frac{1}{(\beta(t))^{n-1}}\,dt \sim c(n)r.
\end{equation}

Combining \eqref{eq6.4'''}--\eqref{eq6.6''} we see that the system of measures $\{\mu_{k}\}_{k \in \mathbb{N}_{0}}$ is $1$-regular on $G^{\beta}$.

Recall item (2) of Example 2.1 and Example 3.1. Thus, we see that the set $G^{\beta}$ is $1$-thick and has the porous boundary. Consider slightly relaxed definition of the trace of $F \in W_{p}^{1}(\mathbb{R}^{n})$ to the set $G^{\beta}$. Namely,  we write $F|_{G^{\beta}}=f$ if $F(x)=f(x)$ for $\mathcal{L}_{n}$-a.e. $x \in G^{\beta}$, i.e. $f$ is the $n$-trace of $F$. Then we clearly derive from Theorem 1.2 the following criterion. 

\textit{Let $p > n-1$, then a function $f: G^{\beta} \to \mathbb{R}^{n}$ belongs to the $n$-trace space $W_{p}^{1}(\mathbb{R}^{n})|_{G^{\beta}}$ if and only if
\begin{equation}
\label{6.8}
\begin{split}
&\Bigl(\int\limits_{G^{\beta}}\bigl(f^{\sharp}_{\{\mu_{k}\}}(x)\bigr)^{p}\,d\mathcal{L}_{n}(x)\Bigr)^{\frac{1}{p}}+\Bigl(\int\limits_{G^{\beta}}\beta_{0}(x)|f(x)|^{p}\,d\mathcal{L}_{n}(x)\Bigr)^{\frac{1}{p}}\\
&+\Bigl(\sum\limits_{k=0}^{\infty}2^{kp(1-\frac{n-1}{p})}\int\limits_{G^{\beta}_{k}}\beta_{k}(x)\bigl(\mathcal{E}_{\mu_{k}}(f,Q(x,2^{-k}))\bigr)^{p}\,d\mathcal{L}_{n}(x)\Bigr)^{\frac{1}{p}} <\infty.
\end{split}
\end{equation}
Furthermore, the left-hand side of \eqref{6.8} gives an equivalent norm in the $n$-trace space $W_{p}^{1}(\mathbb{R}^{n})|_{G^{\beta}}$ and operator $\operatorname{Ext}$ constructed in \eqref{eq4.3} is a bounded linear operator from $W_{p}^{1}(\mathbb{R}^{n})|_{G^{\beta}}$ to $W_{p}^{1}(\mathbb{R}^{n})$ such that $\operatorname{Tr}|_{G^{\beta}} \circ \operatorname{Ext} = \operatorname{Id}$ on  $W_{p}^{1}(\mathbb{R}^{n})|_{G^{\beta}}$.}

As far as we know even this particular result is knew and could not be obtained by previously known techniques.

\textbf{Example 6.4} In this example we restrict ourselves to the case $p \in (n,\infty)$. It is well known that in this case every element $F \in W^{1}_{p}(\mathbb{R}^{n})$ 
has a continuous representative (\cite{Zi}). As a result we can define the $0$-trace (or just trace) of the element $F$ to an \textit{arbitrary set} $S \subset \mathbb{R}^{n}$.  

It is easy to see that every nonempty set $S$ is $0$-thick. Hence, if $S$ is closed there exists a $0$-regular system of measures on $S$. 
Indeed in this case we can avoid the using of delicate arguments of Theorem 3.1. More precisely, for every $k \in \mathbb{N}_{0}$ we consider 
an arbitrary maximal $2^{-k}$ separated set $\{x^{k}_{j}\}_{j \in \mathcal{J}^{k}}$ (recall Definition 2.9). For every $k \in \mathbb{N}$ we set
\begin{equation}
\notag
\mu_{k}:=\sum\limits_{j \in \mathcal{J}^{k}}\delta_{x^{k}_{j}},
\end{equation}
where $\delta_{x}$ is the Dirac measure concentrated on $x \in \mathbb{R}^{n}$. 
It is easy to show that such system of measures $\{\mu_{k}\}=\{\mu_{k}\}_{k \in \mathbb{N}_{0}}$ satisfies requirements (3.4)--(3.6).

Note that for $r \in [2^{-k},2^{-k+1})$ ($k \in \mathbb{N}$)
\begin{equation}
\label{eq6.9'''}
\begin{split}
&\fint\limits_{Q(x,r)}\Bigl|f(y)-\fint\limits_{Q(x,r)} f(z)\,d\mu_{k(r)}(z)\Bigr|d\mu_{k}(y) \sim C(n)\sum\limits_{\substack{j \in \mathcal{J}^{k}\\\|x-x^{k}_{j}\| \le r}}\sum\limits_{\substack{j' \in \mathcal{J}^{k}\\\|x-x^{k}_{j}\| \le r}}|f(x^{k}_{j})-f(x^{k}_{j'})|,\\
& \sim C'(n) \max\limits_{\substack{j,j' \in \mathcal{J}^{k}\\ x^{k}_{j},x^{k}_{j'} \in Q(x,r)}}|f(x^{j}_{k})-f(x^{k}_{j'})|=:C'(n)\mathcal{A}(f,Q(x,r)).
\end{split}
\end{equation}
where $C(n) > 0$ depends only on the amount of points $x^{k}_{j}$ in the cube $Q(x,r)$ (clearly such amount depends only on $n$)

From \eqref{eq6.9'''} we obtain in this special case for every $t \in [0,1)$ 
\begin{equation}
f^{\sharp}_{\{\mu_{k}\}}(x,t) \sim \sup\limits_{r \in (t,1)}\frac{1}{r}\mathcal{A}(f,Q(x,r))=:f^{\natural}(x,t).
\end{equation}
In the case $t=0$ we simply write $f^{\natural}(x)$.

Hence we can reformulate our theorem 1.1 in the case $p > n$ as follows.

\textit{Let $S$ be an arbitrary closed subset of $\mathbb{R}^{n}$. Let $p \in (n,\infty)$. A continuous on $S$ function $f:S \to \mathbb{R}$ belongs to the 
$0$-trace space $W_{p}^{1}(\mathbb{R}^{n})|_{S}$ if and only if for some $\lambda \in (0,1)$
\begin{equation}
\label{eq6.10'''}
\begin{split}
&\Bigl(\sum\limits_{j \in \mathcal{J}^{0}}|f(x^{0}_{j})|^{p}\Bigr)^{\frac{1}{p}}+\Bigl(\int\limits_{S}(f^{\natural}(x))^{p}d\mathcal{L}_{n}(x)\Bigr)^{\frac{1}{p}}\\
&+\Bigl(\sum\limits_{k=0}^{\infty}\sum\limits_{j \in \mathcal{J}^{k}}(f^{\natural}(x^{k}_{j},2^{-k}))^{p}\Bigr)^{\frac{1}{p}} < \infty.
\end{split}
\end{equation}
Furthermore, the left-hand side of \eqref{eq6.10'''} gives an equivalent norm on the $0$-trace space  $W_{p}^{1}(\mathbb{R}^{n})|_{S}$
and operator $\operatorname{Ext}$ constructed in \eqref{eq4.3}  is a bounded linear extension operator $\operatorname{Ext}:W_{p}^{1}(\mathbb{R}^{n})|_{S} \to W_{p}^{1}(\mathbb{R}^{n})$ such that $\operatorname{Ext} \circ \operatorname{Tr}|_{S} = \operatorname{Id}$ on $W_{p}^{1}(\mathbb{R}^{n})|_{S}$.}

Note that this result differs from the main result of \cite{Shv2} and provides an alternative characterization of the trace space in the case $p > n$ (compare with Theorems 1.2 and 2.5 of \cite{Shv2}).

\section{Appendix}

The aim of this section is to give detailed explanations of Example 2.1. Furthermore we include the proof of Theroem 2.2 and 
the proof of our refined version of the Frostman-type Lemma (Theorem 3.1).

\begin{Lm}
\label{Lm6.1}
Let $d \in [0,n]$. Let $E$ be an arbitrary $d$-thick set. Then, the closure $\overline{E}$ of the set $E$ is $d$-thick.
\end{Lm}

\textbf{Proof.} According to the Definition \ref{Def2.3} it is sufficient to prove that

\begin{equation}
\notag
\mathcal{H}^{d}_{\infty}(Q(x,r) \cap \overline{E}) \geq \varepsilon r^{d}
\end{equation}

for every $x \in \partial E$ and every $r \in (0,1]$ with $\varepsilon > 0$ independent on $x$ and $r$.

Let $x \in \partial E$. Take an arbitrary sequence $\{y_{j}\}_{j=1}^{\infty}$ converging to $x$. Fix a number $c > 1$. Given $r > 0$, for sufficiently big indexes $j$

\begin{equation}
\notag
Q(y_{j},r) \subset Q(x,cr).
\end{equation}

Hence, we obtain

\begin{equation}
\notag
\mathcal{H}^{d}_{\infty}(Q(x,cr) \cap \overline{E}) \geq \mathcal{H}^{d}_{\infty}(Q(y_{j},r) \cap E) \geq \frac{\varepsilon}{c^{d}} (cr)^{d}.
\end{equation}

Using the fact that $c > 1$ was chosen arbitrary close to $1$, we conclude.

1) Let $\Omega \subset \mathbb{R}^{n}$ be an open path-connected set. We are going to prove that the set $\overline{\Omega}$ is $1$-thick. Due to the Lemma \ref{Lm4.1} it is sufficient to show that $\Omega$ is $1$-thick.

Fix a point $x \in \Omega$. Let $Q=Q(x,r)$ be a cube with $0 < r \le 1$. Consider two cases.

In the first case $\operatorname{diam}\Omega > r/2$. Then there is a point $y \in \Omega \setminus Q$. Hence there is a curve $\gamma_{x,y}$ which connects $x$ and $y$. Let $\{B_{j}\}_{j \in \mathbb{N}}=\{B(x_{j},r_{j})\}_{j \in
\mathbb{N}}$ be an arbitrary covering of $Q \cap \Omega$ for which

\begin{equation}
\label{eq6.1}
\sum\limits_{j \in \mathbb{N}}r_{j} \le 2\mathcal{H}^{1}_{\infty}(\Omega \cap Q)
\end{equation}

Choose index set $\mathcal{A} \subset \mathbb{N}$ such that $\gamma_{x,y} \cap B_{j} \neq \emptyset$ for every $j \in \mathcal{A}$ and $\gamma_{x,y} \subset \cup_{j \in \mathcal{A}}B_{j}$.

Consider projections $\gamma^{i}_{x,y}$, $i = 1,..,n$ of our curve and projections $B^{i}_{j}$ of balls (from the covering) to coordinate axes. It is clear that there exists $i_{0} \in \{1,...,n\}$ for
which $\mathcal{L}_{1}(\gamma^{i_{0}}_{x,y}) \geq r/2$. By the construction the family of intervals $\{B^{i_{0}}_{j}\}_{j \in \mathcal{A}}$ is a covering of $\gamma^{i_{0}}_{x,y}$. Hence, from \eqref{eq6.1} we derive

\begin{equation}
\label{eq6.2}
\begin{split}
&\mathcal{H}^{1}_{\infty}(\Omega \cap Q) \geq \frac{1}{2}\sum\limits_{j \in \mathbb{N}}r_{j} \geq \frac{1}{2}\sum\limits_{j \in \mathcal{A}}r_{j} \geq \frac{\mathcal{L}_{1}(\gamma^{i_{0}}_{x,y})}{2} \geq
\frac{r}{2}.
\end{split}
\end{equation}

In the second case $\operatorname{diam}\Omega \le r/2$. Hence for every $x \in \Omega$ we have $\Omega \subset Q(x,r)$ and $\mathcal{H}^{1}_{\infty}(Q(x,r) \cap \Omega) \geq \mathcal{H}^{1}_{\infty}(\Omega)$.

We set $\varepsilon:=\min\{\mathcal{H}^{1}_{\infty}(\Omega), 1/2\}$. As a result, for every $x \in \Omega$ and  for every $r \in (0,1]$ we obtain

$$
\mathcal{H}^{1}_{\infty}(Q(x,r) \cap S) \geq \varepsilon r.
$$

This proves the claim.

2) Let $d \in [0,n]$. Let $S$ be a set which is Ahlofrs  $d$-regular. Consider a cube $Q=Q(x,r)$ with $x \in S$ and $0 < r \le 1$. We are going to show that $S$ is $d$-thick. Let $\{B_{j}\}_{j \in
\mathbb{N}}=\{B(x_{j},r_{j})\}_{j \in \mathbb{N}}$ be a covering of $Q \cap S$ such that

\begin{equation}
\label{eq6.3}
\mathcal{H}^{d}_{\infty} (Q \cap S) \geq \frac{1}{2}\sum\limits_{j \in \mathbb{N}}(r_{j})^{d}.
\end{equation}

It is clear that without loss of generality we may assume that $B_{j} \cap S \neq \emptyset$. Hence, for every $j \in \mathbb{N}$ there is a point $\widetilde{x}_{j} \in B_{j} \cap S$. It is clear that
$B_{j} \subset B(\widetilde{x}_{j},2r_{j})$ for every $j \in \mathbb{N}$. Using \eqref{eq6.3}, Ahlofrs $d$-regularity of $S$ and subadditivity of the measure $\mathcal{H}^{d}$ we obtain desirable estimate

\begin{equation}
\label{eq6.4}
\begin{split}
&\mathcal{H}^{d}_{\infty} (Q \cap S) \geq \frac{1}{2^{d+1}}\sum\limits_{j \in \mathbb{N}}(2r_{j})^{d} \geq \sum\limits_{j \in \mathbb{N}}\mathcal{H}^{d}(B(\widetilde{x}_{j},2r_{j}) \cap S)\\
&\geq \mathcal{H}^{d}(Q \cap S) \geq C r^{d}.
\end{split}
\end{equation}

3) Let $\varepsilon, \delta > 0$. Recall (see \cite{J}) that an open  set is called an $(\varepsilon,\delta)$-domain, if for any $x,y \in \Omega$ such that $\|x-y\| < \delta$ there exists a rectifiable path $\gamma_{x,y}$ of length $\frac{\|x-y\|}{\varepsilon}$ connecting $x$ and $y$ such that, for each $z \in \gamma_{x,y}$,

\begin{equation}
\label{eq6.4''}
\operatorname{dist}(z,\partial \Omega) > \varepsilon \frac{\|z-x\|\|z-y\|}{\|y-x\|}.
\end{equation}

Fix an arbitrary path-connected $(\varepsilon,\delta)$-domain $\Omega$. Consider an arbitrary cube $Q=Q(x,r)$ with $x \in \Omega$ and $r < \min\{\delta,\operatorname{diam}\Omega\}$.
Hence there is a point $y \in \Omega \setminus Q(x,r)$ such that \eqref{eq6.4''} holds.
Note that the functions $g_{x}(t)=\|x-\gamma_{x,y}(t)\|$ and $g_{y}(t)=\|y-\gamma_{x,y}(t)\|$ are continuous. Hence, using triangle inequality, it is easy to see that there exists a point $z_{0} \in \gamma_{x,y}$ such that $\min\{\|x-z_{0}\|,\|y-z_{0}\|\} \geq 2^{-1} \|x-y\|$. But this means that $\operatorname{dist}(z_{0}),\partial \Omega) \geq \frac{\|x-y\|}{4}$. Hence,

\begin{equation}
\label{eq6.5''}
Q(z_{0},\varepsilon\frac{\|x-y\|}{8}) \subset \Omega.
\end{equation}

But this mean that for every $x \in \Omega$ and every $r < \min\{\delta,\operatorname{diam}\Omega\}$

\begin{equation}
\notag
\mathcal{L}_{n}(Q(x,r) \cap \Omega) \geq C(\varepsilon)\mathcal{L}_{n}(Q(x,r)).
\end{equation}

And hence, for every $x \in \Omega$ and every $r \in (0,1)$

\begin{equation}
\notag
\mathcal{L}_{n}(Q(x,r) \cap \Omega) \geq C(\varepsilon,\delta,\operatorname{diam}\Omega)\mathcal{L}_{n}(Q(x,r)).
\end{equation}

This proves that $\Omega$ is Ahlfors $n$-regular and thus $n$-thick.

4) Let $G$ be the set, constructed in (4) of the Example 2.1. Let $Q=Q(x,r)$ be an arbitrary closed cube with $x \in G$, $r \in (0,1]$. Let $\{B_{j}\}_{j \in \mathbb{N}_{0}}=\{B(x_{j},r_{j})\}_{j \in \mathbb{N}_{0}}$ be an arbitrary sequence of balls such that $G \cap Q(x,r) \subset \cup_{j \in \mathbb{N}}B_{j}$ and

$$
\mathcal{H}^{n-1}_{\infty}(Q(x,r) \cap G) \geq \frac{1}{2} \sum\limits_{j}(r_{j})^{n-1}.
$$

Without loss of generality we may assume that $B_{j} \cap Q(x,r) \cap G \neq \emptyset$ for every $j \in \mathbb{N}$. Hence, for every $j \in \mathbb{N}$ there is a point $y_{j} \in B_{j} \cap G$. Then, it is clear that $B(y_{j},2r_{j}) \supset B_{j}$ for every $j \in \mathbb{N}$. This gives

\begin{equation}
\label{eq6.1'}
\mathcal{H}^{n-1}_{\infty}(Q(x,r) \cap G) \geq \frac{1}{2^{n}} \sum\limits_{j}(2r_{j})^{n-1}.
\end{equation}

Let $\widehat{Q}$ be the projection of the cube $Q$ to the space $\mathbb{R}^{n-1}=\{(x',x_{n}) \in \mathbb{R}^{n}|x_{n}=0\}$. For every $j \in \mathbb{N}$ let $\widehat{x}_{j}$ be the projection of the center of the ball $B_{j}$ to the space $\mathbb{R}^{n-1}$. Clearly $\widehat{Q} \cap \overline{\Omega} \subset \cup_{j \in \mathbb{N}}B(\widehat{x}_{j},2r_{j})$.

From the previous item we know that $\overline{\Omega}$ is Ahlfors $(n-1)$-regular. Hence, there exists a constant $C >0$ such that

\begin{equation}
\label{eq6.2'}
C r^{n-1} \le \mathcal{H}^{n-1}(Q(x,r) \cap G) \le \sum\limits_{j \in \mathbb{N}}\mathcal{H}^{n-1}(B(\widehat{x}_{j},2r_{j}) \cap \overline{\Omega}) \le  \sum\limits_{j \in \mathbb{N}}(2r_{j})^{n-1}.
\end{equation}

Combining \eqref{eq6.1'} and \eqref{eq6.2'} we conclude.

\textbf{Proof of Theorem \ref{Th2.2}}. We follow the classical scheme of the proof of Theorem \ref{Th2.1} (see section 1.5 of \cite{St2}). Fix $t > 0$, $f \in L_{\gamma}(\mathbb{R}^{n})$ and define

\begin{equation}
\notag
f_{1}(x):=
\begin{cases}
&f(x), \quad |f(x)| \geq \frac{t}{2\min\{1,s^{\alpha}\}},\\
&0, \quad |f(x)| < \frac{t}{2\min\{1,s^{\alpha}\}}.
\end{cases}
\end{equation}

We put $f_{2}(x):=f(x)-f_{1}(x)$ for every $x \in \mathbb{R}^{n}$.

It is clear that

\begin{equation}
\notag
|f(x)| \le |f_{1}(x)|+\frac{t}{2\min\{1,s^{\alpha}\}}.
\end{equation}

Hence, we derive

\begin{equation}
\label{eq6.7}
\operatorname{M}^{< s}[f,\alpha](x) \le \operatorname{M}^{< s}[f_{1},\alpha](x)+\frac{t}{2}.
\end{equation}

Consider the set

$$
E_{t}(f):=\{x \in \mathbb{R}^{n}|\operatorname{M}^{< s}[f,\alpha](x) > t\}.
$$

Using inequality \eqref{eq6.7} we get the inclusion

\begin{equation}
\label{eq2.21'}
E_{t}(f) \subset E_{\frac{t}{2}}(f_{1}).
\end{equation}

In order to estimate $\mathfrak{m}(E_{t}(f))$ we should use more delicate arguments than that of \cite{St2}, because our measure $\mathfrak{m}$ is not assumed to be doubling.
We are going to use Theorem \ref{Th2.6}.

It is clear that for every $x \in E_{\frac{t}{2}}(f_{1})$ there is a Euclidean ball $B_{x}=B(x,r_{x})$ with radius $r_{x} < s$ such that

$$
r^{\alpha}\fint\limits_{B(x,r)}|f_{1}(y)|\,d\mathcal{L}_{n}(y) > \frac{t}{2}.
$$

Hence we have a covering $\mathcal{F}:=\{B_{x}\}$ of the set $E_{\frac{t}{2}}(f_{1})$ with radii bounded above by $s$. Using Besicovitch's covering Theorem, we
obtain finite number of subfamilies of balls $\mathcal{G}_{1},...,\mathcal{G}_{N(n)} \subset \mathcal{F}$ such that each $\mathcal{G}_{i}$ ($i=1,...,N(n)$) is a countable collection of disjoint
balls in $\mathcal{F}$ and

$$
E_{\frac{t}{2}}(f_{1}) \subset \bigcup\limits_{i=1}^{N(n)}\bigcup\limits_{B \in \mathcal{G}_{i}}B.
$$

Note that, using \eqref{eq2.3} and the condition $d \geq n-\alpha$, for every ball $B \in \mathcal{F}$ we have the estimate

\begin{equation}
\label{eq2.18'}
\mathfrak{m}(B) \le Cr^{-\alpha}\mathcal{L}_{n}(B) \le \frac{2C}{t}\int\limits_{B}|f_{1}(y)|\,d\mathcal{L}_{n}(y).
\end{equation}

Using \eqref{eq2.21'}, \eqref{eq2.18'} and the fact that balls in every family $\mathcal{G}_{i}$ are disjoint, we have the key weak-type estimate

\begin{equation}
\begin{split}
\label{eq2.19'}
&\mathfrak{m}(E_{t}(f)) \le \mathfrak{m}(E_{\frac{t}{2}}(f_{1})) \le \frac{2C}{t}\sum\limits_{i=1}^{N(n)} \sum\limits_{B \in \mathcal{G}_{i}} \int\limits_{B}|f_{1}(y)|\,d\mathcal{L}_{n}(y)\\
& \le N(n)\frac{2C}{t}\int\limits_{|f| \geq \frac{t}{2\min\{1,s^{\alpha}\}}}
|f(y)|\,d\mathcal{L}_{n}(y).
\end{split}
\end{equation}

From \eqref{eq2.19'} using standard arguments (see the end of the proof of the corresponding theorem in section 1.5 of the book \cite{St2}) we obtain

\begin{equation}
\|\operatorname{M}^{< s}[f,\alpha]|L_{\gamma}(\mathbb{R}^{n},\mathfrak{m})\|^{\gamma}=\gamma\int\limits_{0}^{\infty}t^{\gamma-1}\mathfrak{m}(E_{t}(f)) \le C(\gamma,n,d,\alpha)\|f|L_{\gamma}(\mathbb{R}^{n},\mathcal{L}_{n})\|^{\gamma}.
\end{equation}

The theorem is proved.

\begin{Remark}
We would like to note that in the proof above we essentially used the fact that $s < +\infty$.
\end{Remark}

\textbf{Proof of Theorem \ref{Th3.1}}. Idea of proof repeat that of Theorem 5.1.12 of \cite{A}. Recall the definition of the mesh of diadic cubes $\mathcal{Q}_{n}$, $n \in \mathbb{N}_{0}$ (see the beginning of section 2). Fix a nonnegative integer $k$ and let $\mu^{k}$ be a measure such that $\mu^{k}$ has constant density and has mass equal to $2^{-kd}$ on each $Q_{k,m}$ that intersects $S$.

We now modify $\mu^{k}$ in the following way. If $\mu^{k}(Q_{k-1,m}) > 2^{-(k-1)d}$ for some $Q_{k-1,m} \in \mathcal{Q}_{k-1}$ we reduce its mass uniformly on $Q_{k-1,m}$
until it equals $2^{-(k-1)d}$. If on the other hand $\mu^{k}(Q_{k-1,m}) \le 2^{-(k-1)d}$, we leave $\mu^{k}$ unchanged on $Q_{k-1,m}$. This way we obtain a new measure $\mu^{k,1}$.
Using the fact that every cube $Q_{k-1,m}$ which has nonempty intersection with $S$ contains $\le 2^{n}$ cubes $Q_{k,m'}$ with the property $Q_{k,m'} \cap S \neq \emptyset$, we have

\begin{equation}
\notag
\mu^{k,0}(Q_{k,m}) \le \mu^{k,1}(Q_{k,m}) \le 2^{n-d}\mu^{k,0}(Q_{k,m})
\end{equation}

We repeat this procedure with $\mu^{k,1}$, obtaining $\mu^{k,2}$, and after $k$ such steps we have obtained $\mu^{k,k}$. It follows from the construction that

\begin{equation}
\label{eq6.12}
\mu^{k,k-j}(Q_{i,m}) \le h(2^{-id}).
\end{equation}

for every $j=0,1,...,k$ and every dyadic cube $Q_{i,m} \in \mathcal{Q}_{i}$ with $i=j,...,k$. Furthermore, it is clear that

\begin{equation}
\label{eq6.13}
\mu^{k,j}(Q_{k,m}) \le \mu^{k,j+1}(Q_{k,m}) \le 2^{n-d}\mu^{k,j}(Q_{k,m}), \quad j = 0,1,...,k-1.
\end{equation}

Using \eqref{eq6.12}, it is easy to see that the sequence $\{\mu^{k,k}(E)\}_{k \in \mathbb{N}_{0}}$ is bounded for every compact set $E \subset S$. Then $\{\mu^{k,k}\}_{k \in \mathbb{N}_{0}}$ has a subsequence that converges weakly to $\mu_{0}$ (see Theorem 2 in section 1.9 of \cite{Evans}), and clearly $\operatorname{supp}\mu \subset S$.

Similarly for every $j \in \mathbb{N}$ we consider the sequence $\{\mu^{k,k-j}\}_{k \geq j}$. This sequence has a subsequence
that converges weakly to $\mu_{j}$.

Fix an arbitrary $j \in \mathbb{N}$. Fix also an arbitrary Borel set $G \subset S$. Compare $\mu_{j}(G)$ and $\mu_{j-1}(G)$. Firstly note, that according to our construction for every dyadic cube $Q_{k,m}$ we have

\begin{equation}
\notag
\mu^{k,k-j+1}(Q_{k,m}) \le \mu^{k,k-j}(Q_{k,m}) \le 2^{n-d}\mu^{k,k-j+1}(Q_{k,m}), \quad k \geq j.
\end{equation}

This gives for every $f \in C_{0}(\mathbb{R}^{n})$

\begin{equation}
\label{eq6.14'}
\int\limits_{\mathbb{R}^{n}}f(x)\,d \mu^{k,k-j+1}(x) \le \int\limits_{\mathbb{R}^{n}}f(x)\,d\mu^{k,k-j}(x) \le 2^{n-d}\int\limits_{\mathbb{R}^{n}}f(x)\,d\mu^{k-j+1}(x), \quad k \geq j.
\end{equation}

Fix an arbitrary $f \in C_{0}(\mathbb{R}^{n})$. Choosing, if required, an appropriate subsequence, and passing to the limit in \eqref{eq6.14'}, we obtain

\begin{equation}
\label{eq6.14}
\int\limits_{\mathbb{R}^{n}}f(x)\,d \mu_{j-1}(x) \le \int\limits_{\mathbb{R}^{n}}f(x)\,d\mu_{j}(x) \le 2^{n-d}\int\limits_{\mathbb{R}^{n}}f(x)\,d\mu_{j-1}(x), \quad j \in \mathbb{N}.
\end{equation}

Using Borel regularity of measures $\mu_{j}$ and taking into account estimate \eqref{eq6.14}, we obtain \eqref{eq3.3}.

Now we show that $\mu_{j}(Q_{i,m}) \le 3^{n}2^{-id}$ for every $i,j \in \mathbb{N}_{0}$, $i \geq j$ and every dyadic cube $Q_{i,m} \in \mathcal{Q}_{i}$. Indeed, if $f_{i,m}$ is a continuous function with $\operatorname{supp}f_{i,m} \subset Q_{i,m}$ such that $\chi_{Q_{i,m}} \le f_{i,m} \le \chi_{3Q_{i,m}}$, then \eqref{eq6.12} gives

\begin{equation}
\notag
\mu_{j}(Q_{i,m}) \le \int\limits_{\mathbb{R}^{n}}f_{i,m}(x)\,d\mu_{j}(x) = \lim\limits_{l \to \infty}\int\limits_{\mathbb{R}^{n}}f_{i,m}(x)\,d\mu^{k_{l},k_{l}-j}(x) \le 3^{n}2^{-id}.
\end{equation}

Hence, using the fact that every closed ball $B(x,r)$ with $x \in \mathbb{R}^{n}$, $r \in (0,2^{-k}]$ has nonempty intersection with $\le 5^{n}$ dyadic cubes $Q_{k(r),m}$ (where $k(r)$ is chosen such that $r \in [2^{-i(r)},2^{-i(r)+1})$), we obtain \eqref{eq3.1}.

It remains to prove \eqref{eq3.2}. Fix an arbitrary nonempty index set $\mathcal{A} \subset \mathbb{Z}^{n}$ and $k \in \mathbb{N}_{0}$. Consider the set $V=\cup_{m \in \mathcal{A}}Q_{k,m}$. Fix an arbitrary $l \in \mathbb{N}$, $l \geq k$ and note that every $x \in V \cap S$ belongs to some (or several) $Q^{(j)} \in \mathcal{Q}_{n_{j}}$, $k \le n_{j} \le l$ such that
$\mu^{l,l-k}(Q_{j,m}) = 2^{-jd}$. We obtain a disjoint covering, $S \cap V \subset \cup_{j}Q^{(j)}$, such that

\begin{equation}
\notag
\mu^{l,l-k}(V \cap S) = \sum\limits_{j}\mu^{l,l-k}(Q^{(j)})=\sum\limits_{j}2^{-n_{j}d} \geq \inf \sum\limits_{i}2^{-n_{i}d},
\end{equation}

where the infimum is taken over all finite or denumerable coverings of $V \cap S$ with $Q^{(i)} \in \cup_{l \geq k}\mathcal{Q}_{l}$. The right hand side is independent
of $l$, and letting $l \to \infty$ it follows that

\begin{equation}
\mu_{k}(V \cap S) \geq \inf \sum\limits_{i}2^{-n_{i}d}.
\end{equation}

To finish the proof it remains to replace cubes by balls. Suppose that $S \cap V \subset \bigcup\limits_{j=1}^{\infty}Q^{(j)}$, where $Q^{(j)} \in \mathcal{Q}_{n_{j}}$ and $n_{j} \geq k$ for every $j \in \mathbb{N}$. Then there is a constant $A_{n}$ such that each $Q^{(j)}$ is contained in the union of $\le A_{n}$ balls with radius $2^{-n_{j}}$. Thus $\mathcal{H}^{d}_{\infty}(V \cap S) \le A_{n} \inf \sum\limits_{j}2^{-n_{j}d}$, where the infimum is taken as above, and thus

\begin{equation}
\notag
\mathcal{H}^{d}_{\infty}(V \cap S) \le A_{n} \mu_{k}(V \cap S).
\end{equation}

The proof of Theorem \ref{Th3.1} is completed.

\textbf{Acknowledgement.}
The authors wish to acknowledge their indebtedness to
Professor P.~Shvartsman,
who read an~early version of the manuscript and made many valuable comments.
Based on those,
we substantially reworked the paper.

\end{document}